\documentclass[11pt]{article}
\usepackage[utf8]{inputenc}

\usepackage[left=3cm, right=3cm, top=3cm]{geometry}
\usepackage{numprint}
\usepackage{float}
\usepackage{graphics}
\usepackage{tikz}
\usepackage{pgfplots} 
\pgfplotsset{compat=1.18} 
\usetikzlibrary{positioning, shapes.geometric, arrows.meta}
\usepackage{xcolor}
\usepackage{amsmath,mathtools,amsthm}
\usepackage{amssymb}
\usepackage{bbm} 
\usepackage{amsfonts}
\usepackage{dsfont}
\usepackage{mathrsfs}

\usepackage{url}
\usepackage{subfig}
\usepackage{enumitem}
\usepackage{appendix}

\usepackage{multicol,multirow}
\usepackage{color}

\usepackage{relsize}
\usepackage[gen]{eurosym}
\definecolor{DarkGreen}{rgb}{0.0, 0.5, 0.0}
\definecolor{mygreen}{RGB}{28,172,0} 
\definecolor{mylilas}{RGB}{170,55,241}

\usepackage{datetime}
\newdateformat{monthyeardate}{\monthname[\THEMONTH], \THEYEAR}
\usepackage{authblk} 

\usepackage{hyperref}
\hypersetup{
    colorlinks,
    citecolor=blue,
    filecolor=blue,
    linkcolor=blue,
    urlcolor=blue
}
\usepackage{booktabs}
\usepackage{siunitx}
\usepackage{actuarialsymbol}
\usepackage{todonotes}

\usepackage{scalerel} 

\newtheorem{theorem}{Theorem}[section]
\newtheorem{lemma}[theorem]{Lemma}
\newtheorem{corollary}[theorem]{Corollary}
\newtheorem{proposition}[theorem]{Proposition}

\theoremstyle{definition}

\newtheorem{remark}[theorem]{Remark}

\newtheorem{example}[theorem]{Example}

\newcommand{\linftyF}{{ \ell^\infty(\Fc) }}

\newcommand{\FotimesG}{{ \Fc \otimes \Gc }}
\newcommand{\linftyFG}{{ \ell^\infty(\FotimesG) }}

\newcommand{\BLoneR}{{ \mathrm{BL}_1(\RR) }}
\newcommand{\BLoneE}{{ \mathrm{BL}_1(E) }}

\newcommand{\BLD}{ {\mathrm{BL}_1(D)}}
\newcommand{\BLDE}{ {\mathrm{BL}_1(D\times E)}}
\newcommand{\BLDRd}{ {\mathrm{BL}_1(D\times \RR^d)}}

\newcommand{\simiid}{{ \, \overset{\mathrm{iid}}{\sim} \, }}
\newcommand{\eqdistr}{{ \, \overset{d}{=} \, }}

\newcommand{\Prob}{{ \operatorname{Prob} }}
\renewcommand{\outer}{\text{o}}
\newcommand{\inner}{\text{o}}

\newcommand{\empProcessGn}{ {\mathbb{G}_n} }
\newcommand{\empProcessGnBoot}{ {\mathbb{G}_n^*} }

\newcommand{\Expec}{ {\mathbb{E} }}
\newcommand{\Var}{\mathop{\rm Var}\nolimits}
\newcommand{\Cov}{\mathop{\rm Cov}\nolimits}

\renewcommand{\to}{{ \, \rightarrow \, }}
\newcommand{\cvweakly}{{ \, \rightsquigarrow \, }}

\newcommand{\norm}[1]{\left\lVert#1\right\rVert}

\newcommand{\Tnst}{{ T_n^* }}
\newcommand{\Tneq}{{ T_n^{*,\mathrm{eq}} }}
\newcommand{\Tncent}{{ T_n^{*,\mathrm{c}} }}

\newcommand{\thetahat}{{ \hat{\theta} }}
\newcommand{\thetahatn}{{ \hat{\theta}_n }}

\newcommand{\thetastMD}{{ \theta_n^{*,MD} }}
\newcommand{\thetast}{{ \theta_n^{*} }}

\newcommand{\hathn}{{ \hat{h}_n }}

\newcommand{\Htheta}{{ H_{\theta} }}

\newcommand{\HHn}{{ \HH_n }}
\newcommand{\HHnst}{{ \HH_n^* }}
\newcommand{\FFn}{{ \FF_{n} }}
\newcommand{\GGn}{{ \GG_{n} }}

\newcommand{\HHnUV}{{ \HH_{n,U,V} }}
\newcommand{\FFnU}{{ \FF_{n,U} }}
\newcommand{\GGnV}{{ \GG_{n,V} }}

\def\argmin{\mathop{\rm argmin}}
\def\argmax{\mathop{\rm argmax}}
\def\lin{\mathop{\rm lin}\nolimits}

\newcommand{\Ac}{{ \mathcal{A} }}
\newcommand{\Bc}{{ \mathcal{B} }}
\newcommand{\Cc}{{ \mathcal{C} }}

\newcommand{\Fc}{{ \mathcal{F} }}
\newcommand{\Gc}{{ \mathcal{G} }}
\newcommand{\Hc}{{ \mathcal{H} }}

\def\Nc{\mbox{$\mathcal N$}}

\def\Xc{{\mathcal X}}
\def\Yc{{\mathcal Y}}

\newcommand{\Zc}{{ \mathcal{Z} }}

\def\FF{{\mathbb F}}
\def\GG{{\mathbb G}}
\def\HH{{\mathbb H}}

\def\MM{{\mathbb M}}
\def\NN{{\mathbb N}}
\def\PP{{\mathbb P}}
\def\QQ{{\mathbb Q}}
\def\RR{{\mathbb R}}

\title{\textbf{Bootstrapping not under the null?}}

\author[1]{Alexis Derumigny}
\author[2]{Miltiadis Galanis}
\author[1]{Wieger Schipper\footnote{This research is (partly) financed by the NWO Spinoza prize awarded to A.W. van der Vaart by the Netherlands Organisation for Scientific Research (NWO).}}
\author[1]{Aad van der Vaart}

\affil[1]{Delft Institute of Applied Mathematics (DIAM), 

Delft University of Technology, Mekelweg 4, Delft, 2628 CD, The Netherlands.}
\affil[ ]{
\href{mailto:A.F.F.Derumigny@tudelft.nl}{a.f.f.derumigny@tudelft.nl}; 
\href{mailto:W.R.Schipper@tudelft.nl}{w.r.schipper@tudelft.nl};
\href{mailto:A.W.vanderVaart@tudelft.nl}{a.w.vandervaart@tudelft.nl}
\vspace{0.2cm}
}

\affil[2]{Department of Informatics and Telecommunications,

National and Kapodistrian University of Athens,
Zografou 161 22, Athens, Greece}
\affil[ ]{\href{mailto:mgalan@di.uoa.gr}{mgalan@di.uoa.gr}}

\date{\today}

\begin{document}

\maketitle

\begin{abstract}
We propose a bootstrap testing framework for a general class of hypothesis tests, which allows resampling under the null hypothesis as well as other forms of bootstrapping. We identify combinations of resampling schemes and bootstrap statistics for which the resulting tests are asymptotically exact and consistent against fixed alternatives. We show that in these cases the limiting local power functions are the same for the different resampling schemes. We also show that certain naive bootstrap schemes do not work.
To demonstrate its versatility, we apply the framework to several examples: independence tests, tests on the coefficients in linear regression models, goodness-of-fit tests for general parametric models and for semi-parametric copula models. Simulation results confirm the asymptotic results and suggest that in smaller samples non-traditional bootstrap schemes may have advantages. This bootstrap-based hypothesis testing framework is implemented in the R package BootstrapTests.

\medskip

\noindent
\textbf{Keywords:} Conditional weak convergence, parametric bootstrap, empirical bootstrap, nonparametric bootstrap, null bootstrap, hypothesis test, independence test, goodness-of-fit test.

\smallskip

\noindent
\textbf{MSC2020:} Primary 62F40, 62G10, 62G09, 62F03;
secondary 62E20.
\end{abstract}


\section{Introduction}
\label{SectionIntroduction}

Many problems of hypothesis testing on a probability distribution $H$ can be written in the form
\begin{align}
\label{eq_hypothesis_problem}
    \Hc_0: \phi(H) = 0 \quad \text{ versus }\quad
    \Hc_1: \phi(H) \neq 0,
\end{align}
for a functional $\phi$ with values in some normed space $(E, \| \cdot \|_E)$, see e.g. \cite{bickel2001bootstrap}.
Examples include independence tests and goodness-of-fit tests. 
Given an i.i.d.\ sample $X_1, \dots, X_n$ following $H$, a natural test statistic is
\begin{align}
\label{EqNaturalTestStatistic}
    T_n := \sqrt{n} \| \phi(\HHn) \|_E,
\end{align}
where $\HHn$ is the empirical distribution of $X_1, \dots, X_n$, 
the discrete measure that puts mass $1/n$ at each of the observations.
Often, the asymptotic distribution of $T_n$ is difficult to compute, and bootstrap methods are used to estimate its limiting distribution and determine critical values for the test. These consist of creating 
a sample of bootstrap observations $X_1^*, \dots, X_n^*$ according to a measure $R_n = R_n(X_1, \dots, X_n)$ 
that depends
on the original observations, and use the conditional distribution given $X_1,\ldots, X_n$ of an appropriate
bootstrap counterpart $T_n^*$ of $T_n$  as an estimate of the distribution of $T_n$. 

Bootstrap samples can be created in many ways and so can the bootstrap test statistics $T_n^*$. 
A main interest in the present paper is the question whether ``to bootstrap under the null hypothesis or not'', i.e.\ choosing 
the measures $R_n$ to belong to the null hypothesis or not.
The former is often recommended in the literature. However, we shall see that bootstrapping under a general measure, for instance
the empirical measure $\HHn$, may also work provided the bootstrap statistics $T_n^*$
are defined properly. 

Let $\HHnst$ denote the empirical distribution of the bootstrap values $X_1^*, \dots, X_n^*$, which are assumed to be i.i.d.\ according to $R_n$, for given $X_1,\ldots, X_n$.
Two possible bootstrap statistics are the \textit{equivalent bootstrap test statistic}  and
the \textit{centred bootstrap test statistic}, given by
\begin{align}
\label{eq_T_n_eq}
    \Tneq &:= \sqrt{n} \| \phi(\HHnst) \|_E,\\
\label{eq_T_n_c}
    \Tncent &:= \sqrt{n} \| \phi(\HHnst) - \phi(\HHn) \|_E.
\end{align}
In the literature, it is sometimes recommended to use $\Tncent$ (see e.g. \cite{hall1991two,bickel2001bootstrap}) while some others recommend $\Tneq$, especially for goodness-of-fit testing  (\cite{beran1988prepivoting,stute1993bootstrap,genest2009goodness}).
We shall show that both statistics may work, but in different situations, depending on the bootstrap scheme $R_n$.
This has already been remarked, see the discussions in \cite{bickel2001bootstrap} and \cite{derumigny2017tests}.
The intuition is that under the null hypothesis the original test statistic can be rewritten as $\sqrt{n} \| \phi(\HHn) - \phi(H)\|_E$, since $\phi(H)=0$ under
$\Hc_0$. Since $\HHn$ is the empirical distribution of a sample from $H$, and $\HHnst$ is the empirical
distribution of a sample from $R_n$, the natural bootstrap counterpart of this statistic is
\begin{align}
\label{EqCorrectBootstrap}
    \Tnst := \sqrt{n} \| \phi(\HHnst) - \phi(R_n) \|_E.
\end{align}
This reduces to \eqref{eq_T_n_eq} in the case of a bootstrap under the null hypothesis (that is, $\phi(R_n)=0$),
and to \eqref{eq_T_n_c} in the case of the empirical bootstrap (i.e.\ $R_n=\HH_n$). We shall show that
\eqref{EqCorrectBootstrap} is typically correct, which explains that both \eqref{eq_T_n_eq} 
and \eqref{eq_T_n_c} can be correct, depending on $R_n$. 

Here, correctness refers first to the level of the test: using the quantiles of the
bootstrap distribution of $\Tnst$ as cutoff for the test based on $T_n$ will yield level $\alpha$,
at least asymptotically as $n\rightarrow\infty$. One might think that the power of the test is
more sensitive to the type of bootstrap scheme. However, we shall show that relative to the usual
ways in which (asymptotic) power is evaluated, there is no difference between the various bootstrap
schemes, if based on the statistic \eqref{EqCorrectBootstrap}.

The statistic \eqref{EqCorrectBootstrap} allows for other bootstrap schemes $R_n$ in addition to \eqref{eq_T_n_eq} and \eqref{eq_T_n_c}, such as
a bootstrap from a parametric or semiparametric model. In this paper, we develop a general theory that includes various types of bootstraps as special cases. We apply this theory to 
the examples of independence testing, goodness-of-fit testing and testing a regression coefficient,
exhibiting several correct bootstrap schemes in every example. A general assumption is that
the function $\phi$ is Hadamard differentiable. To also include some other examples, we 
extend the theory to test statistics of the more general form
$$T_n=\sqrt n\|\phi(\HHn,\hat\theta_n)\|_E,$$ for estimators $\hat \theta_n=\hat\theta_n(X_1,\ldots, X_n)$,
where the bootstrap test statistic $T_n^*$ is formed
using an appropriate bootstrap estimator $\theta_n^*$.

The article is organised as follows. Section~\ref{sec:general_framework} describes the framework and the main results in a general formulation. In Section~\ref{sec:extension_framework_parameters} this is extended to include parameter estimators,
with special attention for parametric bootstrap schemes.
Applications are developed in Section~\ref{SectionIndependenceTesting} (testing for independence),
Section~\ref{sec:GoF_Setting} (goodness-of-fit testing), Section~\ref{SectionLinearRegression} 
(testing a slope) and Section~\ref{SectionCopulas} (testing goodness-of-fit of copula models).
A simulation study is performed in Section~\ref{SectionSimulation} to investigate the finite-sample performance of the tests in various settings. We provide an \texttt{R} package, called BootstrapTests \cite{BootstrapTestPackage}, which implements our bootstrap-based hypothesis testing procedures. Technical results are collected in an appendix. 

\medskip

\textbf{Notation}: Denote by $P \otimes Q$ the product distribution of two distributions $P$ and $Q$.
For a measure $P$, we write $P^\outer$ and $P_\inner$ respectively for the inner and the outer measures of $P$.
The shorthand $Qf$, for a measurable function $f$ and (signed) measure $Q$, is notation for $Qf := \int f\, \mathrm{d}Q$. For a set $\Fc$, the space $\linftyF$ is the set of all bounded functions $z: \Fc\to \mathbb{R}$ equipped with the uniform norm
$\|z\|_\Fc = \sup_{f \in \Fc} |z(f)|$. A signed measure $Q$ can be identified with the
 the map $f\mapsto Qf$, which is contained in  $\linftyF$ provided the map is bounded.

\section{General framework}
\label{sec:general_framework}

For $n \in \NN$, let  $X_1, \dots, X_n$ be  i.i.d.\ random variables
taking values in a measurable space $(\Xc, \Ac)$, following the distribution $H$.
We embed the probability distributions $H$ into some normed space $(D,\|\cdot\|_D)$ and for a
given map $\phi: D\to E$ with values in another normed space $(E, \| \cdot \|_E)$,
we consider the testing problem \eqref{eq_hypothesis_problem}.

The empirical distribution $\HHn=n^{-1}\sum_{i=1}^n\delta_{X_i}$ 
of the observations gives another element of $D$, and
is used to form the test statistic $T_n$ given by \eqref{EqNaturalTestStatistic}.
The null hypothesis is rejected for large values of $T_n$. 

\begin{example}[Normed space, $\linftyF$]
In many examples the spaces $D$ or $E$ are spaces of uniformly bounded functions.
For instance, measures on Euclidean space may be identified with their cumulative distribution
functions and viewed as elements of Skorohod space, or the set of uniformly
bounded functions $z: \RR\to \RR$. 

A fairly general setup is to identify
a probability measure $H$ with the map $f\mapsto Hf:=\int f\,dH$ from a
a given set $\Fc$ of measurable functions $f: \Xc\to\RR$, in which case
the empirical measure is identified with the map
$f\mapsto \HHn f:= n^{-1}\sum_{i=1}^n f(X_i)$. If $H$ and $\Fc$ are such that 
$\sup_{f\in\Fc}|Hf|<\infty$, then $H$ and $\HH_n$ are elements of the space $\linftyF$
of uniformly bounded functions $z: \Fc\to \RR$, which can be equipped with the
uniform norm $\|z\|_\Fc=\sup_{f\in\Fc}|z(f)|$.
\end{example}

We employ a bootstrap scheme to set a critical value for the test.
Given a probability measure $R_n=R_n(X_1,\ldots, X_n)$ and for
given observations $X_1,\ldots, X_n$, we
draw an i.i.d.\ sample $X_1^*, \dots, X_n^*$ from $R_n$, and form their
empirical measure $\HHnst$. Next, for given $\alpha\in(0,1)$ and still for
given $X_1,\ldots, X_n$, we determine 
the $(1-\alpha)$-quantile $\xi^*_{n, 1-\alpha}=\xi^*_{n, 1-\alpha}(X_1,\ldots, X_n)$ 
of the (conditional) distribution of the bootstrap statistic \eqref{EqCorrectBootstrap}.
The null hypothesis is then rejected if $T_n\ge \xi^*_{n, 1-\alpha}$. In practice,
these quantiles are computed as the empirical quantiles of a large number of simulated
bootstrap values $T_n^*$.

The asymptotic level and power of these test are the limits
as $n\rightarrow\infty$ of the probabilities $\Pr_H\bigl(T_n\ge \xi^*_{n, 1-\alpha}\bigr)$
computed under distributions $H$ that belong to the null or alternative hypotheses, respectively. We wish
to investigate these in their dependence on the choice of bootstrap scheme $R_n$ and bootstrap statistic.
Some possibilities are given by the following examples.

\begin{example}[Empirical bootstrap]
The choice $R_n(X_1,\dots,X_n)=\HHn$ is known as the empirical bootstrap. It is the original
choice of the bootstrap (see \cite{Efron1979}), and corresponds to redrawing the bootstrap
values $X_1^*,\ldots, X_n^*$ with replacement from the original observations. As it does not
refer to the testing problem, this type of bootstrap is not often recommended for testing.
We shall see that it can actually work well in combination with the bootstrap statistic \eqref{eq_T_n_c},
but will typically fail with the statistic \eqref{eq_T_n_eq}.
\end{example}

\begin{example}[Independent bootstrap]
For two-dimensional observations $(X_i, Y_i)_{i=1}^n$, the measure $R_n$ can be set equal to
the product measure of the empirical measures on $X$ and $Y$ separately, i.e.\
$R_n = \PP_n^X \otimes \PP_n^Y$. If the null hypothesis asserts that  
$X_i$ and $Y_i$ are independent, then this gives an example of bootstrapping under
the null hypothesis. We study this further in Section~\ref{SectionIndependenceTesting}
along with the empirical bootstrap.
\end{example}

\begin{example}[Parametric bootstrap]
For a given parametrised family of distributions $(H_\theta: \theta\in\Theta)$ and
an estimator $\hat\theta_n=\hat\theta_n(X_1,\ldots, X_n)$, we can set $R_n=H_{\hat\theta_n}$. 
We consider this further in Sections~\ref{sec:GoF_Setting} and~\ref{SectionCopulas}.
\end{example}

In the remainder of this section we derive general results, which cover these and other examples.
In this generality, a proper formulation requires a precise description of the bootstrap
scheme, including a specification of the joint measurability structure underlying the 
variables $X_1, \ldots, X_n, X_1^*,\ldots, X_n^*$ and measures $R_n$. In order not to burden
the message, we defer a description of the most technical details to Appendix~\ref{AppendixMeasurability}.

\subsection{Central limit theorems and their bootstrap counterparts}

The results are based on convergence in distribution of the
processes $\sqrt n(\HH_n-H)$ and $\sqrt n(\HHnst-R_n)$ in the metric space $D$. For the first we
assume ordinary convergence in distribution
\begin{equation}
\label{EqConvergenceHHn}
\sqrt n(\HH_n-H)\cvweakly \GG_H,\qquad \text{in } D.
\end{equation}
This means that $\Expec^\outer h\bigl(\sqrt n(\HH_n-H)\bigr)\to \Expec h(\GG_H)$, for every
continuous, bounded function $h: D\to\RR$. For greater flexibility the maps $\sqrt n(\HH_n-H)$
need not be Borel measurable, but the definition is understood in the sense of 
Hoffmann-J\o rgensen, employing outer expectation $\Expec^\outer$. 
See e.g.\ \cite{Billingsley} for weak convergence theory, or 
Part~1 of \cite{van2023weak}, for this theory extended to non-measurable elements.
The limit $\GG_H$ is always assumed to be a tight Borel measurable map in $D$.

For the bootstrap process $\empProcessGnBoot := \sqrt{n}(\HHnst - R_n)$ we assume convergence in distribution,
conditionally given $X_1,\ldots, X_n$. Usually the bootstrap scheme stabilizes in the limit in the sense that the sequence of measures $R_n = R_n(X_1, \ldots, X_n)$
converges to a deterministic limit $R(H)$. For instance, the empirical bootstrap scheme $R_n = \HHn$ converges to $R(H) = H$.
In all our examples the sequence $\empProcessGnBoot$ converges conditionally
in distribution in $D$ to the variable $\GG_{R(H)}$.
A convenient notation for this convergence is 
\begin{equation}
\label{eq_conditional_bootstrap}
    \empProcessGnBoot
    = \sqrt{n}(\HHnst - R_n)
    \, | \, X_1, \dots, X_n
    \cvweakly
    \GG_{R(H)},\qquad \text{in } D,  \quad \text{a.s. or in probability.} 
\end{equation}
To give a rigorous meaning to this type of convergence, which also takes care of measurability, we use a version of the bounded Lipschitz metric. It is known that the convergence \eqref{EqConvergenceHHn} is equivalent to the convergence
to zero $\sup_{h\in \BLD}\bigl|\Expec^\outer h\bigl(\sqrt n(\HH_n-H)\bigr)- \Expec h(\GG_H)\bigr|\to0$, where $\BLD$ is the
set of 1-Lipschitz functions $h: D\to [-1, 1]$. In agreement with this, the preceding display is formally understood to mean
\begin{equation}
    \sup_{h \in \BLD} \Big|
    \Expec^\outer[h(\empProcessGnBoot) | X_1, \dots, X_n]    - \Expec h(\GG_{R(H)}) \Big|
    \to 0, \qquad
    H^\infty\text{-outer a.s. or prob}.
\label{eq:cond_cv_Gnboot_precise}
\end{equation}
Here $H^\infty$ denotes the joint distribution of $(X_n)_{n \geq 1}$,
and outer a.s.\ or in probability convergence means that the left side of the equation is bounded above 
by measurable random variables that converge a.s.\ or in probability to zero (see Chapter~I.9 in \cite{van2023weak}). 
To allow for the possibility that the process $\GG_n^*$ is not
measurable in the bootstrap variables $X_1^*,\ldots,X_n^*$, we use the
outer (conditional) expectation  $\Expec^\outer[h(\empProcessGnBoot) | X_1, \dots, X_n]$ rather
than the ordinary conditional expectation. See Appendix~\ref{AppendixMeasurability} for a precise
definition of the bootstrap scheme and this expected value.

In the following examples we specialise these definitions to the case that $D=\linftyF$ is a space of bounded functions (defined under ``Notation'' in Section~\ref{SectionIntroduction}),
and note that for the important case of the empirical bootstrap $R_n=\HHn$,
 the bootstrap convergence \eqref{eq_conditional_bootstrap} (or its precise version \eqref{eq:cond_cv_Gnboot_precise})
is implied by the convergence \eqref{EqConvergenceHHn} of the ordinary empirical process and hence automatic.

\begin{example}[Donsker class, $\linftyF$]
For $D=\ell^\infty(\Fc)$ and the empirical measure $\HHn$ identified with the map
$f\mapsto \HH_nf=n^{-1}\sum_{i=1}^n f(X_i)$, the process $\sqrt n (\HH_n-H)$ in \eqref{EqConvergenceHHn}
is the empirical process $(\empProcessGn f)_{f \in \Fc}$, given by
\begin{equation}
\label{eq_empprocfunc}
    f \mapsto \empProcessGn f
    := \sqrt{n}(\HHn - H) f
    = \frac{1}{\sqrt{n}} \sum_{i=1}^n \bigl(f(X_i) - Hf\bigr).
\end{equation}
The set of functions $\Fc$ is called an $H$-\emph{Donsker class} if and only if the
convergence \eqref{EqConvergenceHHn} is valid. 
The limit process $\GG_H = (\GG_H f)_{f \in \Fc}$ is an $H$-Brownian bridge. In view of the multivariate
central theorem, it can be seen to be a zero-mean Gaussian process with covariance function
$\Expec[\GG f_1 \GG f_2]= H(f_1 f_2) - Hf_1 Hf_2$ for $f_1, f_2 \in \Fc.$

There is a considerable literature on empirical processes, giving many examples of Donsker classes
(see \cite{DudleyBook}, \cite{van1996}, \cite{van2023weak} and references).
For use in a nonparametric testing setup, classes of functions that are $H$-Donsker for every probability measure $H$, called \emph{universal Donsker classes},
are most attractive. These include the classical Donsker class of indicator functions of cells in $\RR^d$, and more generally
all bounded, suitably measurable, Vapnik-Chervonenkis classes.
\end{example}

\begin{example}[Empirical bootstrap, $\linftyF$]
\label{ExampleEmpiricalBootstrapConvergence}
For $D=\linftyF$ and the empirical bootstrap $R_n = \HHn$, the conditional
convergence \eqref{eq_conditional_bootstrap} is satisfied
for every Donsker class $\Fc$ with square integrable envelope function,
with $R(H) = H$ (See \cite{GineZinn1990,Gine1997} or \cite{van2023weak}, Section~3.7.1).
Furthermore, the sequence  $R_n = \HHn$ tends in $\linftyF$ to $H$, outer almost surely
(every Donsker class is a Glivenko-Cantelli class).

This means that for $D=\linftyF$ and the empirical bootstrap,
the convergence assumptions \eqref{EqConvergenceHHn} and \eqref{eq_conditional_bootstrap} of Theorem~\ref{thm:asymptotic_behavior_T}, below,
reduce to the single assumption that $\Fc$ is $H$-Donsker.
\end{example}

\begin{example}[Empirical bootstrap, Banach space]
\label{ExampleEmpiricalBootstrapConvergenceBanachSpace}
Probability measures $H$ on $\RR^p$ can be identified with their cumulative distribution functions
$x\mapsto F(x)= H(-\infty,x]$, which in turn can be viewed as elements of $L_2(\RR^p,\mu)$, for a finite
Borel measure $\mu$ on $\RR^d$. In particular, the empirical measure $\HH_n$ can be identified with
the empirical distribution function $x\mapsto \FF_n(x)=n^{-1}\sum_{i=1}^n 1_{X_i\le x}$, which is the
average of the i.i.d.\ random elements $1_{X_1\le \cdot},\ldots, 1_{X_n\le \cdot}$ in $L_2(\RR^p,\mu)$.
By the central limit theorem in $L_2(\RR^p,\mu)$, the sequence  $\sqrt n(\FF_n-F)$ converges in distribution in 
$L_2(\RR^p,\mu)$.

The empirical bootstrap yields the random elements  $1_{X_1^*\le \cdot},\ldots, 1_{X_n^*\le \cdot}$ in $L_2(\RR^p,\mu)$,
with average the bootstrap  empirical distribution function $x \mapsto \FF_n^*(x)=n^{-1}\sum_{i=1}^n 1_{X_i^*\le x}$. 
By the bootstrap central limit theorem in $L_2(\RR^p,\mu)$,  the sequence of processes $\sqrt n(\FF_n^*-\FF_n)$ converges
in distribution to the same limit as the original processes, conditionally given almost every sequence $X_1,X_2,\ldots$.

Thus for $D=L_2(\RR^p,\mu)$ and the empirical bootstrap,
the convergence assumptions \eqref{EqConvergenceHHn} and \eqref{eq_conditional_bootstrap} of Theorem~\ref{thm:asymptotic_behavior_T}, below,
are satisfied automatically. 

The restriction to a finite measure $\mu$ ensures that every cumulative distribution function $F$ is an element
of  $L_2(\RR^p,\mu)$. By restricting to $F$ with sufficiently light left and right tails and considering differences
$F-F_0$ for a suitable fixed cumulative distribution $F_0$ (to control the right tail), this finding can be extended
to more general measures $\mu$, including Lebesgue measure. 

Actually, the present example is a special case of Example~\ref{ExampleEmpiricalBootstrapConvergence}, as
the central limit theorem in a separable Banach space, such as $L_2(\RR^p,\mu)$, holds if and only if the unit ball
of the dual space is a Donsker class. (See e.g.\ \cite{van2023weak}, Section~2.1.4. This is based on the identification between an element $x$ of a separable Banach space $D$ and the mapping $x^{**} : f \mapsto f(x)$, which satisfies that $x^{**} \in \linftyF$ for $\Fc$ the unit ball of the dual space of $D$.)
By the same identification, we can extend the bootstrap central limit theorems proved for $\linftyF$-valued random elements to $D$-valued random elements.
This is true for any identification of the measures $\delta_x$ with elements in any separable Banach space such that the central limit theorem holds for $\delta_{X_1},\ldots,\delta_{X_n}$ in this
Banach space.
For a Hilbert space the central limit theorem holds provided the second moment of the norm
is finite, which is easily the case for the bounded variables $1_{X_i\le\cdot}$ in $L_2(\RR^p,\mu)$ for a finite
measure $\mu$, and also true 
for the random elements $1_{X_i\le\cdot}-F_0$ for a general measure $\mu$ whenever $\int \bigl[(F-F_0)^2+ F(1-F)\bigr]\,d\mu <\infty$.
\end{example}

\subsection{Asymptotic results for general bootstrap-based testing: level, power and local power}

We derive the limiting distributions of $T_n$ and $\Tnst$ under the assumption that the map $\phi: D \to E$ is Hadamard differentiable. 
A map $\phi: D_{\phi} \subset D \to E$ from a subset of a normed space $D$ into
a normed space $E$ is said to be \emph{uniformly Hadamard-differentiable} tangentially to a set $D_0 \subset D$ 
at $H \in D_\phi$ if there exists  a continuous, linear map $\phi_H' : \lin D_0 \to E$ such that
     $t_n^{-1}\bigl(\phi(H_n + t_n h_n) - \phi(H_n)\bigr)\rightarrow \phi_H'(h)$,
for all converging sequences $t_n \to 0$ in $\RR$, and $H_n \to H$ in $D$ and $h_n \to h\in D_0$ 
such that $H_n\in  D_\phi$ and $H_n + t_n h_n \in D_\phi$ for every $n$.
If the convergence is verified only for fixed $H_n=H$ independent of $n$, then the ``uniformly'' is dropped and the
map is said to be \emph{Hadamard-differentiable} tangentially to the set $D_0$. The phrase 
``tangentially to $D_0$'' is omitted if $D_0=D$ (see \cite{van2023weak}, Chapter~3.10, or \cite{Gill1989}). 

\begin{theorem}
\label{thm:asymptotic_behavior_T}
Suppose that $\phi: D \to E$ is Hadamard-differentiable at $H$ and at $R(H)$ tangentially to  a measurable linear subspace $D_0\subset D$.
If  \eqref{EqConvergenceHHn} and \eqref{eq_conditional_bootstrap} hold  in outer probability and the sequence $\sqrt n\bigl(R_n-R(H)\bigr)$ is asymptotically tight in $D$, where $\GG_H$ and $\GG_{R(H)}$ take their values in $D_0$, then,  under $H$,
\begin{align}
 T_n &\cvweakly \|\phi_{H}'(\GG_H)\|_E,\qquad\text{if }\phi(H)=0,\label{EqConvergenceTn}\\
     T_n &\cvweakly \infty, \ \qquad\qquad\qquad\text{if }\phi(H)\neq 0,\label{EqTntoinfty}\\
\noalign{\noindent while}
       \label{eq:cond_cv_Tnstar}
        \Tnst | X_1, \dots, X_n &\cvweakly
        \|\phi_{R(H)}' (\GG_{R(H)}) \|_E, \qquad \text{in outer prob.}
    \end{align}
If $\phi$ is uniformly Hadamard differentiable at $R(H)$, then the
condition of asymptotic tightness of the sequence $\sqrt n\bigl(R_n-R(H)\bigr)$ can be relaxed to the condition that $R_n\to R(H)$ in outer probability. If, moreover, \eqref{eq_conditional_bootstrap} holds outer almost surely and $R_n\to R(H)$ outer almost surely, then \eqref{eq:cond_cv_Tnstar} is valid also outer almost surely.
\end{theorem}

\begin{proof}
Assumption \eqref{EqConvergenceHHn}  and the Delta-method applied to $\phi$ followed by the continuous mapping theorem applied to the function  $z\mapsto \|z\|_E$
(Theorems~3.10.4 and~1.3.6 in \cite{van2023weak}) give
that $S_n:=\sqrt n\|\phi(\HHn)-\phi(H)\|_E\cvweakly 
\|\phi_H'(\GG_H)\|_E$. Assertion \eqref{EqConvergenceTn} is the special case that $\phi(H)=0$,
since $T_n=S_n$ in that case. If $\phi(H)\neq0$, then
$\sqrt n\|\phi(H)\|_E\rightarrow\infty$.
Since the sequence of variables $S_n$ is bounded in probability (by Lemma 1.3.8(ii) in \cite{van2023weak}), and 
$T_n = \sqrt{n}\|\phi(H_n)\|_E   \geq \sqrt{n} \|\phi(H)\|_E - S_n$
by the reverse triangle inequality, we see that
$\Prob_\inner(T_n > t)
\geq \Prob_\inner(S_n<\sqrt{n} \|\phi(H)\|_E - t)
\rightarrow 1$, for every $t$, which is equivalent to \eqref{EqTntoinfty}.

The proof of assertion \eqref{eq:cond_cv_Tnstar} follows the same lines as the proof of the first assertion, 
but must take proper care of the measurability issues involved in conditioning.
We apply the conditional Delta-method (see Lemma~\ref{lemma:Functional_Delta_Method_boot} or Theorem~3.10.13 in \cite{van2023weak}) to the convergence \eqref{eq_conditional_bootstrap} (or \eqref{eq:cond_cv_Gnboot_precise}) in combination with tightness of the
sequence $\sqrt n\bigl(R_n-R(H)\bigr)$, and next the conditional continuous mapping theorem (see
Lemma~\ref{lemma:Continuous_mapping_theorem_boot}).
\end{proof}

\medskip

We now show asymptotic properties of the bootstrap-based hypothesis test itself.
Remember that the bootstrap test at level $\alpha$ rejects $\Hc_0$ if
$T_n > \xi^*_{n, 1-\alpha}$ 
where
$\xi^*_{n, \alpha} = \xi^*_{n, \alpha}(X_1, \dots, X_n)$
is the $\alpha$-quantile of the conditional distribution of $T^*_n$, given $X_1, \dots, X_n$.%
\footnote{We shall silently assume that the variables $T_n$ and $T_n^*$ and the conditional distribution of the latter
can be defined in a measurable way, so that these quantiles are well defined. Otherwise, we might
theoretically use measurable majorants, along the lines of Problem~1.10.1 in \cite{van2023weak}, but this would
be hard to implement in practice.}

\begin{corollary}\label{cor_power}
Under the conditions of Theorem~\ref{thm:asymptotic_behavior_T}, and for (i) under
the further condition that the distribution of $\|G_H\|_E$ does not have an atom at
its $(1-\alpha)$-quantile,
    \begin{enumerate}[label=(\roman*)]
        \item If $\phi(H)=0$ and
        $\| \phi_{R(H)}'(\GG_{R(H)}) \|_E
        \eqdistr \| \phi'_{H}(\GG_{H}) \|_E$,
        then $P^\outer(T_n \geq \xi^*_{n, 1-\alpha})\to \alpha$.
        \item If $\phi(H)\neq 0,$
        then $P_\inner(T_n \ge \xi^*_{n, 1-\alpha}) \to1$.
    \end{enumerate} 
More precisely, (i) assumes \eqref{EqConvergenceHHn} and \eqref{eq_conditional_bootstrap} and $R_n\to R(H)$ under $H\in\Hc_0$,
whereas (ii) assumes \eqref{eq_conditional_bootstrap} and \eqref{EqTntoinfty} under $H\in \Hc_1$.
\end{corollary}

\begin{proof}
(i).  By \eqref{eq:cond_cv_Tnstar}, the sequence $\xi^*_{n,1- \alpha}$ tends
in probability to the $(1-\alpha)$-quantile of the variable $\|\phi' _{R(H)}(\GG_{R(H)})\|_E$.   
Under the assumption in (i), this is equal to  the $(1-\alpha)$-quantile of the variable $\|\phi'_{H}(\GG_{H})\|_E$.
Combined with \eqref{EqConvergenceTn}, this gives $P^\outer(T_n\geq \xi^*_{1-\alpha})\rightarrow \alpha$.
(See Problem~1.10.1 in \cite{van2023weak} to handle possible non-measurable maps $T_n$.)

(ii). In view of \eqref{EqTntoinfty} it suffices to show that the sequence $\xi^*_{n,1- \alpha}$  is bounded
in probability. This is true under \eqref{eq:cond_cv_Tnstar} even if the variable $\|\phi' _{R(H)}(\GG_{R(H)})\|_E$ has an atom at 
the $(1-\alpha)$-quantile of the variable  $\|\phi' _{R(H)}(\GG_{R(H)})\|_E$.
\end{proof}

The choice of bootstrap scheme $R_n$ will determine the measure $R(H)$ that appears in the bootstrap limit in \eqref{eq_conditional_bootstrap}. For (asymptotically) correct type 1 error probabilities, (i) of the preceding lemma imposes the condition
\begin{equation}
    \label{EqBootstrapLimit}
\| \phi_{R(H)}'(\GG_{R(H)}) \|_E \eqdistr \| \phi'_{H}(\GG_{H}) \|_E.
 \end{equation}
This is trivially satisfied if $R(H)=H$, which will be seen to be true for $H$ in the null hypothesis for most of our examples.

The second assertion of the corollary shows that the bootstrap test is consistent under every  fixed alternative, in the sense that the power converges to 1. This is true without further conditions on the bootstrap scheme. One might hope that differences between various bootstrap schemes would become apparent by the performance under sequences of alternatives $H_n$ approaching the null hypothesis, but this is not the case. 

For the usual \emph{contiguous alternatives} (see \cite{LeCam1960,LeCamYang} or \cite{van2000asymptotic}, Chapters~6 and~14), this can be seen without further calculations. By the definition of contiguity a sequence of random variables $\xi_n(X_1,\ldots, X_n)$ converges
in probability to a fixed value under $H$ if and only if it converges in probability under
a contiguous sequence $H_n$, to the same value. In particular, the sequence of bootstrap quantiles $\xi^*_{n, 1-\alpha}=\xi^*_{n, 1-\alpha}(X_1,\ldots, X_n)$ converges in probability to the same value under $H$ and under any contiguous alternatives $H_n$. By \eqref{EqConvergenceTn} the limit is the $(1-\alpha)$-quantile of the variable $\|\phi'_{R(H)}(\GG_{R(H)})\|_E$. This does depend on the bootstrap scheme $R_n$ through $R(H)$, but this dependence disappears if we choose
the bootstrap scheme to satisfy \eqref{EqBootstrapLimit}, which is necessary to obtain correct type 1 error. Under \eqref{EqBootstrapLimit} the bootstrap quantiles tend in probability to the same value under all contiguous alternatives, and the power of the test only depends on the differential behaviour of $T_n$ under $H$ and $H_n$, which is independent of the choice of bootstrap scheme.

The following corollary derives an explicit expression for the power
under local alternatives $H_n$ such that $\sqrt n\,\phi(H_n)\rightarrow\tau$, 
for some $\tau\in E$.
Formally, we now have a triangular array of observations,
where for each $n\in\NN$ the variables $X_{1,n}, \dots, X_{n,n}$ are an i.i.d.\ sample following a distribution $H_n$.
To keep track of all samples coming from different distributions $H_n$, set $S_1:=(X_{1,1}), S_2:= (X_{1,2},X_{2,2}),\dots, S_n:= (X_{1,n},\dots,X_{n,n})$.

\begin{corollary}
\label{cor_power_alternative_hypotheses}
Suppose that $\phi: D \to E$ is  Hadamard-differentiable at $H$ and at $R(H)$ tangentially to a  measurable linear subspace $D_0\subset D$.
Assume that $\sqrt n(\HHn-H_n)\cvweakly \GG_H$ under $H_n$ and that $\sqrt{n}(\HHnst - R_n) \, | \, S_1, \dots, S_n 
    \cvweakly \GG_{R(H)}$, in outer probability  in $D$, where $\GG_H$ and $\GG_{R(H)}$ take their values in $D_0$.
If $\sqrt{n}\,\phi(H_n) \rightarrow \tau$,  for some $\tau \in E$ and the sequence $\sqrt n\bigl(R_n-R(H)\bigr)$ is asymptotically tight in $D$, then for every $\alpha$
such that the $(1-\alpha)$-quantile $\xi_{1-\alpha}(R(H))$ of $\|\phi'_{R(H)}(\GG_{R(H)})\|_E$ is a continuity point of both
$\|\phi'_{R(H)}(\GG_{R(H)})\|_E$ and $\|\phi'_{H}(\GG_{H}) + \tau\|_E$,
\begin{align}
P_{H_n}(T_n > \xi^*_{n,1-\alpha})
        \rightarrow P \big( \|\phi'_{H}(\GG_{H}) + \tau \|_E
        > \xi_{1-\alpha}(R(H)) \big).
       \label{eq:limiting_local_power}
\end{align}
\end{corollary}

\begin{proof}
The assumption of conditional convergence in distribution of $\sqrt{n}(\HHnst - R_n)$ replaces assumption 
\eqref{eq:cond_cv_Gnboot_precise} in Theorem~\ref{thm:asymptotic_behavior_T}, and leads in the same way 
by application of the Delta-method and continuous mapping theorem to
convergence in probability of the bootstrap quantiles $\xi^*_{n,1-\alpha} $ to the $(1-\alpha)$-quantiles 
$\xi_{1-\alpha}(R(H))$ of the variable $\|\phi'_{R(H)}(\GG_{R(H)})\|_E$, whenever this is a continuity point of the latter variable.

Similarly, the assumption combined with the Delta-method give that the sequence of variables $G_n=\sqrt{n}(\phi(\HHn)-\phi(H_n))$ tends in distribution in $\linftyF$ to $\phi_H'(\GG_H)$. By the decomposition $T_n = \sqrt{n}\|\phi(\HHn)\|=\|G_n+\sqrt{n}\phi(H_n)\|$, Slutsky's lemma and the continuous mapping theorem, we see that $T_n \cvweakly \|\phi'_H(\GG_H)+\tau\|_E$.

By Slutsky's lemma we have weak convergence of $T_n -\xi^*_{n,1-\alpha}$ to $\|\phi'_{H}(\GG_{H}) + \tau\|_E- \xi_{1-\alpha}(R(H))$. 
This implies convergence of the cumulative distribution function at $0$ provided this is a continuity point.
\end{proof}

\begin{remark}
The cumulative distribution function of the norm of a Borel measurable Gaussian variable in a separable Banach space (such as
$\|\phi'_{H}(\GG_{H}) +\tau\|_E$ or $\|\phi'_{R(H)}(\GG_{R(H)}) \|_E$) is continuous except possibly at the left end point of its support
(and strictly increasing). Moreover, if the variable
is centered, then the cumulative distribution function is continuous everywhere. Hence with rare exceptions, quantiles are continuity points.
\end{remark}


\subsection{Bootstraps that do not work}
\label{proof:rem:unconsistent_combinations_indep}
The bootstrap statistic \eqref{EqCorrectBootstrap} reduces to the equivalent bootstrap statistic \eqref{eq_T_n_eq} for a bootstrap 
$R_n$ under the null hypothesis (i.e. $\phi(R_n)=0$) and to the centred bootstrap statistic \eqref{eq_T_n_c} in case of the empirical bootstrap $R_n=\HHn$. Two other possibilities would be:
\begin{enumerate}[label=(\roman*)]
\item the equivalent statistic $\Tnst=\sqrt n\|\phi(\HHnst)\|_E$ combined with the empirical  bootstrap $R_n=\HHn$,
\item the centred statistic $\Tnst=\sqrt n\|\phi(\HHnst)-\phi(\HHn)\|_E$ combined with a bootstrap $R_n$ under the null hypothesis.
\end{enumerate}
Neither of these possibilities leads to a correct type 1 error, in general, but much worse, the powers of the resulting tests
tend to zero at any alternative.

\begin{theorem}
Assume that the conditions of Theorem~\ref{thm:asymptotic_behavior_T} hold, where the variables
$\GG_H$ and $\GG_{R(H)}$ are tight centered Gaussian random variables.
Assume that the distribution of $\|\GG_H\|_E$ does not have any atoms and that $\alpha<1/2$.
In case (ii) assume, moreover, that \eqref{EqBootstrapLimit} holds for $H$ with $\phi(H)=0$. 
Then 
$\limsup_{n\to\infty} P^\outer(T_n \geq \xi^*_{n, 1-\alpha})=0$, for every $H$ with $\phi(H)$ contained in the support of $\GG_H$ for (i) or 
$\phi(H)$ contained in the support of $\GG_{R(H)}$ for (ii).
\end{theorem}

\begin{proof}
(i). We can write $T_n=\sqrt n\|\phi(\HH_n)\|_E=\|W_n-\tau_n\|_E$, for
$W_n=\sqrt n\bigl(\phi(\HH_n)-\phi(H)\bigr)$ and $\tau_n=-\sqrt n\,\phi(H)$, where $W_n\cvweakly W=\phi_H'(\GG_H)$.
Similarly, we can write $T_n^*=\|W_n^*+W_n-\tau_n\|_E$, for $W_n^*=\sqrt n\bigl(\phi(\HHnst)-\phi(\HH_n)\bigr)$.
We have $W_n^* | X_1, \dots, X_n \cvweakly W^*$,
where $W^*\sim W$, and hence $W_n^* | W_n \cvweakly W^*$, since $W_n$ is a function of $X_1, \dots, X_n$.
For given $\tau_n$, define maps $\xi_n: E\to \RR$ by letting
$\xi_n(w)$ be the smallest solution to the equation $P^\outer\bigl(\|\bar W_n+w-\tau_n\|_E\le \xi_n(w) \bigr)\ge 1-\alpha$,
for $\bar W_n\sim W_n^*| W_n=w$ and every $w$.
Then the bootstrap critical values  are $\xi_{n,1-\alpha}^*=\xi_n(W_n)$ and the probability in the assertion is
$P^\outer\bigl(T_n \geq \xi^*_{n, 1-\alpha})=P^\outer(\|W_n-\tau_n\|_E \geq \xi_n(W_n)\bigr)$.

For every $e^*$ in the dual space $E^*$ of norm $\|e^*\|_{E^*}=1$, the variable 
$\|\bar W_n+w-\tau_n\|_E$ is  bounded below by the variable $e^*(\tau_n-w)-e^*(\bar W_n)$.
Combined with the definition of $\xi_n(w)$ this implies that $\inf_{\|e^*\|_{E^*}=1}P^\outer\bigl(-e^*(W_n^*)\le \xi_n(W_n)-e^*(\tau_n-W_n) | W_n \bigr)\ge 1-\alpha$. 

For every $\xi\in\RR$ and $\delta>0$, there exists a Lipschitz function $h_{\xi,\delta}: \RR\to [0,1]$ with 
$1_{(-\infty,\xi]}\le h_{\xi,\delta}\le 1_{(-\infty,\xi+\delta]}$, where the Lipschitz constants can be independent of $\xi$ (but depend on $\delta$). Then the functions
$h_{\xi,\delta}\circ e^*: E\to [0,1]$ are Lipschitz, with Lipschitz constants uniformly bounded in $\xi$ and $\|e^*\|_{E^*}\le1$. Therefore 
the convergence $W_n^*|W_n\cvweakly W^*$ in outer probability implies that
\begin{align*}
    &\sup_{\|e^*\|_{E^*}=1}\Bigl|\Expec^\outer \Bigl(h_{\xi_n(W_n)-e^*(\tau_n-W_n),\delta}\bigl(-e^*(W_n^*)\bigr)| W_n\Bigr)
-\Expec_{\bar W} h_{\xi_n(W_n)-e^*(\tau_n-W_n),\delta}\bigl(-e^*(\bar W)\bigr)\Bigr|\\
&\qquad\lesssim \sup_{h\in \BLoneE}\Bigl|\Expec^\outer \bigl(h(W_n^*)| W_n\bigr)-\Expec h(\bar W)\Bigr|\to 0,
\end{align*}
in outer probability. Here $\bar W$ denotes a random element with $\bar W\sim W$ that is independent of $W_n$ and the
expectation $\Expec_{\bar W}$ is evaluated with respect to $\bar W$ for fixed $W_n$.

Combining the preceding paragraphs and using that $1_{(-\infty,\xi]}\le h_{\xi,\delta}$, we conclude that
$\inf_{\|e^*\|_{E^*}=1}\Expec_{\bar W} h_{\xi_n(W_n)-e^*(\tau_n-W_n),\delta}\bigl(-e^*(\bar W)\bigr)\ge 1-\alpha +o_P(1)$. Next using that \\
$h_{\xi,\delta}\le 1_{(-\infty,\xi+\delta]}$, we conclude that 
$\inf_{\|e^*\|_{E^*}=1}\Expec_{\bar W} 1\{-e^*(\bar W)\le \xi_n(W_n)-e^*(\tau_n-W_n)+\delta\}\ge 1-\alpha +o_P(1)$.
Because the variables $-e^*(\bar W)$ are
one-dimensional normal with mean zero and variances 
$\sigma_{e^*}^2:=\Expec e^*(W)^2$, this is equivalent to 
$$\inf_{\|e^*\|_{E^*}=1}\Phi\Bigl(\frac{\xi_n(W_n)-e^*(\tau_n-W_n)+\delta}{\sigma_{e^*}}\Bigr)\ge 1-\alpha +o_P(1).$$
We conclude that $\xi_n(W_n)\ge \underline\xi(\tau_n-W_n)-\delta+o_P(1)$, for $$\underline\xi(w):=\sup_{\|e^*\|_{E^*}=1} \bigl(e^*(w)+ \sigma_{e^*}\Phi^{-1}(1-\alpha)\bigr).$$
In view of the last line of the first paragraph,
the proof can be completed by showing that 
$P^\outer(\|W_n-\tau_n\|_E \geq \underline\xi(\tau_n-W_n)-\delta\bigr)\to 0$, as $n\to\infty$ followed by $\delta\to0$. 
Because $w\mapsto \|w\|_E-\underline \xi(w)$ is uniformly continuous, we can replace $W_n$ by $W$  and it suffices to show that 
\begin{equation}
\label{EqPowerNotWorking}
\lim_{\delta\to0}\limsup_{n\to\infty} 
P\bigl(\|W-\tau_n\|_E \geq \underline\xi(\tau_n-W)-\delta\bigr)=0.
\end{equation}
We prove this separately in the cases that $\tau=\tau_n/\sqrt n=-\phi(H)$ is 
zero or non-zero. 


Assume that $\tau=0$. By the Hahn-Banach theorem there exists $e_w^*\in E^*$ with 
$\|e_w^*\|_{E^*}=1$ such that $e_w^*(-w)=\|w\|_E$, for every $w\in E$. Thus
$\underline\xi(w)\ge \|w\|_E+\sigma_{e_w^*}\Phi^{-1}(1-\alpha)$ and it suffices
to show that $P(0\ge \sigma_{e_W^*}\Phi^{-1}(1-\alpha)- \delta)\to 0$ as $\delta\to0$, i.e.
$\sigma_{e_W^*}>0$ almost surely. Now for any $e^*\in E^*$, the 
equality $\sigma_{e^*}=0$ implies that $e^*(W)=0$ almost surely. By continuity of
$e^*$ this gives that $e^*$ vanishes on the support $S$ of $W$ (the smallest closed set in $E$ with $P(W\in S)=1$). 
If $W: \Omega\to E$ is defined on $\Omega$, define $\Omega_0=\{\omega\in\Omega: \sigma_{e^*_{W(\omega)}}=0\}$ and $\Omega_1=\{\omega\in\Omega: W(\omega)\in S\}$. Then  $\Omega_0=\{\omega\in\Omega: e_{W(\omega)}^*(s)=0,\forall s\in S\}$ by the preceding remark and
hence $e^*_{W(\omega)}(W(\omega))=0$ for all  $\omega\in\Omega_0\cap \Omega_1$.
Because by construction $e^*_{W(\omega)}(W(\omega))=\|W(\omega)\|_E$ for every $\omega\in\Omega$ and
$\|W\|_E>0$ almost surely by assumption, it follows that $P(\Omega_0)=P(\Omega_0\cap \Omega_1)=0$. This concludes the proof of \eqref{EqPowerNotWorking} in the case
that $\tau=0$.

Assume that $\tau\not=0$  and that $\tau$ is in the support of $W$, as assumed in the theorem. By the Hahn-Banach theorem, for every $n$ and $w$ 
there exists $e_{n,w}^*\in E^*$ with  $\|e_{n,w}^*\|_{E^*}=1$ 
such that $e_{n,w}^*(\tau_n-w)=\|\tau_n-w\|_E$. Then
$e_{n,w}^*(\tau)=e_{n,w}^*(\tau_n-w)/\sqrt n+e_{n,w}^*(w)/\sqrt n
=\|\tau_n-w\|_E/\sqrt n+e_{n,w}^*(w)/\sqrt n\to\|\tau\|_E$, uniformly in $w$ such that $\|w\|_E$ is bounded. Because $\tau$ is in the support of $W$,
there exists $g$ in the reproducing kernel Hilbert space of $W$ such that $\|g-\tau\|_E<\|\tau\|_E/2$ (see e.g.\ \cite{vdVRKHS}) and hence $|e_{n,w}^*(g)-e_{n,w}^*(\tau)|\le \|g-\tau\|_E< \|\tau\|_E/2$, so that $\liminf_{n\to\infty} e_{n,w}^*(g)\ge \|\tau\|_E/2$, uniformly in $\|w\|_E\le K$, for any given $K$.
Given an orthonormal basis $(h_i)$ of its reproducing Hilbert space, the variable
$W$ and $g$ can be represented as $W=\sum_{i=1}^\infty Z_ih_i$ and $g=\sum_{i=1}^\infty g_ih_i$, for $(Z_i)$ a sequence of i.i.d.\ standard normal variables and $(g_i)\in\ell_2$. Then $\sigma_{e^*}^2=\sum_i e^*(h_i)^2$ and hence $e^*(g)=\sum_i g_ie^*(h_i)\le \|g\|_{\ell_2}\sigma_{e^*}$, by the Cauchy-Schwarz inequality, for any $e^*\in E^*$.
It follows that $\liminf_{n\to\infty}\inf_{\|w\|_E\le K}\sigma_{e_{n,w}^*}^2\ge \|\tau\|_E/(2\|g\|_{\ell_2})$, for every $K$.
We conclude that there exists $c>0$ so that $\underline\xi(\tau_n-w)\ge \|\tau_n-w\|_E+ c \Phi^{-1}(1-\alpha)$, for every $w$ with $\|w\|_E\le K$, for sufficiently large $n$. Then the probability in \eqref{EqPowerNotWorking} is bounded
above by $P(\|W-\tau_n\|_E \geq \|W-\tau_n\|_E +c\Phi^{-1}(1-\alpha)-\delta,\|W\|_E\le K)
+P(\|W\|_E>K)$. The second probability on the right can be made arbitrarily small, while
the first probability vanishes for sufficiently small $\delta$ for given $K$.
This finishes the proof of \eqref{EqPowerNotWorking}.


(ii). We have $T_n=\|W_n-\tau_n\|_E$, for $W_n$ and $\tau_n$ as before, and since $\phi(R_n)=0$, we have
$T_n^*=\|W_n^*-W_n+\tau_n\|_E$, for $W_n^*=\sqrt n\bigl(\phi(\HHnst)-\phi(R_n)\bigr)$ with
 $W_n^*| W_n\cvweakly W^*$, for $W^*=\phi_{R(H)}'(\GG_{R(H)})$.
Apart from the signs in  $\|W_n^*-W_n+\tau_n\|_E$, this is the same as before, where 
$\|W^*\|_E\sim \|W\|_E$ if \eqref{EqBootstrapLimit} holds.
\end{proof}

\section{General framework with parameter estimators}\label{sec:extension_framework_parameters}
While statistics of the type \eqref{EqNaturalTestStatistic} cover many interesting examples, in some situations more natural statistics take the form
$$T_n=\sqrt n \|\phi(\HHn,\hat\theta_n)\|_E.$$
Here $\phi: D\times \RR^d\to E$ is a differentiable map,
$\HHn$ is the empirical measure of $X_1,\ldots, X_n$ as before,
 and $\hat\theta_n=\hat\theta_n(X_1,\ldots, X_n)$ are given statistics with values in $\RR^d$. In this case
the analogue of the bootstrap statistic \eqref{EqCorrectBootstrap} is given by 
\begin{align}
\label{EqCorrectBootstrapTheta}
    \Tnst := \sqrt{n} \| \phi(\HHnst,\theta_n^*) - \phi(R_n,\hat\theta_n) \|_E,
\end{align}
where $\HHnst$ is the  empirical measure of the bootstrap sample $X_1^*,\ldots, X_n^*$, taken from $R_n$, as before, and 
$\theta_n^*=\theta_n^*(X_1^*,\ldots, X_n^*; X_1,\ldots, X_n)$  are appropriate maps.
The results of the preceding section readily extend to statistics of this type. The proof of the following theorem is similar to the proofs in the preceding section and is omitted.

Instead of \eqref{EqConvergenceHHn} and \eqref{eq_conditional_bootstrap}, assume that there exist vectors $\theta(H)\in\RR^d$
such that, for certain limit variables $U_H$ and $V_{R(H)}$,
\begin{align}
\label{EqJointConvergenceEPtheta}
    \sqrt n\bigl(\HHn-H,\hat\theta_n-\theta(H)\bigr)&\cvweakly (\GG_H,U_H),\ \qquad\quad\text{in } D\times\RR^d,\\
\label{EqJointConvergenceBootstrapEPtheta}
    \sqrt n\bigl(\HHnst-R_n,\theta_n^*-\hat\theta_n\bigr) | X_1,\ldots, X_n&\cvweakly (\GG_{R(H)},V_{R(H)}),\quad\text{in } D\times\RR^d,\quad
\text{in outer prob.}
\end{align}
The condition for consistency under the null hypothesis now becomes
\begin{equation}
\label{EqBootstrapLimitWithTheta}
    \|\phi_{R(H)}' (\GG_{R(H)},V_{R(H)}) \|_E\eqdistr \|\phi_{H}'(\GG_H,U_H)\|_E.
\end{equation}

\begin{theorem}
\label{TheoremBootstrapExtension}
Assume that the map $\phi: D \times \RR^d \to E$ is
Hadamard differentiable at $\bigl(H,\theta(H)\bigr)$ and at $\bigl(R(H),\theta(R(H)\bigr)$, both tangentially to the same measurable linear space $D_0\times\RR^d\subset D \times \RR^d$.
If  \eqref{EqJointConvergenceEPtheta} and
\eqref{EqJointConvergenceBootstrapEPtheta} hold, where $\GG_H$ and $\GG_{R(H)}$ take their values in $D_0$, and the sequence $\sqrt n\bigl(R_n-R(H),\hat\theta_n-\theta(H)\bigr)$ is asymptotically tight in $D\times \RR^d$, then, under $H$,
\begin{align}
 T_n &\cvweakly \|\phi_{H}'(\GG_H,U_H)\|_E,\qquad\text{if }\phi\bigl(H,\theta(H)\bigr)=0,\\
T_n &\cvweakly \infty,\qquad\qquad\qquad\qquad\text{if }\phi\bigl(H,\theta(H)\bigr)\neq 0,\label{EqPowerTestWithParameter}\\
\noalign{\noindent while}\Tnst | X_1, \dots, X_n &\cvweakly
\|\phi_{R(H)}' (\GG_{R(H)},V_{R(H)}) \|_E, \quad \text{in outer probability}.
\end{align}
Consequently, for every $\alpha$ so that the variable $\|\phi_{H}'(\GG_H,U_H)\|_E$  does not have an atom at its $(1-\alpha)$-quantile,
\begin{enumerate}[label=(\roman*)]
 \item If $\phi\bigl(H,\theta(H)\bigr)=0$ and \eqref{EqBootstrapLimitWithTheta} holds,
        then $P(T_n \geq \xi^*_{n, 1-\alpha})\rightarrow \alpha$.
 \item If $\phi\bigl(H,\theta(H)\bigr)\neq 0$,  then $P_\inner(T_n > \xi^*_{n, 1-\alpha}) \rightarrow 1$.
\end{enumerate} 
If $\phi$ is uniformly Hadamard differentiable at $\bigl(R(H),\theta(R(H)\bigr)$, then the asymptotic tightness of 
$\sqrt n\bigl(R_n-R(H),\hat\theta_n-\theta(H)\bigr)$ can be relaxed to the convergence
$(R_n,\hat\theta_n)\to \bigl(R(H),\theta(R(H))\bigr)$ in outer probability in $D\times\RR^d$.
\end{theorem}

\begin{proof}
This follows the same lines as the proofs of Theorem~\ref{thm:asymptotic_behavior_T} and Corollary~\ref{cor_power}.
\end{proof}

\begin{remark}
    The conditions needed for assertions (i) and (ii) in Theorem~\ref{TheoremBootstrapExtension} can be more precisely stated as follows. Assertion (i) is true under conditions \eqref{EqJointConvergenceEPtheta} and \eqref{EqJointConvergenceBootstrapEPtheta} and the convergence $(R_n,\hat\theta_n)\to \bigl(R(H),\theta(R(H))\bigr)$ under $H\in \Hc_0$, whereas (ii) needs only \eqref{EqJointConvergenceBootstrapEPtheta} and  \eqref{EqPowerTestWithParameter} for $H\in\Hc_1$.
\end{remark}

\begin{remark}
If the estimators $\hat\theta_n$ are of the special form
$\hat\theta_n=\theta(\HHn)$ for a Hadamard differentiable map
$\theta: D\to \RR^d$, then natural bootstrap versions are $\theta_n^*=\theta(\HHnst)$. In this case
Theorem~\ref{TheoremBootstrapExtension} is a corollary of
Theorem~\ref{thm:asymptotic_behavior_T}, applied with the Hadamard differentiable map $H\mapsto \phi\bigl(H, \theta(H)\bigr)$.
\end{remark}

\subsection{Parametric bootstrap}
In several examples the null hypothesis consists of a parametrised set $\{H_\theta: \theta\in\Theta\}$ of distributions, for
$\Theta$ an open subset of $\RR^d$. In that case a \emph{parametric bootstrap} $R_n=H_{\hat\theta_n}$, for given
estimators $\hat\theta_n=\hat\theta_n(X_1,\ldots, X_n)$ is natural, together with bootstrap values $\theta_n^*=\hat\theta_n(X_1^*,\ldots, X_n^*)$
calculated from the bootstrap sample $X_1^*,\ldots, X_n^*$ from $R_n$ in the same way as $\hat\theta_n$ is calculated from the
original observations.
The following theorem shows that a slight strengthening of \eqref{EqJointConvergenceEPtheta} (at null distributions) then implies 
\eqref{EqJointConvergenceBootstrapEPtheta}. Suppose that for every converging sequence $h_n\to h$ in $\RR^d$ and
for $X_1,\ldots,X_n$ a sample from $H_{\theta_0+h_n/\sqrt n}$,
\begin{equation}
\label{EqJointConvergenceEPthetaRegular}
    \sqrt n\bigl(\HHn-H_{\theta_0+h_n/\sqrt n},\hat\theta_n-\theta_0-h_n/\sqrt n\bigr)\cvweakly (\GG_{H_{\theta_0}},U_{H_{\theta_0}}),\quad\text{in } D\times\RR^d,
\text{ under } H_{\theta_0+h_n/\sqrt n}.
\end{equation}
The convergence of the second coordinates $\sqrt n(\hat\theta_n-\theta_0-h_n/\sqrt n)$ under $H_{\theta_0+h_n/\sqrt n}$ with the
same limit distribution for every sequence $h_n\to h$, is known as the \emph{regularity} of the estimator sequence
$\hat\theta_n$. Regularity can be seen as a form of local robustness, in that small, vanishing perturbations of the underlying
distribution do not change the limit behaviour.
Many estimator sequences are regular  (see \cite{Hajek70} or \cite{van2000asymptotic}), including the empirical
measure $\HH_n$ as an estimator of $H$ (see \cite{van2023weak}, Section~3.12.1).
Regularity may be ascertained directly, or can be derived using Le Cam's lemma
given local asymptotic normality of the parametric model. The latter approach is taken in \cite{GenestRemillard2008}.

Regularity can be combined with estimators that are known to converge at $\sqrt n$-rate. The following stronger assumption can handle
estimators that are just consistent. Suppose that for every sequence $\theta_{0,n}\to\theta_0$,
\begin{equation}
\label{EqJointConvergenceEPthetaRegularNoRate}
    \sqrt n\bigl(\HHn-H_{\theta_{0,n}},\hat\theta_n-\theta_{0,n}\bigr)\cvweakly (\GG_{H_{\theta_0}},U_{H_{\theta_0}}),\quad\text{in } D\times\RR^d,
\text{ under } H_{\theta_{0,n}}.
\end{equation}

\begin{theorem}[Parametric bootstrap]
\label{TheoremParametricBootstrap}
Let $R_n=H_{\hat\theta_n}$ and $\theta_n^*=\hat\theta_n(X_1^*,\ldots, X_n^*)$. If 
the sequence $\sqrt n\bigl(\hat\theta_n-\theta(H)\bigr)$ is tight under $H$ and \eqref{EqJointConvergenceEPthetaRegular} holds, 
for $\theta_0:=\theta(H)$,
then \eqref{EqJointConvergenceBootstrapEPtheta} holds under $H$ with $R(H)=H_{\theta_0}$ and $V_{R(H)}=U_{H_{\theta_0}}$. 
If \eqref{EqJointConvergenceEPthetaRegular} can be strengthened to \eqref{EqJointConvergenceEPthetaRegularNoRate}, then
the tightness of the sequence $\sqrt n\bigl(\hat\theta_n-\theta(H)\bigr)$ can be relaxed to the consistency $\hat\theta_n\to \theta(H)$. Finally, if $\theta \mapsto H_\theta$ is continuous in $D$ at $\theta(H)$ and $\hat\theta_n\to \theta(H)$, then $R_n\to R(H)$ in $D$ in outer probability under $H$, and if $\theta \mapsto H_\theta$ is Hadamard differentiable at $\theta(H)$ and the sequence $\sqrt n\bigl(\hat\theta_n-\theta(H)\bigr)$ is tight, then the sequence $\sqrt n\bigl(R_n-R(H)\bigr)$ is asymptotically tight.
\end{theorem}

\begin{proof}
For $\theta\in\Theta$ consider the bounded Lipschitz distance given by
$$L_n(\theta)=\sup_{h\in\BLDRd}\Bigl| \Expec^\outer h\bigl(\sqrt n(\HH_n-H_{\theta},\hat\theta_n-\theta)\bigr)
-\Expec h(\GG_{H_{\theta_0}},U_{H_{\theta_0}})\Bigr|.$$
For $R_n=H_{\hat\theta_n}$, the bootstrap empirical measure $\HH_n^*$ behaves conditionally given $X_1,\ldots, X_n$ as the
ordinary empirical measure of a sample from $H_\theta$, for $\theta=\hat\theta_n$. Therefore, 
the claim of the theorem is equivalent to the convergence $L_n(\hat\theta_n)\to0$ in outer probability, under $H_{\theta_0}$.
Since $R_n=H_{\hat\theta_n}$ and $\theta_n^*=\hat\theta_n(X_1^*,\ldots, X_n^*)$,
assumption \eqref{EqJointConvergenceEPthetaRegular} says that $g_n(h_n):=L_n(\theta_0+h_n/\sqrt n)\to0$, for every sequence $h_n\to h$.
If the sequence $\hat h_n:=\sqrt n(\hat\theta_n-\theta_0)$ is asymptotically tight,
the extended continuous mapping theorem, as given in Theorem~1.11.1 in \cite{van2023weak}, shows that
$L_n(\hat\theta_n)=g_n(\hat h_n)\cvweakly 0$.

If \eqref{EqJointConvergenceEPthetaRegularNoRate} holds, then we follow the same line of argument but directly applied to the functions $\theta\mapsto L_n(\theta)$ without the localisation to the functions $g_n$.

That $R_n\to H_{\theta_0}$ in probability follows  similarly from the assumption that $H_\theta\to H_{\theta_0}$. That the sequence $\sqrt n\bigl(R_n-R(H)\bigr)$ is asymptotically tight follows by the Delta-method.
\end{proof}

\begin{remark}
For $h_n=0$ condition \eqref{EqJointConvergenceEPthetaRegular} reduces to \eqref{EqJointConvergenceEPtheta} at $H=H_{\theta_0}$,
which is in the null hypothesis. For the power of the test, Theorem~\ref{TheoremBootstrapExtension} uses 
condition \eqref{EqJointConvergenceBootstrapEPtheta} for $H$ in the alternative hypothesis, which is unrelated to 
\eqref{EqJointConvergenceEPthetaRegular}.
Alternatively, the power can often be easily analysed by a direct verification of  \eqref{EqPowerTestWithParameter}.
\end{remark}

\begin{example}[Donsker class, $\linftyF$]
\label{ExampleAsymptoticLinearity}
Many estimator sequences $\hat\theta_n$ are \emph{asymptotically linear} in the sense that 
\begin{equation}
       \sqrt{n}\bigl(\thetahat_n - \theta(H)\bigr) = \frac{1}{\sqrt{n}}\sum_{i=1}^n\psi_H (X_i) + o_{P_H}(1),
       \label{EqAsymptoticLinearity}
\end{equation}
for measurable functions $\psi_H: \Xc\to\RR^d$ with $H\psi_H = 0$ and $H \| \psi_H \|^2 < \infty$ (called influence functions). In that case the joint convergence \eqref{EqJointConvergenceEPtheta} with $D=\linftyF$ holds if and only if $\Fc$ is an $H$-Donsker class, with the limit variable $(\GG_H, U_H)=(\GG_H, \GG_H\psi_H)$, where
on the right $\GG_H$ is an $H$-Brownian bridge indexed by the class
of functions $\Fc\cup \{\psi_H\}$. 

This solves many examples, although below we also consider cases where the estimators $\hat\theta_n$ are not asymptotically Gaussian, and the joint convergence is obtained otherwise.

Similar remarks can be made for the convergence \eqref{EqJointConvergenceBootstrapEPtheta} of the bootstrap process. If $\theta_n^*$ 
is constructed similarly to $\hat\theta_n$, then it is reasonable to expect that in analogy to \eqref{EqAsymptoticLinearity},
\begin{equation}
       \sqrt{n}\bigl(\theta_n^* - \hat\theta_n\bigr) = \frac{1}{\sqrt{n}}\sum_{i=1}^n\psi_{R(H)} (X_i^*) + \hat\epsilon_n^*,
       \label{EqAsymptoticLinearityBootstrap}
\end{equation}
where $\hat\epsilon^*| X_1,\ldots, X_n\cvweakly 0$ in probability. Together with asymptotic linearity of the bootstrap process $\sqrt n(\HHnst-R_n)$, this will readily give the joint convergence \eqref{EqJointConvergenceBootstrapEPtheta}.

Alternatively, for the parametric bootstrap \eqref{EqJointConvergenceBootstrapEPtheta} may be obtained 
from Theorem~\ref{TheoremParametricBootstrap} under condition \eqref{EqJointConvergenceEPthetaRegular} at $\theta_0=\theta(H)$. 
(Given asymptotic linearity
\eqref{EqAsymptoticLinearity} and local asymptotic normality of the model $\{H_\theta: \theta\in\Theta\}$ at $\theta_0$,
the latter condition itself can be shown to be equivalent to $H_{\theta_0}( \psi_{H_{\theta_0}}\dot\ell_{\theta_0}^T)=I$,
in view of Le Cam's third lemma, where $\dot\ell_{\theta_0}$ is the score function of the model (see \cite{van2000asymptotic}, Chapters~6, and~7).)
\end{example}

\section{Independence testing}
\label{SectionIndependenceTesting}

As a first application, we study the classical problem of testing the independence of two random variables.
Consider a random pair $(X, Y)$ following a distribution $H$ on a product  measurable space $(\Xc \times \Yc,\Ac\times\Bc)$. Denote by $P$ and $Q$ the marginal distributions of $X$ and $Y$, respectively.
Assume that we observe an i.i.d.\ sample $(X_1, Y_1), \dots, (X_n, Y_n)$ from a distribution $H$.
We want to test for the independence of $X$ and $Y$, i.e.
\begin{equation}\label{eq_indep_hypothesis}
    \Hc_0: H = P \otimes Q
    \quad \text{versus} \quad 
    \Hc_1: H \neq P \otimes Q.
\end{equation}
We put this in the context of the testing problem \eqref{eq_hypothesis_problem} by
considering the map $\phi(H)= H - P \otimes Q$, where $P$ and $Q$ are the marginal distributions of $H$. 

To formalise this, let $\Fc$ and $\Gc$ be sets of  measurable functions  $f: \Xc\to\RR$ 
and $g: \Yc\to \RR$, and denote by $\Fc \otimes \Gc$ the set of functions
$f\otimes g: \Xc\times\Yc\to \RR$ given by 
 $(x,y) \mapsto f(x)g(y)$, when $f$ and $g$ vary over $\Fc$ and $\Gc$. Assume that the constant function
$x\mapsto 1$ and $y\mapsto 1$ are contained in $\Fc$ and $\Gc$, respectively, so that
$\Fc\otimes\Gc$ contains the functions  $1\otimes g$ and $f\otimes 1$, given by
$(x,y)\mapsto g(y)$ and $(x,y)\mapsto f(x)$, for all $f\in\Fc$ and $g\in\Gc$. Then consider $H$ and $P\otimes Q$ 
as elements of the space $\linftyFG$, and consider the map $\phi: \linftyFG\to\linftyFG $ defined by 
$$\phi(H)(f \otimes g)
= H(f \otimes g)- H(f \otimes 1) H(1 \otimes g)= H(f \otimes g) -Pf Qg=(H - P \otimes Q)(f\otimes g) .$$
Then $\phi(H)=0$, for every $H$ in the null hypothesis \eqref{eq_indep_hypothesis}, 
and the converse is true if the classes $\Fc$ and $\Gc$ are sufficiently rich.
The norm $\|\phi(H)\|_{\Fc\otimes\Gc}$ can be considered a distance of $H$ to independence.

With this notation the testing problem \eqref{eq_indep_hypothesis} is 
equivalent to \eqref{eq_hypothesis_problem}.  For $E=\linftyFG$, the test statistic
\eqref{EqNaturalTestStatistic} becomes 
\begin{align*}
    T_n
    &= \sqrt{n}\norm{\phi(\HHn)}_\FotimesG
    = \sqrt{n}\norm{\HHn
    - \PP_n \otimes \QQ_n}_\FotimesG.
\end{align*}
Two bootstrap resampling schemes are natural. The bootstrap under the null resamples from
$R_n=\PP_n\otimes \QQ_n$, where $\PP_n$ and $\QQ_n$ are the empirical distributions of 
$X_1,\ldots, X_n$ and $Y_1,\ldots, Y_n$, respectively.  The empirical bootstrap uses the joint empirical
distribution $R_n=\HHn$. Whereas the empirical
bootstrap resamples from the pairs $(X_i,Y_i)_{i = 1, \dots, n}$, the null bootstrap independently 
resamples observations from $(X_i)_{i=1}^n$ and $(Y_i)_{i=1}^n$.
In the two cases, the bootstrap test statistic \eqref{EqCorrectBootstrap} 
reduces to
\begin{align*}
       \Tneq
     &= \sqrt{n} \| \phi(\HHnst) - \phi(\PP_n\otimes\QQ_n) \|_\FotimesG
     = \sqrt{n}\norm{\HHnst - \PP_n^*\otimes \QQ_n^*}_\FotimesG,\\
 \Tncent
    &= \sqrt{n} \| \phi(\HHnst) - \phi(\HHn) \|_\FotimesG
    = \sqrt{n} \| \HHnst-\PP_n^*\otimes\QQ_n^* - (\HHn - \PP_n\otimes \QQ _n) \|_\FotimesG.
\end{align*}
These statistics are of the equivalent or centred types \eqref{eq_T_n_eq} and \eqref{eq_T_n_c}, respectively. We shall apply
the general results to see that both types of bootstraps are
consistent. On the other hand, each of the other two combinations, the non-centred statistic with the empirical bootstrap or the centred  statistic with the null bootstrap, are inconsistent (see Table~\ref{tab:indepdence_pairs} and
Section~\ref{proof:rem:unconsistent_combinations_indep}).

The consistency is a consequence of Corollary~\ref{cor_power},
and the fact that in both cases the sequence $R_n$ tends to $R(H)=H$ under
the null hypothesis, so that \eqref{EqBootstrapLimit} is trivially satisfied.
In fact, for the empirical bootstrap $R_n=\HHn\rightarrow H$, for every $H$,
whereas for the null bootstrap $R_n\rightarrow R(H)=P\otimes Q$, which is
equal to $H$ if $\phi(H)=0$. 

\renewcommand{\arraystretch}{1.5}

\begin{table}[H]
    \centering
    \begin{tabular}{l|cc}
    \multirow{2}{*}{
    \tikz{
    \node[below left, inner sep=1pt] (Rn) {$R_n$};
    \node[above right,inner sep=1pt] (Tnstar) {$\Tnst$};
    \draw (Rn.north west |-Tnstar.north west) -- (Rn.south east -| Tnstar.south east);
    }
    }
    & Centered & Equivalent \\ 
    & $
    \sqrt{n}\norm{\HHnst - \PP_n^*\otimes \QQ_n^* - (\HHn - \PP_n\otimes \QQ _n)}_\FotimesG$
    \hspace{1cm}
    & $
    \sqrt{n}\norm{\HHnst - \PP_n^*\otimes \QQ_n^*}_\FotimesG$
    \\
    \hline
    $\HHn$ & consistent & not consistent \\
    $\PP_n\otimes \QQ_n$ & not consistent & consistent
    \end{tabular}
    \caption{
    Consistency or lack thereof of combinations of resampling schemes $R_n$ and corresponding bootstrap test statistics $\Tnst$.
    }
    \label{tab:indepdence_pairs}
\end{table}

For a formal proof, we verify the conditions of
Theorem~\ref{thm:asymptotic_behavior_T}. Hadamard differentiability is
provided by the following lemma.

\begin{lemma}[Hadamard differentiability]
\label{lemma_phi_indep}
Let $\Fc$ and $\Gc$ be classes of  measurable functions  that contain the constant function 1, and let
$D_0$ be the set of measures with $\|H\|_{\Fc\otimes\Gc}+\|P\|_\Fc+\|Q\|_\Gc<\infty$. Then
the map $\phi: D_0\subset \linftyFG \to \linftyFG$ defined by $\phi(H)=H-P\otimes Q$ is 
uniformly Hadamard differentiable at any $H \in \linftyFG$ with derivative given by
$h\mapsto \bigl(f\otimes g\mapsto h(f\otimes g)-h(f\otimes 1)H(1\otimes g)-H(f\otimes 1)h(1\otimes g)\bigr)$.
\end{lemma}

\begin{proof}
Fix sequences $H_n\to H$ and $h_n\to h$ in $\linftyFG$, and $t_n\to0$ in $\RR$. By some algebra,
\begin{align*}
    \frac{\phi(H_n + t_nh_n) - \phi(H_n)}{t_n}(f\otimes g)
     &= h_n(f\otimes g)    - h_n(f\otimes1)H_n(1\otimes g)- H_n(f\otimes1)h_n(1\otimes g)\\
&\qquad\qquad- t_n h_n(f\otimes1)h_n(1\otimes g).
\end{align*}
The first term on the right converges to $h(f\otimes g)$ uniformly in $f$ and $g$, by the assumption that $h_n\to h$ in
 $\linftyFG$.
Because $1\in\Gc$, the functions $f\otimes 1$ are contained in $\Fc\otimes\Gc$, and hence
$h_n(f\otimes1)$ tends to $h(f\otimes 1)$, uniformly in $\Fc$, again by the assumption that $h_n\to h$ in
 $\linftyFG$. Similarly, the sequences $h_n(1\otimes g)$, $H_n(1\otimes g)$, and $H_n(f\otimes1)$ tend
to the limits $h(1\otimes g)$, $H(1\otimes g)$, and $H(f\otimes1)$ uniformly in $f$ and $g$. Since the product
$(a,b)\to ab$ is continuous, it follows that the right side of the display tends to
$h(f\otimes g)    - h(f\otimes1)H(1\otimes g)- H(f\otimes1)h(1\otimes g)$ uniformly in
$f$ and $g$. Thus $t_n^{-1}\bigl(\phi(H_n + t_nh_n) - \phi(H_n)\bigr)$ converges in
$\linftyFG$ to the derivative as given.

This limit is linear in $h$, and continuous by the same arguments as before.
\end{proof}

\begin{lemma}[Null bootstrap]
\label{LemmaIndependenceNullBootstrap}
Let $\Fc$ and $\Gc$ be separable classes of measurable functions 
that contain the constant function such that $\Fc\times\Gc$ satisfies the uniform
entropy condition for envelope functions $F$, $G$ and $F\otimes G$ that are $H$-square integrable.
Then
$R_n = \PP_n\otimes \QQ_n
\to R(H) = P \otimes Q$
in $\linftyFG$,
outer almost surely and the bootstrap empirical measure $\HHnst$ corresponding to $R_n$ satisfies
$\sqrt{n}(\HHnst - \PP_n \otimes\QQ_n)    \, | \, (X_1,Y_1) \dots, (X_n,Y_n)    \cvweakly    \GG_{P\otimes Q}$, almost surely.
\end{lemma}

\begin{proof}
Because the constant function is contained in $\Fc$ and $\Gc$, the classes of functions $\Fc\equiv\Fc\otimes 1$ and
$\Gc\equiv 1\times\Gc$ are subclasses of $\Fc\otimes\Gc$. Because the latter class 
satisfies the uniform entropy condition, so  do $\Fc$ and $\Gc$. Because they have integrable envelopes
and are suitably measurable, they are Glivenko-Cantelli. Thus
$(\PP_n \otimes \QQ_n)(f \otimes g) =
(\PP_n f)(\QQ_ng)$ 
satisfies
$\|\PP_n \otimes \QQ_n - P\otimes Q\|_{\Fc \otimes \Gc}
\le \|\PP_n-P\|_{\Fc} 
\|\QQ_n\|_{\Gc}
+ \|P\|_{\Fc} \|\QQ_n - Q\|_{\Gc}$, which tends to zero outer almost surely.

The second assertion on the convergence of
$\sqrt{n}(\HHnst - \PP_n \otimes \QQ_n)$
is Theorem~3.9.3 in \cite{van2023weak}.
\end{proof}

Lemmas~\ref{lemma_phi_indep} and~\ref{LemmaIndependenceNullBootstrap} show that the conditions of Theorem~\ref{thm:asymptotic_behavior_T} are satisfied for the null bootstrap $R_n = \PP_n\otimes \QQ_n$ with $R(H) = P \otimes Q$, for classes of functions $\Fc$ and $\Gc$ satisfying the conditions. Condition \eqref{EqBootstrapLimit}  is trivially satisfied, as
$R(H) = P \otimes Q$ is equal to $H$ under the null hypothesis.

Lemma~\ref{lemma_phi_indep} and Example~\ref{ExampleEmpiricalBootstrapConvergence}
show the same for the empirical bootstrap $R_n = \HH_n$ with $R(H) = H$, provided the class $\Fc \otimes \Gc$ is universally Donsker.

Thus in view of Corollary~\ref{cor_power}, both bootstrap procedures are consistent, for many examples of
classes of functions $\Fc$ and $\Gc$.

\begin{example}[Kolmogorov-Smirnov]\label{example_KS}
When $\Xc = \RR^p$, $\Yc = \RR^q$, and the classes $\Fc$ and $\Gc$ consist of the indicator functions of the cells $(-\infty,a]$,
for $a$ varying over $\RR^p$ or $\RR^q$, 
respectively, then the test statistic $T_n$ is the Kolmogorov-Smirnov statistic for independence
\begin{align*}
    T_n
    = \sqrt{n}
    \sup_{(x, y) \in \RR\rule{0pt}{0.6em}^{p+q}}
    \big| \FF_{(X,Y),n}(x , y)
    - \FF_{X,n}(x) \FF_{Y,n}(y) \big|,
\end{align*}
where $\FF_{(X,Y),n}$, $\FF_{X,n}$ and $\FF_{Y,n}$  are the empirical cumulative distribution functions of $(X,Y)$, $X$ and $Y$, respectively. These sets of functions satisfy the conditions of the preceding lemmas and hence lead to consistent bootstraps.
\end{example}

\begin{remark}
The results in this section can be generalised straightforwardly to the case of the joint independence test 
between $d$ random variables $X_1, \dots, X_d$. We  sketch this generalisation.
Assume that $\mathbf{X} = (X_1, \dots, X_d)$ follows the distribution $H$ on $\bigtimes_{j = 1}^d \Xc_j$ with marginal distributions $H_1, \dots, H_d$.   We observe an i.i.d.\ sample  $(\mathbf{X}_i    = (X_{i,1}, \dots, X_{i,d})    )_{i=1, \dots, n}$ from $H$.
Denote by $\HHn$ the empirical distribution of $\mathbf{X}$, and by $\HH_{j,n}$ the empirical distribution of the $j$-th marginal, for $j = 1, \dots, d$.  The null hypothesis of joint independence can be rewritten as $\phi(H) = 0$, where  $\phi(H) = H - \bigotimes_{j = 1}^d H_j$. The possible bootstrap schemes are $\HHn$ and $\bigotimes_{j = 1}^d \HH_{j,n}$, with bootstrap test statistics
\begin{align*}
    \Tncent
    = \sqrt{n} \bigg\|\HHnst    - \bigotimes_{j = 1}^d \HH_{j,n}^*  
    - \bigg(\HHn - \bigotimes_{j = 1}^d \HH_{j,n} \bigg)\bigg\|_{\bigotimes_{j = 1}^d\Fc_j},
\end{align*}
and
\begin{align*}
    \Tneq
    = \sqrt{n}\,\Bigl\|\HHnst  - \bigotimes_{j = 1}^d \HH_{j,n}^*  \Bigr\|_{\bigotimes_{j = 1}^d\Fc_j},
\end{align*}
for suitable sets $\Fc_j$ of real-valued measurable functions, respectively defined on $\Xc_j$.
\end{remark}

\section{Goodness-of-fit testing}
\label{sec:GoF_Setting}

In goodness-of-fit testing we wish to test the null hypothesis that the distribution $H$ of the observations $X_1,\ldots, X_n$ belongs to a given parametrised family $\{H_\theta:\theta\in\Theta\}$ of distributions. The null and alternative hypotheses are given by 
\begin{align*}
    \Hc_0: \:  H \in \{ \Htheta : \theta\in\Theta \}
    \quad\text{ versus }\quad
    \Hc_1: \:  H \notin \{ \Htheta : \theta\in\Theta \}.
\end{align*} 
A natural test statistic takes the form
\begin{align}
\label{EqGoodnessOfFitStatistics}
    T_n = \sqrt{n} \| \HHn - H_{\thetahat_n} \|_E,
\end{align}
for given estimators $\thetahat_n$ of the unknown parameter $\theta$ and some norm $\|\cdot\|_E$.
We can fit this in the setup of Section~\ref{sec:extension_framework_parameters}
with the map $\phi$ defined as $\phi(H,\theta)=H-H_\theta$.
Alternatively, we can use the setup of Section~\ref{sec:general_framework} if the
estimators $\thetahat_n$ can be viewed as Hadamard differentiable functionals 
$\hat\theta_n=\theta(\HHn)$ of $\HHn$, with the map $\phi$ defined by $\phi(H)=H-H_{\theta(H)}$.

\begin{example}[Kolmogorov-Smirnov]\label{example:kolmogorov-smirnov GOF}
Choosing $E=\linftyF$, for some class $\Fc$ of functions, gives the statistic $T_n=\sqrt{n}\| \HHn - H_{\hat\theta_n}\|_\Fc$.
In particular, for observations in $\RR^d$ and $\Fc$ the class of indicator functions of cells $(-\infty,a]$, for
$a\in\RR^d$, we find the classical Kolmogorov-Smirnov statistic $\sqrt n\|\FF_n-F_{\hat\theta_n}\|_\infty$, which is the
uniform distance between the cumulative distribution functions corresponding to the measures $\HHn$ and $H_{\hat\theta_n}$.
\end{example}

\begin{example}[Cram\'er-von Mises]
The weighted Cram\'er-von Mises statistic   $\int_{-\infty}^\infty \bigl(\FF_n(x) - F_{\hat\theta_n}(x)\bigr)^ 2 \,\mathrm{d}\mu(x)$ is the square
$L_2(\RR, \mu)$-distance between the empirical cumulative distribution function $\FF_n$ and the
cumulative distribution function $F_{\hat\theta_n}$ estimated according to the parametric model, for some measure $\mu$.
This fits the setup with $\|\cdot\|_E$ the $L_2(\RR,\mu)$-norm. 

The classical Cram\'er-von Mises  statistic uses $\mu=F_{\hat\theta_n}$ or $\mu=\FF_n$, which depend on $X_1,\ldots, X_n$, while
$d\mu(x)=dx/\bigl(1-F_{\hat\theta_n}(x)\bigr)$ gives the Anderson-Darling statistic. We could fit these in our general
setup by redefining $\phi$ as $\phi(H,\theta)=\int (H-H_\theta)^2\,dH$. 
\end{example}

The estimators $\hat\theta_n$ in \eqref{EqGoodnessOfFitStatistics} can take various forms. For instance, they might be maximum
likelihood estimators,  general $M$-estimators or minimum distance estimators. 
A special case are the minimum distance estimators
relative to the criterion used to define the measure of fit $T_n$, i.e.\ 
\begin{equation}
    \label{EqMDestimator}
\hat\theta_n^{MD}=\argmin_\theta \| \HHn - H_{\theta} \|_E.
\end{equation}
In this case the test statistic \eqref{EqGoodnessOfFitStatistics} reduces to the
distance $\inf_\theta\|\HHn-H_\theta\|_E$ of the empirical measure to the parametric model.

It is natural to perform a bootstrap under the null hypothesis, which presently 
means that the bootstrap observations are constructed as an  i.i.d.\ sample from the estimated distribution
$R_n=H_{\hat\theta_n}$. We shall see that also
the empirical bootstrap $R_n=\HH_n$ gives correct results, provided the bootstrap 
test statistic is constructed properly, according to \eqref{EqCorrectBootstrapTheta}, with corresponding
bootstrap estimators $\theta_n^*$. In the present situation the latter definition reduces to 
\begin{align}
    \Tnst = \sqrt{n} \| \phi(\HHnst, \theta_n^*) - \phi(R_n,\hat\theta_n) \|_E
    = \sqrt{n} \|
    \HHnst - H_{\theta_n^*} - R_n + H_{\hat\theta_n} \|_E .\label{eq_Tnst_GoF}
\end{align}
Here $\theta_n^*$ are appropriately constructed bootstrap versions of the estimator.
For the null bootstrap $R_n=H_{\hat\theta_n}$  and empirical bootstrap $R_n=\HHn$,
the statistic \eqref{eq_Tnst_GoF} reduces to the equivalent and centered statistics, respectively, given by
    \begin{align*}
        \Tneq &:= \sqrt{n} \| \HHnst - H_{\theta_n^*}\|_E,\\
        \Tncent &:= \sqrt{n} \| \HHnst - H_{\theta_n^*}- \HH_n +  H_{\hat\theta_n} \|_E.
    \end{align*}
The term $R_n-H_{\hat\theta_n}$ in the general statistic \eqref{eq_Tnst_GoF} may be viewed as a ``correction term'',
which vanishes in the case of the null bootstrap, but not for the empirical bootstrap.

The bootstrap values $\theta_n^*$ used in \eqref{eq_Tnst_GoF} may also need correction, depending on the choice of $\hat\theta_n$ and the
bootstrap scheme. To construct $\theta_n^*=\hat\theta_n(X_1^*,\ldots, X_n^*)$ from the bootstrap values
in the same way as $\hat\theta_n=\hat\theta_n(X_1,\ldots, X_n)$ is defined as a function of the original observations may seem
natural. If $\hat\theta_n=\theta(\HH_n)$ is a Hadamard differentiable function of the empirical distribution,
then this will work, as in that case the goodness-of-fit statistic is the norm of a
Hadamard differentiable function of $\HH_n$ and Theorem~\ref{thm:asymptotic_behavior_T} applies under just conditions on $\HHn$ and $R_n$.
However, in the case of the minimum distance estimator \eqref{EqMDestimator}, this
would lead to $\theta_n^*=\argmin_\theta \| \HHnst - H_{\theta} \|_E$, which is not necessarily a good choice.
We shall see that for bootstrap $R_n$ a correct match to  \eqref{EqMDestimator} is
\begin{equation}
\label{EqBootstrapMDEstimator}
\theta_n^{*,MD}=\argmin_\theta \| \HHnst - H_{\theta}- R_n+H_{\hat\theta_n} \|_E.
\end{equation}
Unless $R_n$ follows the null bootstrap, this includes the correction term $R_n-H_{\hat\theta_n}$, and because of this
$\theta_n^*=\theta_n^*(X_1^*,\ldots, X_n^*; X_1,\ldots, X_n)$ depends on both the bootstrap values and the original observations.
Some intuition for this choice is that it turns the bootstrap value $\Tnst$ in \eqref{eq_Tnst_GoF} into the corrected (centered) distance 
$\argmin_\theta \| \HHnst - H_{\theta} -R_n+H_{\hat\theta_n}\|_E$ to the model. 

It is known that minimum distance estimators based on the Kolmogorov-Smirnov distance
lack robustness for small deviations from the model (\cite{DonohoLiu1988pathologies}), in contrast to minimum distance
estimators based on smoother distances (see \cite{Millar1984}, \cite{DonohoLiu1988}). This may explain our finding below, that a correction
is needed in the bootstrap scheme for the Kolmogorov-Smirnov distance, but is not essential for the Cram\'er-von Mises
statistic.

We finish this general discussion with a lemma on the Hadamard differentiability of the present functional $\phi$.
Consider probability measures $H$ as elements of  appropriate normed spaces $D$ and $E$.

\begin{lemma}[Hadamard differentiability]
\label{LemmaGoodnessHadamard}
Let $\Theta$ be  an open set in $\RR^d$ and suppose that the map $\theta\mapsto H_\theta$ 
from $\Theta$ to $E$ is continuously (Hadamard) differentiable with derivative $\dot H_\theta: \RR^d\to E$. Furthermore,
assume that $\|H\|_E\le C \|H\|_D$, for some constant $C$ and every $H$. Then
the map $\phi: D\times\Theta\to E$ defined 
by $\phi(H,\theta)=H-H_\theta$ is uniformly Hadamard differentiable at every $(H,\theta)\in D\times\Theta$.
The derivative is given by $(g,h)\mapsto g-\dot H_\theta h$.
\end{lemma}

\begin{proof}
The map $H\mapsto H$ from $D$ to $E$ is linear, and continuous by assumption, and hence is uniformly differentiable.
The map  $\theta\mapsto H_\theta$ from $\Theta$ to $E$ is continuously differentiable by assumption and hence
uniformly differentiable.
Therefore the map $(H,\theta)\mapsto (H,H_\theta)$ from $D\times\Theta$ to $E\otimes E$ is uniformly differentiable.
The map $(H,G)\mapsto H-G$ from $E\otimes E$  to $E$ is linear and continuous and hence uniformly differentiable.
Thus the result follows by the chain rule.
\end{proof}

\subsection{Minimum distance estimators}
\label{SectionMinimumDistanceInformal}

The minimum distance estimator \eqref{EqMDestimator} renders the goodness-of-fit statistic \eqref{EqGoodnessOfFitStatistics} into
the distance of the empirical measure to the model, and hence is a natural choice. In this section we give an informal discussion
of the behaviour of this estimator and its bootstrap version \eqref{EqBootstrapMDEstimator}, leaving precise results to Appendix~\ref{appendix:minimum-distance-estimators}.
(Rigorous results for the minimum distance estimator itself go back to at least \cite{Millar1984}.)
We assume that $X_1,\ldots,X_n$ are an i.i.d.\ sample from some distribution $H$, with special attention for measures
$H=H_{\theta_0}$ belonging to the null hypothesis.

As $\|\HHn-H_\theta\|_E\to \|H-H_\theta\|_E$ almost surely, the minimum distance estimator \eqref{EqMDestimator}  will typically
tend to the point of minimum $\theta(H)=\argmin_\theta\|H-H_\theta\|_E$. More refined properties of \eqref{EqMDestimator}
depend on the distance $\|\cdot\|_E$, but for $H=H_{\theta_0}$  from the null hypothesis a general analysis is possible.

Under $H=H_{\theta_0}$ the minimum distance estimator inherits 
the $\sqrt n$-consistency of the empirical distribution (relative to $\|\cdot\|_E$) under general
identifiability conditions (see Lemma~\ref{lemma_argmin}). Thus
we can focus on the asymptotic behaviour of the sequence
$\sqrt n(\hat\theta_n^{MD}-\theta_0)$. An informal derivation of its limit behaviour is (see Lemma~\ref{lemma:joint_cv_H_theta} for a precise statement)
\begin{align*}
\sqrt n(\hat\theta_n ^{MD}-\theta_0)
&=\argmin_{h}
\|\HHn-H_{\theta_0+h/\sqrt n}\|_E, \\
&\doteq \argmin_{h\in\RR^d}
\|\sqrt n(\HHn-H_{\theta_0}) - \dot H_{\theta_0}h\|_E
\cvweakly \argmin_{h\in\RR^d}
\|\GG_{H_{\theta_0}}-\dot H_{\theta_0}h\|_E.
\end{align*}
Here $\GG_{H_{\theta_0}}$ is the weak limit in $E$ of the sequence $\sqrt n(\HHn-H_{\theta_0})$ and $\dot H_{\theta_0}: \RR^d\to E$ is the derivative of $\theta\mapsto H_\theta$ at $\theta_0$.
This derivation remains valid if $\theta_0$ is replaced by a converging sequence $\theta_{0,n}\to\theta_0$ (and $H_{\theta_0}$ by the corresponding sequence of
measures $H_{\theta_{0,n}}$), with the same limit variable. In particular, the minimum distance estimator is regular in the sense of \eqref{EqJointConvergenceEPthetaRegularNoRate}.
The joint convergence of the sequence $\sqrt n(\HH_n-H_{\theta_{0,n}},\hat\theta_n^{MD_n}-\theta_{0,n})$ under $H_{\theta_{0,n}}$ can be proved too,
and hence \eqref{EqJointConvergenceEPthetaRegularNoRate} is true (under the conditions given in Lemma~\ref{lemma:joint_cv_H_theta}). 
This verifies the most important condition for consistency of the parametric bootstrap $R_n=H_{\hat\theta_n^{MD}}$.

For the Kolmogorov-Smirnov distance, the limit variable in the preceding display is typically intractable, and non-normal. For a Hilbertian norm
$\|\cdot\|_E$, the point of minimum $\hat h$ can be computed from the stationary equation $\langle \GG_{H_{\theta_0}}-\dot H_{\theta_0}\hat h, \dot H_{\theta_0}\rangle_E=0$,
obtained by differentiation of the squared distance 
$\|\GG_{H_{\theta_0}}-\dot H_{\theta_0}h\|_E^2$. 
This gives
$$\sqrt n(\hat\theta_n^{MD}-\theta_0)\cvweakly \langle \dot H_{\theta_0},\dot H_{\theta_0}^T\rangle_E^{-1}\,\langle \GG_{H_{\theta_0}}, \dot H_{\theta_0}\rangle_E.$$
In this formula $\dot H_{\theta_0}$ is identified with the vector $(\dot H_{\theta_0}e_1,\ldots, \dot H_{\theta_0}e_d)^T$ in $E^d$ 
obtained by evaluating the derivative $\dot H_{\theta_0}: \RR^d\to E$ at the unit vectors in $\RR^d$ (the gradient of $\theta\mapsto H_\theta$). 
The limit variable is linear in $\GG_{H_{\theta_0}}$ and hence normally distributed.

For the bootstrap version \eqref{EqBootstrapMDEstimator}, a similar analysis suggests that, under every $H$ such that
$\sqrt n(\HHnst-R_n)| X_1,\ldots, X_n\cvweakly \GG_{R(H)}$ and $\hat\theta_n\to \theta(H)$ for some $\theta(H)$, both in probability,
\begin{align*}
\sqrt n(\theta_n^{*,MD}-\hat\theta_n)&=\argmin_h\|\HHnst-R_n-(H_{\hat \theta_n+h/\sqrt n}-H_{\hat\theta_n})\|_E,\\
&\doteq \argmin_{h\in\RR^d}\|\sqrt n(\HHnst-R_n)-\dot H_{\hat\theta_n}h\|_E\cvweakly \argmin_{h\in\RR^d}\|\GG_{R(H)}-\dot H_{\theta(H)}h\|_E,
\end{align*}
conditionally given $X_1,\ldots, X_n$, in probability. This argument works for general estimator sequences $\hat\theta_n$
and applies to any distribution $H$, also from the alternative hypothesis. The joint convergence \eqref{EqJointConvergenceBootstrapEPtheta} 
is often true too, for $H$ from both hypotheses (see Lemma~\ref{lemma_convergence_theta_rescaled}).
Equality $R(H)=H$ is sufficient for the limit variables of the test statistic
and its bootstrap version to agree. This is true for $H=H_{\theta_0}$ in the null hypothesis for both the null and empirical bootstraps $R_n$.

The correction term $R_n-H_{\hat\theta_n}$ in \eqref{EqBootstrapMDEstimator} is important in the preceding analysis. 
For the bootstrap value $\theta_n^*=\argmin_\theta \|\HHnst-H_\theta\|_E$ without correction, the analogous analysis is
\begin{align*}
\sqrt n(\theta_n^*-\hat\theta_n)&=\argmin_h\|\HHnst-H_{\hat \theta_n+h/\sqrt n}\|_E,\\
&= \argmin_{h}\|\sqrt n(\HHnst-R_n)+\sqrt n(R_n-H_{\hat\theta_n})-\sqrt n(H_{\hat\theta_n+h/\sqrt n}-H_{\hat\theta_n})\|_E.
\end{align*}
Typically the process $\sqrt n(\HHnst-R_n)$ converges conditionally given $X_1,\ldots, X_n$, but the process $\sqrt n(R_n-H_{\hat\theta_n})$
is a function of these observations and may remain random. (The null bootstrap is an exceptional case, where 
this term is zero.)  This appears to cause divergence in the case of the Kolmogorov-Smirnov statistic. However, for a Hilbertian distance and the empirical bootstrap, the extra
term may be harmless if the estimator $\hat\theta_n=\hat\theta_n^{MD}$ is the minimum distance estimator. In that case the right side of the preceding display is equal to 
\begin{align*}
&\argmin_{h}\Bigl[\bigl\|(\HHnst-R_n)-(H_{\hat\theta_n+h/\sqrt n}-H_{\hat\theta_n})\bigr\|_E^2+\|R_n-H_{\hat\theta_n}\|_E^2\\
&\qquad\qquad\qquad+2\bigl\langle (\HHnst-R_n)
-(H_{\hat\theta_n+h/\sqrt n}-H_{\hat\theta_n}),R_n-H_{\hat\theta_n}\bigr\rangle_E\Bigr].
\end{align*}
The second term does not depend on $h$ and hence can be omitted from the argmin. For the same reason the inner product in the third term can be reduced
to $-\langle H_{\hat\theta_n+h/\sqrt n}-H_{\hat\theta_n},R_n-H_{\hat\theta_n}\rangle_E$. 
If $R_n=\HH_n$ and $\hat\theta_n=\hat\theta_n^{MD}$, then
the linear approximation $-\langle \dot H_{\hat\theta_n}h,R_n-H_{\hat\theta_n}\rangle_E/\sqrt n$ to this term
vanishes for all $h$, by the stationary equation for the minimisation problem that defines $\hat\theta_n^{MD}$. In that case the preceding display is asymptotic to, for $\ddot H_\theta$ the second derivative of
$\theta\mapsto H_\theta$,
$$\argmin_{h\in\RR^d}\Bigl[\bigl\|\sqrt n(\HHnst-R_n)-\dot H_{\theta_0}h\bigr\|_E^2- \langle h^T\ddot H_{\theta_0}h,R_n-H_{\theta_0}\rangle_E\Bigr].$$
If $R_n\to H_{\theta_0}$, then the second term is negligible and the resulting expression is the same term as obtained with the corrected bootstrap minimum distance estimator \eqref{EqBootstrapMDEstimator}.
(See Lemma~\ref{LemmaMDEquivalenceCentered} for a precise expression of this argument.) 

Thus the correction term is important for the Kolmogorov-Smirnov distance, but may be omitted for the combination of a Hilbertian distance and the empirical bootstrap. We also see that under the alternative hypothesis, when $R_n\to H$ not of the form $H_{\theta_0}$, the sequence $\sqrt n(\theta_n^*-\hat\theta_n)$ will typically converge in distribution, but with a different variance
than the sequence $\sqrt n(\theta_n^{*,MD}-\hat\theta_n)$ (see Lemma~\ref{LemmaMDEquivalenceCentered}).

A similar argument shows that, for a Hilbertian norm,
the minimum distance estimator \eqref{EqMDestimator} will be $\sqrt n$-consistent for the point of minimum $\theta(H)=\argmin_\theta\|H-H_\theta\|$ and be asymptotically normal, also for $H$ from the alternative hypothesis.
The joint convergence \eqref{EqJointConvergenceEPtheta} will typically also be valid at most distributions $H$ from both hypotheses (see Lemma~\ref{LemmaMDHilbertian}).

\subsection{Parametric bootstrap}

Consistency of the parametric bootstrap $R_n=H_{\hat\theta_n}$ together with the natural bootstrap values  $\theta_n^*=\hat\theta_n(X_1^*,\ldots, X_n^*)$
can be obtained by combining Theorems~\ref{TheoremBootstrapExtension} and~\ref{TheoremParametricBootstrap}. For simplicity, assume $D=E$.
The phrase ``every $H$'' in the following refers to every probability distribution $H$ under consideration, typically all probability distributions.

\begin{lemma}[Null bootstrap]
\label{LemmaGoodnessNull}
Let $\Theta$ be  an open set in $\RR^d$, suppose that the map $\theta\mapsto H_\theta$ from $\Theta$ to $E$ is continuously Hadamard differentiable,
and assume that $\inf_{\theta\in\Theta}\|H-H_\theta\|_E>0$, for every $H\in\Hc_1$.
Let $R_n=H_{\hat\theta_n}$ and $\theta_n^*=\hat\theta_n(X_1^*,\ldots, X_n^*)$ for a sequence of estimators
$\hat\theta_n$ such that for every $H$ there exist $\theta(H)\in\Theta$ such that \eqref{EqJointConvergenceEPtheta} holds
 and such that \eqref{EqJointConvergenceEPthetaRegular} holds with $\theta_0=\theta(H)$. 
Then $P(T_n \geq \xi^*_{n, 1-\alpha})\rightarrow \alpha$ if $H\in\Hc_0$ and $P_\inner(T_n > \xi^*_{n, 1-\alpha}) \rightarrow 1$
if $H\in \Hc_1$. This remains valid if \eqref{EqJointConvergenceEPtheta} holds for $H\in\Hc_0$
and \eqref{EqJointConvergenceEPthetaRegularNoRate} holds and $\hat\theta_n\to\theta(H)$, in probability, for every $H$.
\end{lemma}

\begin{proof}
Condition \eqref{EqJointConvergenceEPtheta} is satisfied by assumption and implies that $\sqrt n\bigl(\hat\theta_n-\theta(H)\bigr)=O_P(1)$.
In view of Theorem~\ref{TheoremParametricBootstrap}, assumption \eqref{EqJointConvergenceEPthetaRegular} 
implies condition \eqref{EqJointConvergenceBootstrapEPtheta} of  Theorem~\ref{TheoremBootstrapExtension}
with $R(H)=H_{\theta_0}$, where $\theta_0=\theta(H)$, and also gives that $R_n\to R(H)$, in probability.
The assumption $\inf_{\theta\in\Theta}\|H-H_\theta\|_E>0$ ensures that $\Hc_1$ is indeed the set of $H$ with $\phi\bigl(H,\theta(H)\bigr)\not=0$.
Therefore the assertions follow from Theorem~\ref{TheoremBootstrapExtension}.

For the final assumption we apply the triangle inequality to see that $T_n\ge \sqrt n\|H-H_{\hat\theta_n}\|_E-\sqrt n\|\HH_n-H\|_E$, where
the first term is bounded below by $\sqrt n\inf_\theta\|H-H_\theta\|_E\to\infty$ and the second is bounded in probability, so that 
\eqref{EqPowerTestWithParameter}  is valid, and we can apply the preciser part of Theorem~\ref{TheoremBootstrapExtension}.
\end{proof}

Condition \eqref{EqJointConvergenceEPthetaRegularNoRate} is verified for general minimum distance estimators $\hat\theta_n$ in Lemma~\ref{lemma:joint_cv_H_theta}. 
Condition \eqref{EqJointConvergenceEPtheta} is verified for general minimum distance estimators for $H$ in the null hypothesis
in Lemma~\ref{lemma:joint_cv_H_theta} and for general $H$ for Hilbertian minimum distance estimators in 
Lemma~\ref{LemmaMDHilbertian}. 

For estimators that are Hadamard differentiable functionals $\hat\theta_n=\theta (\HHn)$ and $E=\linftyF$ with a Donsker class $\Fc$, 
conditions \eqref{EqJointConvergenceEPtheta} and \eqref{EqJointConvergenceEPthetaRegular} follow from the regularity of the empirical distribution 
(see \cite{van2023weak}, Section~3.12.1).

This covers many examples.

\subsection{Empirical bootstrap}

If $\hat\theta_n=\theta (\HHn)$ and $\theta_n^*=\theta (\HHnst)$ for a Hadamard differentiable map $\theta$
and $D=\linftyF$ for a Donsker class $\Fc$ or $D=L_2(\RR^p,\mu)$, then Theorem~\ref{thm:asymptotic_behavior_T} applies,
where the conditions are verified in Examples~\ref{ExampleEmpiricalBootstrapConvergence} and~\ref{ExampleEmpiricalBootstrapConvergenceBanachSpace}.

More generally, we verify conditions \eqref{EqJointConvergenceEPtheta}  and \eqref{EqJointConvergenceBootstrapEPtheta}.
Here we may combine the convergence of the bootstrap empirical process, as before,
with asymptotic linearity as explained in Example~\ref{ExampleAsymptoticLinearity},
where it may help to choose $D=\ell^\infty(\Fc\cup\{\psi\})$ with the classs $\Fc$ from $E=\ell^\infty(\Fc)$ enlarged with the
influence functions of the estimators $\hat\theta_n$. 
For minimum distance estimators the desired results can be found in  Lemmas~\ref{lemma:joint_cv_H_theta} and~\ref{LemmaMDHilbertian}.

\section{Testing the slope in linear regression}
\label{SectionLinearRegression}

Suppose that we observe an i.i.d.\ sample $(X_1,Y_1),\ldots,(X_n,Y_n)$ of observations following the regression model
\begin{align}
\label{EqRegressionModel}
    Y_i&=a+bX_i+\epsilon_i,\qquad   \Expec \epsilon_i=\Expec \epsilon_i  X_i=0.
\end{align}
We are interested in testing the hypothesis that the slope of the regression line is zero:
\begin{equation}\label{eq_hypothesis_test}
    \Hc_0:\: b=0,\quad \text{versus} \quad \Hc_1:\: b\neq0.
\end{equation}
Under the given condition on the errors, the slope is identified as 
$$b=\frac{\Cov(X_i,Y_i)}{\Var X_i}= \frac{H (xy)- P(x) Q(y)}{ P(x^2)-(Px)^2},$$
for $H$ the distribution of $(X_i,Y_i)$, $P$ and $Q$ the marginal distributions of $X_i$ and $Y_i$, respectively, and terms such as $H(xy)$ are understood as $\Expec_H[X Y]$.
Thus the testing problem is of type \eqref{eq_hypothesis_problem}, for $\phi(H)$ defined by
\begin{equation}
\label{EqPhiRegression}
\phi(H)=\frac{H (xy)- P(x) Q(y)}{ P(x^2)-(Px)^2}.
\end{equation}
For $\HH_n$ the empirical distribution of the pairs $(X_1,Y_1),\ldots,(X_n,Y_n)$ and $\PP_n$ and $\QQ_n$ the empirical distributions of $X_1,\ldots, X_n$ and $Y_1,\ldots, Y_n$, respectively, 
the test statistic \eqref{EqNaturalTestStatistic} becomes
$$T_n=\sqrt{n}\, |\phi(\HHn)| =\frac{|\HHn(xy)-\PP_n(x)\QQ_n(y)|}{\PP_n(x^2)-\PP_n(x)^2}
=\frac{|\sum_{i=1}^n(X_i-\bar X_n)(Y_i-\bar Y_n)|}{\sum_{i=1}^n(X_i-\bar X_n)^2}.$$
The statistic $T_n$ is simply the absolute value of the least squares estimator for $b$.
Various bootstrap schemes are possible, each forming bootstrap pairs $(X_1^*,Y_1^*),\ldots, (X_n^*,Y_n^*)$.
\begin{enumerate}[label=(\roman*)]
\item The null bootstrap $R_n=\PP_n\otimes \QQ_n$ independently resamples $X_1^*,\ldots,X_n^*$ from $X_1,\ldots,X_n$
and $Y_1^*,\ldots,Y_n^*$ from $Y_1,\ldots,Y_n$.
\item The empirical bootstrap $R_n=\HH_n$ resamples the pairs $(X_1^*,Y_1^*),\ldots, (X_n^*,Y_n^*)$ from
the pairs  $(X_1,Y_1),\ldots,(X_n,Y_n)$.
\end{enumerate}
Under the null hypothesis, the regression equation gives $Y_i=a+\epsilon_i$, which is not necessarily independent of $X_i$ under the (minimal) error conditions in the regression model \eqref{EqRegressionModel}.
Therefore, the terminology ``null bootstrap'' is misleading, although it 
would be appropriate under the commonly made assumption that $\epsilon_i$ and $X_i$ are independent. We shall see that
the empirical bootstrap is correct under the general model \eqref{EqRegressionModel}, whereas the null bootstrap
is only correct given additional moment restrictions (which are implied by but weaker than independence of $\epsilon_i$ and $X_i$).

For the bootstraps (i) and (ii) the test statistic \eqref{EqCorrectBootstrap} becomes
\begin{align*}
\Tneq&=\sqrt n\,|\phi(\HH_n^*)-\phi(\PP_n\otimes \QQ_n)|
= \frac{|\HHnst(xy)-\PP_n^*(x)\QQ_n^*(y)|}{\PP_n^*(x^2)-(\PP_n^*x)^2}
=\frac{|\sum_{i=1}^n(X_i^*-\bar X_n^*)(Y_i^*-\bar Y_n^*)|}{\sum_{i=1}^n(X_i^*-\bar X_n^*)^2},\\
  \Tncent &= \sqrt n\,|\phi(\HH_n^*)-\phi(\HHn)|= \biggl|\frac{\HHnst(xy)-\PP_n^*(x)\QQ_n^*(y) }{\PP_n^*(x^2)-(\PP_n^*x)^2}
-\frac{\HHn(xy)-\PP_n(x)\QQ_n(y) }{\PP_n(x^2)-(\PP_nx)^2}\biggr|.
\end{align*}
These are of the equivalent and centered types, as in \eqref{eq_T_n_eq} and \eqref{eq_T_n_c}.
Let $\Fc$ be the set of four functions $(x,y)\mapsto xy$, $(x,y)\mapsto x^2$,  $(x,y)\mapsto x$ and $(x,y)\mapsto y$, and consider the measure $H$ to be an element of $D = \linftyF$. Then the map
$\phi: \linftyF\to \RR$ is Hadamard differentiable.

\begin{lemma}[Hadamard differentiable]
\label{LemmaSlopeHadamard}
Let $\Fc$ be the set of four functions as indicated. Then the  map $\phi: \linftyF\to \RR$ defined in \eqref{EqPhiRegression} is uniformly Hadamard differentiable at every 
$H$ such that $H(x^2)-(Hx)^2>0$.
\end{lemma}

\begin{proof}
Each of the four maps $H\mapsto Hf$, for $f\in \Fc$, is linear and continuous and hence uniformly differentiable.
The map $\phi$ composes the four maps by the map $(a,b,c,d)\mapsto (a-cd)/(b-c^2)$, which is continuously
differentiable provided $b-c^2$ is bounded away from zero. The lemma follows by the chain rule.
\end{proof}

\begin{lemma}[Null bootstrap]
\label{LemmaSlopeNull}
Let $\Fc$ be the set of four functions as indicated and assume that $Hf^2<\infty$, for every $f\in\Fc$.
Then  $R_n=\PP_n\otimes \QQ_n\to P\otimes Q$ in $\linftyF$, almost surely, and the bootstrap empirical measure $\HHnst$ corresponding to $R_n$ satisfies
$$\sqrt{n}(\HHnst - \PP_n\otimes\QQ_n)    \, | \, (X_1,Y_1) \dots, (X_n,Y_n)    \cvweakly    \GG_{P\otimes Q},\quad \text{ in }\linftyF,
\text{ almost surely}.$$
\end{lemma}

\begin{proof}
The almost sure convergence $R_n\to P\times Q$ in $\linftyF$ is equivalent to $R_n f \to (P \times Q) f$,  almost surely, for every of the four functions $f\in\Fc$.  This follows from the
law of large numbers. For instance, for $f(x,y)=xy$, the assertion becomes $R_n(xy)=\PP_n(x)\QQ_n(y)\to P(x)Q(y)$,
almost surely.

Similarly, the conditional convergence in distribution follows from the multivariate (Lindeberg) central limit theorem.
Under $R_n$, given $ (X_1,Y_1) \dots, (X_n,Y_n) $, the variables 
$\sqrt{n}(\HHnst - \PP_n\otimes\QQ_n)f=n^{-1/2}\sum_{i=1}^n\bigl(f(X_i^*,Y_i^*)-(\PP_n\otimes\QQ_n)f\bigr)$ 
are centered at mean zero and hence it suffices to compute the covariance for every pair of $f\in\Fc$ and
verify the Lindeberg condition for every $f\in\Fc$. For instance, for the
function $f(x,y)=xy$, the Lindeberg condition becomes convergence to zero for every $\eta>0$ of the moments
$$\Expec_{R_n} \bigl((X_i^*Y_i^*)^21_{|X_i^*Y_i^*|>\eta\sqrt n}|  (X_1,Y_1) \dots, (X_n,Y_n)\bigr)
=\frac1n\sum_{i=1}^n (X_iY_i)^21_{|X_iY_i|>\eta \sqrt n}.$$
The right side tends to zero almost surely by the law of large numbers. Similarly, the covariances can be shown
to tend to the covariances of the Brownian bridge process $\GG_{P\otimes Q}$, almost surely. For instance,
$\Cov_{R_n}\bigl(X_i^*Y_i^*,Y_i^*\bigr)=\PP_n(x)\QQ_n(y^2)-\PP_n(x)\bigl(\QQ_n(y))^2$ tends almost surely to
$P(x)Q(y^2)-P(x)(Q(y))^2=\Cov\bigl(\GG_{P\otimes Q}(xy),\GG_{P \otimes Q}(y)\bigr)$.
\end{proof}

Lemmas~\ref{LemmaSlopeHadamard} and~\ref{LemmaSlopeNull} show that the conditions of Theorem~\ref{thm:asymptotic_behavior_T} are satisfied for the
null bootstrap $R_n=\PP_n\times \QQ_n$. Condition \eqref{EqBootstrapLimit} for correct type 1 error of the bootstrap becomes 
$|\phi'_{P\otimes Q}(\GG_{P \otimes Q})|\sim |\phi'_H(\GG_H)|$, under the null hypothesis. If the errors $\epsilon_i$
are independent from the $X_i$, then $Y_i$ and $X_i$ are also independent under the null hypothesis, 
and hence $H=P\otimes Q$, so that \eqref{EqBootstrapLimit}  is trivially satisfied. Under the less strict
condition $\Expec \epsilon_i=\Expec \epsilon_i  X_i=0$ imposed in \eqref{EqRegressionModel}, this is not necessarily the case. However,
the two processes $\GG_{P \otimes Q}$ and  $\GG_H$ involved are the four-dimensional vectors
$$\bigl(\GG_{P\otimes Q}(xy), \GG_{P\otimes Q}(x^2), \GG_{P\otimes Q}(x), \GG_{P\otimes Q}(y)\bigr)\quad\text{ and }\quad
\bigl(\GG_H(xy), \GG_H(x^2), \GG_H(x), \GG_H(y)\bigr).$$
Both vectors possess multivariate normal distributions with mean zero. If their covariance matrices are equal,
then they are equal in distribution and hence \eqref{EqBootstrapLimit}  is satisfied. Because marginal quantities
are the same under $H$ and $P\otimes Q$, this equality can be reduced to equality of the mixed moments
$\Cov(X_iY_i,X_i^2)$, $\Cov(X_iY_i,X_i)$, $\Cov(X_iY_i,Y_i)$, $\Cov(X_i^2,Y_i)$ and $\Cov(X_i,Y_i)$.
Here the last one is zero under both measures, by \eqref{EqRegressionModel}, and under the null hypothesis
$Y_i=a+\epsilon_i$, equality of the first four moments can be reduced to equality of these
same moments but with $Y_i$ replaced by $\epsilon_i$, since marginal moments are the same.
This can finally be reduced to equality of 
$\Cov(X_i^3,\epsilon_i)$, $\Cov(X_i^2,\epsilon_i)$ and $\Cov(X_i,\epsilon_i^2)$. These all vanish under
$P\otimes Q$, and hence we conclude that \eqref{EqBootstrapLimit}  is satisfied for every $H$ under the null
hypothesis such that $\Cov_H(X_i^3,\epsilon_i)=\Cov_H(X_i^2,\epsilon_i)=\Cov_H(X_i,\epsilon_i^2)=0$,
next to the assumption $\Cov_H(X_i,\epsilon_i)=0$ already in place in \eqref{EqRegressionModel}.

Thus the null bootstrap is consistent under the latter moment condition, which is considerably
weaker than independence of $X_i$ and $\epsilon_i$. In particular, the mixed moment
condition is satisfied if $\Expec_H(\epsilon_i | X_i)=0$ and $\Var_H(\epsilon_i | X_i)$ does not
depend on $X_i$.

The condition can be further relaxed by not requiring
equality in distribution of the two full vectors in the preceding display, but only of the
induced variables $\phi'_{P\otimes Q}(\GG_{P \otimes Q})$ and $\phi'_H(\GG_H)$. As these are both one-dimensional
Gaussian with zero mean, this reduces to a single equation. However, this is harder to interpret, by its
dependence on $H$.

Lemma~\ref{LemmaSlopeHadamard} and Example~\ref{ExampleEmpiricalBootstrapConvergence}
show that the conditions of Theorem~\ref{thm:asymptotic_behavior_T} are
satisfied for the empirical bootstrap $R_n=\HH_n$, with $R(H)=H$. By the latter equation,
condition \eqref{EqBootstrapLimit} is trivially satisfied for every distribution 
(with $\Expec \epsilon_i=\Expec \epsilon_i  X_i=0$ as in \eqref{EqRegressionModel}). Thus the empirical bootstrap has a
wider range of application.

\begin{remark}
The choice of the least squares estimator as test statistic leads to the map \eqref{EqPhiRegression}. The test based
on the covariance corresponds to the map $\phi(H)=H(xy)-P(x)Q(y)$. An analogous analysis then shows that the null
bootstrap is consistent under every $H$ such that $\Cov_H(X_i^2,\epsilon_i)=\Cov_H(X_i,\epsilon_i^2)=0$, thus requiring
one fewer mixed moment.
\end{remark}

\subsection{Studentised test statistic}
It is common to standardise the test statistic so that its variance is approximately one. The standardised least squares estimator 
corresponds to the map 
\begin{equation}
\label{EqPhiRegressionStudentised}
\phi(H)=\frac{H (xy)- P(x) Q(y)}{ \sqrt{P(x^2)-(Px)^2}}.
\end{equation}
The test statistic $T_n=\sqrt{n}\, |\phi(\HHn)| $ of \eqref{EqNaturalTestStatistic}  is the studentised least squares estimator
and is asymptotically standard normal under the null hypothesis. In practice, quantiles from the $t$-distribution may be used
instead of normal quantiles, but in our (first-order) asymptotic framework this does not make a difference. However, using 
\eqref{EqPhiRegressionStudentised} instead of \eqref{EqPhiRegression} gives a different bootstrap scheme and in case of the
null bootstrap this gives different asymptotic behaviour. 

Slight adaptations of Lemmas~\ref{LemmaSlopeHadamard} and~\ref{LemmaSlopeNull} show that the conditions of Theorem~\ref{thm:asymptotic_behavior_T} are again satisfied for the null bootstrap
$R_n =\PP_n \otimes \QQ_n$. Because the standardisation renders all limit variables
$\phi'_{P\otimes Q}(\GG_{P\otimes Q})$ and $\phi'_H(\GG_H)$ standard normal for $H$ in the null hypothesis,
condition \eqref{EqBootstrapLimit} is now satisfied without further conditions. Thus the null bootstrap is consistent under the minimal conditions of model \eqref{EqPhiRegression}.

This shows that studentising a test statistic may have beneficial effects on a bootstrap procedure, even at first order asymptotic level.
The benefits for higher-order asymptotic correctness are well studied (see e.g.\ \cite{HallBook}).

\subsection{Residual bootstraps} \label{sec:residual_bootstraps}
Besides bootstrap schemes (i) and (ii), we might use a bootstrap based on resampling residuals. For given estimators $\hat a_n$ and $\hat b_n$,
define the (estimated) residuals as $\hat \epsilon_i=Y_i-\hat a_n-\hat b_n X_i$, and let
$\HH_n^{res}=n^{-1}\sum_{i=1}^n\delta_{X_i,\hat\epsilon_i}$ be the empirical measure of the pairs
$(X_1,\hat\epsilon_1),\ldots,(X_n,\hat\epsilon_n)$, with marginal distributions $\PP_n$ and $\QQ_n^{res}$.
\begin{enumerate}[label=(\roman*)]
\setcounter{enumi}{2}
\item The residual bootstrap $R_n^{res}=\HH_n^{res,*}$ resamples pairs $(X_1^*,\epsilon_1^*),\ldots, (X_n^*,\epsilon_n^*)$ from
the pairs  $(X_1,\hat\epsilon_1),\ldots,(X_n,\hat\epsilon_n)$.
\item The fixed design residual bootstrap $R_n^{f,res}$ resamples $\epsilon_1^*,\ldots, \epsilon_n^*$ from
the pairs  $\hat\epsilon_1,\ldots,\hat\epsilon_n$, but sets $X_i^*=X_i$ for every $i=1,\ldots, n$.
\end{enumerate}
Given $(X_1^*,\epsilon_1^*),\ldots, (X_n^*,\epsilon_n^*)$, we can form bootstrap values 
$Y_i^*=\hat a_n+\hat b_n X_i^*+\epsilon_i^*$. Bootstrap (iii) is then identical to 
the bootstrap scheme (ii) as before, and hence the centred bootstrap statistic $\Tncent$ will work.

The residual bootstraps, in particular scheme (iv), are motivated by fixed design regression,
where the dependent variables $Y_1,\ldots, Y_n$ are not i.i.d.\, whence an empirical bootstrap involving these
values seems less natural. If the (true) residuals $\epsilon_1,\ldots, \epsilon_n$ are i.i.d.\, then the
estimated residuals $\hat\epsilon_1,\ldots,\hat\epsilon_n$ should be close to i.i.d.\ and resampling them is natural.
In the fixed design setting, the variables $X_i$ are obviously independent of the errors $\epsilon_i$, and the
centred bootstrap $\Tncent$ will work, for instance if $\max_{1\le i\le n}|X_i|=o(\sqrt n)$, $n^{-1}\sum_{i=1}^n(X_i-\bar X_n)^2\to \tau^2>0$
and given errors with a finite moment of order $>2$.

For a null bootstrap, it might be natural to construct bootstrap values 
as $Y_i^*=\hat a_n+\epsilon_i^*$. Combined with the least squares estimator $\hat a_n=\bar Y_n$ under the
null hypothesis rather than $\hat a_n$ as before, scheme (iii) would give yet another description of bootstrap scheme (ii).

\section{Goodness-of-fit for copulas}
\label{SectionCopulas}

Suppose we observe a random sample $(X_1,Y_1),\ldots,(X_n,Y_n)$ from the distribution of a two-dimensional random vector $(X,Y)$
with cumulative distribution function $H$ and univariate marginal distribution functions $F$ and $G$, which we assume to be continuous.
The \emph{copula function} corresponding to $H$ is the cumulative distribution function $C: [0,1]^2\to[0,1]$ of the pair $\bigl(F(X), G(Y)\bigr)$.
For a given parametric family $\{C_\theta: \theta\in\Theta\}$ of copula functions, we consider the testing problem
\begin{equation}
\label{eq:hypothesis_copulaGoF}
    \Hc_0: C \in \{C_\theta: \theta \in \Theta\} \quad \text{ versus } \quad \Hc_1: \:  C \notin \{ C_\theta : \theta\in\Theta \}.
\end{equation}
By Sklar's theorem, we have $H(x,y) = C\bigl( F(x), G(y)\bigr)$, for every $x,y\in \RR$. Thus we can put this in the framework of Section~\ref{sec:extension_framework_parameters}
with the map $\phi_1$ defined by
$$\phi_1(H,\theta)= H-C_\theta\circ(F,G).$$
Here $(F,G)$ is shorthand for the map $(x,y)\mapsto \bigl(F(x),G(y)\bigr)$.
Alternatively, as is often done in practice, we can first transform to uniform marginals and use the map
$$\phi(H,\theta)=H\circ(F^{-1},G^{-1})-C_\theta.$$
For $\HHn$ the empirical distribution function, the map $\phi(\HHn,\hat\theta_n)$ yields the \emph{empirical copula function} $\hat C_n:= \HHn\circ(\FF_n^{-1},\GG_n^{-1})$ minus the estimated parametric copula function $C_{\hat\theta_n}$. For the usual norms $\|\cdot\|_E$, for instance the uniform norm, the corresponding test statistics $T_{n,1}=\|\phi_1(\HHn,\hat\theta_n)\|_E$ and
$T_{n}=\|\phi(\HHn,\hat\theta_n)\|_E$ are nearly the same.
We shall restrict ourselves to the second statistic, given by the map $\phi$.

We shall see (again) that the empirical bootstrap $R_n=\HHn$ combined with the centred bootstrap test statistic gives good results. 

More involved, but perhaps more popular and natural, is a bootstrap under the null hypothesis. Because presently the null hypothesis specifies one of the distributions $C_\theta\circ (F,G)$, for varying $\theta$ and $(F,G)$, this is a \emph{semiparametric bootstrap} if the distributions $(F,G)$ are completely unspecified. However, by transforming the observations to uniform variables, this bootstrap scheme can actually be reduced to a parametric bootstrap, which resamples from the estimated copula $C_{\hat\theta_n}$, as follows.

The pairs of variables $(U_1, V_1), \dots, (U_n, V_n)$ defined by 
$U_i = F(X_i)$ and $V_i = G(X_i)$, are an i.i.d.\ sample from the copula distribution $C$.
Because $F$ and $G$ are assumed unknown, these variables are not observed, but
pretend for the moment that they were available. We could then form their empirical distribution $\HHnUV$ and estimators $\hat\theta_{n,U,V}$, and next the test statistic $\|\phi(\HHnUV,\hat\theta_{n,U,V})\|_E$ for the goodness-of-fit problem \eqref{eq:hypothesis_copulaGoF}, which specifies the parametric family $\{C_\theta: \theta\in\Theta\}$ for the distribution of the variables $(U_i,V_i)$. In this setting a bootstrap under the null hypothesis would consist of redrawing a sample of variables $(U_1^*, V_1^*), \dots, (U_n^*, V_n^*)$ from $C_{\hat\theta_{n,U,V}}$ and then forming the bootstrap test
values $\|\phi(\HH_{n,U,V}^*, \theta_{n,U,V}^*)\|_E$, for $\HH_{n,U,V}^*$ the empirical distribution of $(U_1^*, V_1^*), \dots, (U_n^*, V_n^*)$ and $\theta_{n,U,V}^*$ the
estimator of $\theta$ computed on the latter values. This procedure would entail an ordinary parametric bootstrap $R_n=C_{\hat\theta_n}$, in the sense of Section~\ref{sec:GoF_Setting}, and the results obtained there apply directly.

Parameter estimators $\hat\theta_n=\hat\theta_n(X_1,Y_1,\ldots, X_n,Y_n)$ in copula models are typically based on the ranks of the observations $(X_1, Y_1), \dots, (X_n, Y_n)$, and as such they take the same value if we used the unobservable variables $(U_1, V_1), \dots, (U_n, V_n)$ instead. In other words, we may assume that the estimators $\hat\theta_{n,U,V}$ in the preceding paragraph are equal to the original $\hat\theta_n$.

Under this condition the test statistic $\|\phi(\HHnUV,\hat\theta_{n,U,V})\|_E$ thus obtained from the unobservable variables $(U_i,V_i)$ is exactly the same as the test statistic $\|\phi(\HHn,\hat\theta_n)\|_E$ obtained from the original observations $(X_i,Y_i)$. In fact, the empirical copula $\hat C_{n,U,V}$, of the variables $(U_i, V_i)$ is identical
to the empirical copula $\hat C_n$ obtained from the observations $(X_i,Y_i)$. For $\FFnU$ and $\GGnV$ the marginal cumulative distribution functions of $U_1,\ldots,U_n$ and $V_1,\ldots, V_n$, respectively:
$$\hat C_n= \HHn\circ(\FF_n^{-1},\GG_n^{-1})= \HHnUV \circ (\FF_{n,U}^{-1}, \GG_{n,V}^{-1})=\hat C_{n,U,V}$$ 
(See \cite[First paragraph, page 538]{van2023weak}.)
It follows that the test statistic can be written in two equivalent ways, using either the observed data $(X_1, Y_1), \dots, (X_n, Y_n)$ or the unobserved data $(U_1, V_1), \dots, (U_n, V_n)$:
\begin{align*}
      T_n &= \| \phi(\HHn, \thetahatn) \|_E
    = \| \hat C_{n} - C_\thetahatn \|_E
    = \| \hat C_{n,U,V} - C_\thetahatn \|_E
    = \| \phi(\HHnUV, \thetahatn) \|_E.
\end{align*}
In practice we use the formula on the left, but theory may be based
on the formula on the right.
Given the estimate $\hat\theta_n$, we draw an i.i.d.\ sample $(U_1^*, V_1^*), \dots, (U_n^*, V_n^*)$ from the estimated copula $C_{\hat\theta_n}$, compute $\theta_{n}^*=\hat\theta_n(U_1^*,V_1^*,\ldots, U_n^*,V_n^*)$ and form the bootstrap value
\begin{align*}
    \Tneq = \| \phi(\HH_{n,U,V}^*, \theta_n^*) \|_E
    = \| \hat C_{n,U,V}^* - C_{\theta_n^*} \|_E
= \bigl\|\HH_{n,U,V}^* \circ (\FF_{n,U}^{*,-1}, \GG_{n,V}^{*,-1})- C_{\theta_n^*} \bigr\|_E.
\end{align*}
It may be noted that although by construction the $(U_i^*,V_i^*)$ possess uniform marginals, the bootstrap test statistic transforms these variables by the marginal quantiles functions before computing the distance of their empirical distribution to the estimated copula. This is necessary to mimic the formation of the test statistic, and follows naturally from our description through the map $\phi$.

The following examples give popular methods of estimating the copula parameter $\theta$. Both are based on the ranks of the observations.

\begin{example}[Pseudo-maximum likelihood estimator]
    If every copula $C_\theta$ admits a density $c_\theta$, then the pseudo-maximum likelihood estimator \cite{tsukahara2005semiparametric} is defined as
    \begin{align*}
        \thetahatn
        = \argmax_\theta
        \sum_{i = 1}^n \log c_\theta \big( \FFn(X_i), \GGn(Y_i) \big)
        = \argmax_\theta
        \sum_{i = 1}^n \log c_\theta \big( \FFnU(U_i), \GGnV(V_i) \big).
    \end{align*}
    In \cite{genest1995semiparametric} it is shown that this estimator is asymptotically normal. Using the literature on multivariate rank order statistics and in particular chapter 3 in \cite{ruymgaartshorackvanzwet1972}, it can be shown that under suitable regularity conditions the pseudo-maximum likelihood estimator is also asymptotically linear. This puts the pseudo-maximum likelihood estimator in the setting of Example \ref{ExampleAsymptoticLinearity} and hence the goodness-of-fit problem for copulas fits the general framework for parameter estimators as in Section \ref{sec:extension_framework_parameters}.
\end{example}

\begin{example}[Inversion of Kendall's tau]
    For many one-dimensional copula families, the map 
    $\theta\mapsto \psi(\theta) := \tau(C_\theta) := 4 \iint C_\theta\, dC_\theta - 1$ is one-to-one between $\Theta$ and $[-1, 1]$. The quantity $\tau=\tau(C)$ can be estimated by Kendall's tau statistic $\hat\tau_n = (\#\text{Concordant Pairs} - \#\text{Discordant Pairs}) / (n (n-1))$.
    Next, $\theta$ is estimated by inversion, i.e. $\thetahatn = \psi^{-1} (\hat\tau_n)$. 
    It can be shown that $\hat\tau_n +1$ is a U-statistic (see Example 12.5 in \cite{van2000asymptotic}).
    U-statistics can also be shown to be asymptotically linear (Theorem 12.3 in \cite{van2000asymptotic}) and we obtain that $\psi(\thetahatn) +1$ is also asymptotically linear. If we assume $\psi$ is also continuous, then by means of the continuous mapping theorem the limiting distribution can be found and the setting of Example \ref{ExampleAsymptoticLinearity} is recovered.
    
\end{example}

The function $\phi$ is Hadamard differentiable (combining Lemma 3.10.30 in \cite{van2023weak} and the regularity of the copula model $\theta \mapsto C_\theta$).
Therefore, we can obtain similar results as the ones from Section \ref{sec:GoF_Setting}. In particular, we can apply Lemma \ref{LemmaGoodnessNull} to show that the parametric bootstrap also works for copulas.

\section{Simulation study}
\label{SectionSimulation}

We illustrate the developed theory in a simulation study of the hypothesis tests for independence, the slope in the linear regression setting and the goodness-of-fit setting. The theory gives correction terms for the bootstrap test statistic that depend on the bootstrap scheme used to perform the resampling. We determine the power and rejection rates of different combinations of bootstrap resampling schemes and (corrected) bootstrap test statistics, for varying sample sizes. The goal is to illustrate which combinations work and which do not work, and compare this to the developed theory. In Appendix \ref{app:comparison} we compare the performance of the theoretically valid bootstrap schemes.

\medskip

We have implemented the bootstrap-based hypothesis testing procedures as an \texttt{R} package called BootstrapTests \cite{BootstrapTestPackage}, in the settings of independence testing, testing the slope in a linear regression setting, and goodness-of-fit testing.

\medskip

The p-values are computed using 100 bootstrap samples.
For each setting, we approximate the power and the level of the test using 200 simulations.
The simulations were run on the DelftBlue supercomputer \cite{DHPC2024}. To give an impression of the runtime: around 30 CPU-hours were needed for the independence testing simulations.

\subsection{Simulation: independence testing}

We use the set-up as in Example \ref{example_KS}, with $\Xc = \Yc = \RR$ and $\Fc, \Gc$ consisting of the indicator functions of the cells $(-\infty, a],$ for $a$ varying over $\RR$.
For the bootstrap resampling scheme $R_n$, we either use the empirical bootstrap $R_n= \HH_n$ or the independence bootstrap $R_n = \PP_n \otimes \QQ_n$ (which we also call the null bootstrap). Then, we choose the bootstrap test statistic $T_n^*$. We have two options here: the equivalent bootstrap test statistic $\Tneq$, or the centered bootstrap test statistic $\Tncent$. As summarized in Table \ref{tab:indepdence_pairs}, the theory shows that consistent combinations of resampling scheme $R_n$ and bootstrap test statistics $T_n^*$ are given by $(\HHn, \Tncent)$ and $(\PP_n \otimes \QQ_n, \Tneq)$. 

\medskip

In our simulation study, we consider a variety of data-generating processes. We create a sample $X_1,\dots, X_n \simiid \Nc (0,1) $ and construct $Y_i = b X_i + \epsilon_i$, for $i=1,\dots,n$ and $\epsilon_i \simiid \Nc(0,1).$ Together, we create pairs $(X_i,Y_i)$ for $i=1,\dots,n$ for which we want to test independence.
Remark that $b = 0$ corresponds to the null hypothesis of independence between $X$ and $Y$. The value $b = 0$ is used in the empirical level analysis. Higher values of $b$ make the data-generating process further away from the null hypothesis of independence and for this reason we varied the $b$-values in our power analysis simulations.

\medskip

Combining the results of Figures~\ref{fig:power_test_independence} and~\ref{fig:level_test_independence}, we observe, as predicted by the theory, that only two combinations of resampling scheme and bootstrap test statistic show a high power and are relatively well-calibrated in terms of the level.
Furthermore, we observe that as the value of $b$ increases, the power also rises. However, for very small values of $b$, the power remains low for all sample sizes that we tested. This is in line with the intuition that $b$ represents the perturbation from independence.
The two other combinations have zero power and a level of zero, independent of the sample size, meaning that these combinations are forming inconsistent and invalid tests.


\begin{figure}[p]
  \centering
  \resizebox{0.75\textwidth}{!}{\input{independence_power_plot}}
  \caption{
    Power of the independence test as a function of the sample size, 
    for different values of $b$,
    different resampling schemes (empirical or independence)
    and different bootstrap test statistics
    ($\Tncent$ and $\Tneq$).
  }
  \label{fig:power_test_independence}


  \centering
  \resizebox{0.7\textwidth}{!}{\input{independence_level_plot}}
  \caption{
    Level of the independence test as a function of the sample size, 
    for different resampling schemes 
    and different bootstrap test statistics
    ($\Tncent$ and $\Tneq$).
  }
  \label{fig:level_test_independence}
\end{figure}

\subsection{Testing the slope in linear regression}

We want to illustrate the results from Section \ref{SectionLinearRegression}. To this end, we perform bootstraps on i.i.d samples $(X_1,Y_1),\dots, (X_n,Y_n)$, following the model
\begin{equation*}
    Y_i = b X_i + \epsilon_i, \quad i=1,\dots,n,
\end{equation*}
with $X_i \simiid \Nc(0,1)$, $\epsilon_i \simiid \Nc(0,1)$ and $X_i$ independent from $\epsilon_i$. Observe that the minimal conditions in the regression model \eqref{EqRegressionModel} are satisfied under independence of $\epsilon_i$ and $X_i$. Under this stronger assumption of independence, the empirical and null bootstrap are both yielding valid tests, when paired with the correct bootstrap test statistic. Here, the null and empirical bootstrap correspond to bootstrap schemes (i) and (ii) in Section \ref{SectionLinearRegression}, respectively. From the theory (similar to the indepence testing) we expect that $(\HHn, \Tncent)$ and $(\PP_n \otimes \QQ_n, \Tneq)$ should work. We also consider the residual and hybrid null bootstrap. The residual bootstrap corresponds to bootstrap scheme (iii) in Section \ref{sec:residual_bootstraps} and described there, it will yield the same results as the empirical bootstrap scheme. The hybrid null bootstrap performs the same resampling procedure as the residual bootstrap, but forms $Y_i^* =  \hat{a}_n + \epsilon_i^*$, with $\hat{a}_n$ the least-squares estimator for the intercept $a.$ So here we actually estimate the model $Y_i = a +b X_i + \epsilon_i$, with the true and unknown $a=0$ in our simulation setup. If $\hat{a}_n$ is reasonably close to the true $a$, then the hybrid null bootstrap should perform similarly to the independence bootstrap, while, strictly speaking, not performing the same bootstrap procedure.

\medskip

In our simulation study, we consider a variety of data-generating processes, similar to the independence testing simulation. We create a sample $X_1,\dots, X_n \simiid \Nc (0,1) $ and construct $Y_i = b X_i + \epsilon_i$, with $\epsilon_i \simiid \Nc(0,1),$ for $i=1, \dots, n$. To distinguish between the null and alternative hypothesis, we vary the value for $b$. Higher values of $b$ make the data-generating process further away from the null hypothesis.

\medskip

\begin{figure}[hp]
  \centering
  \resizebox{0.8\textwidth}{!}{\input{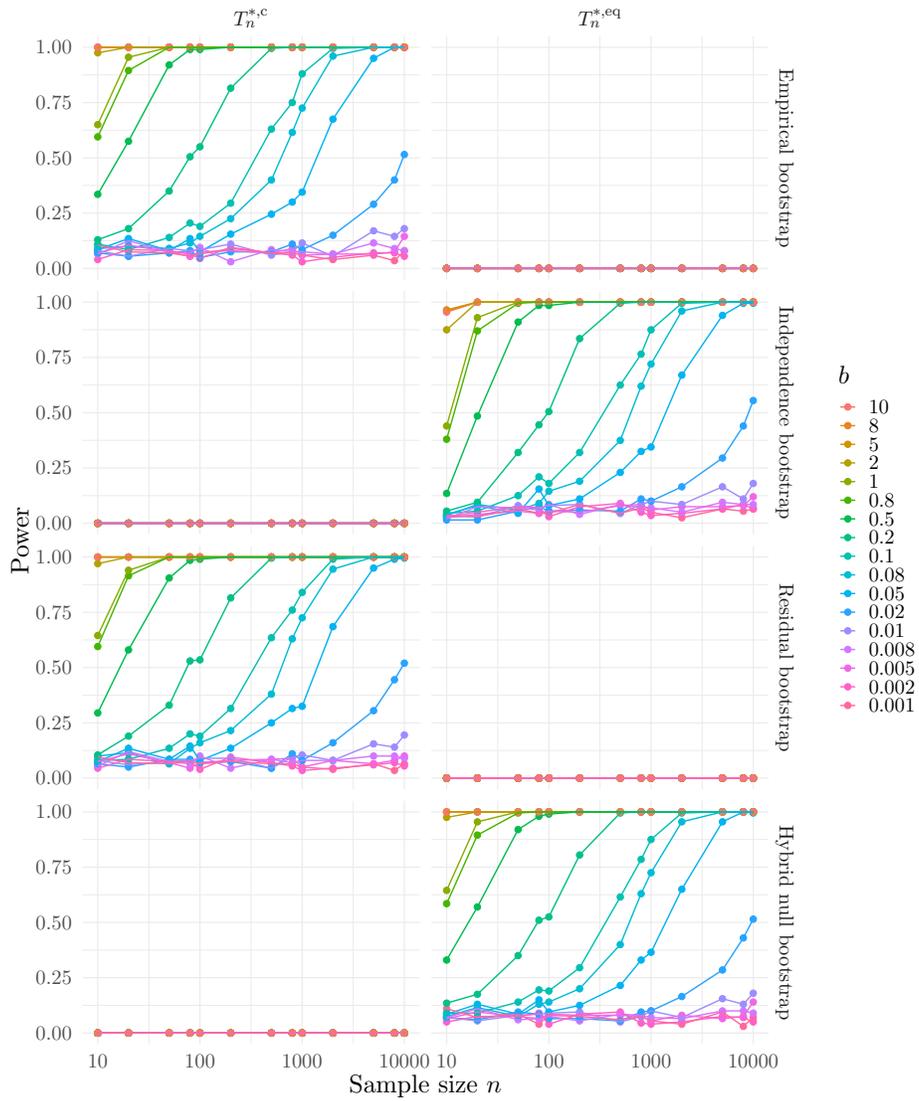}}
  \caption{
    Power as a function of the sample size in the regression setting, for different value of the coefficient $b$ and for different combinations of bootstrap resampling schemes $R_n$ and bootstrap test statistic ($\Tncent$ or $\Tneq$).
  }
  \label{fig:power_test_regressiontest}
\end{figure}

\begin{figure}[hp]
    \centering
    \resizebox{0.8\textwidth}{!}{\input{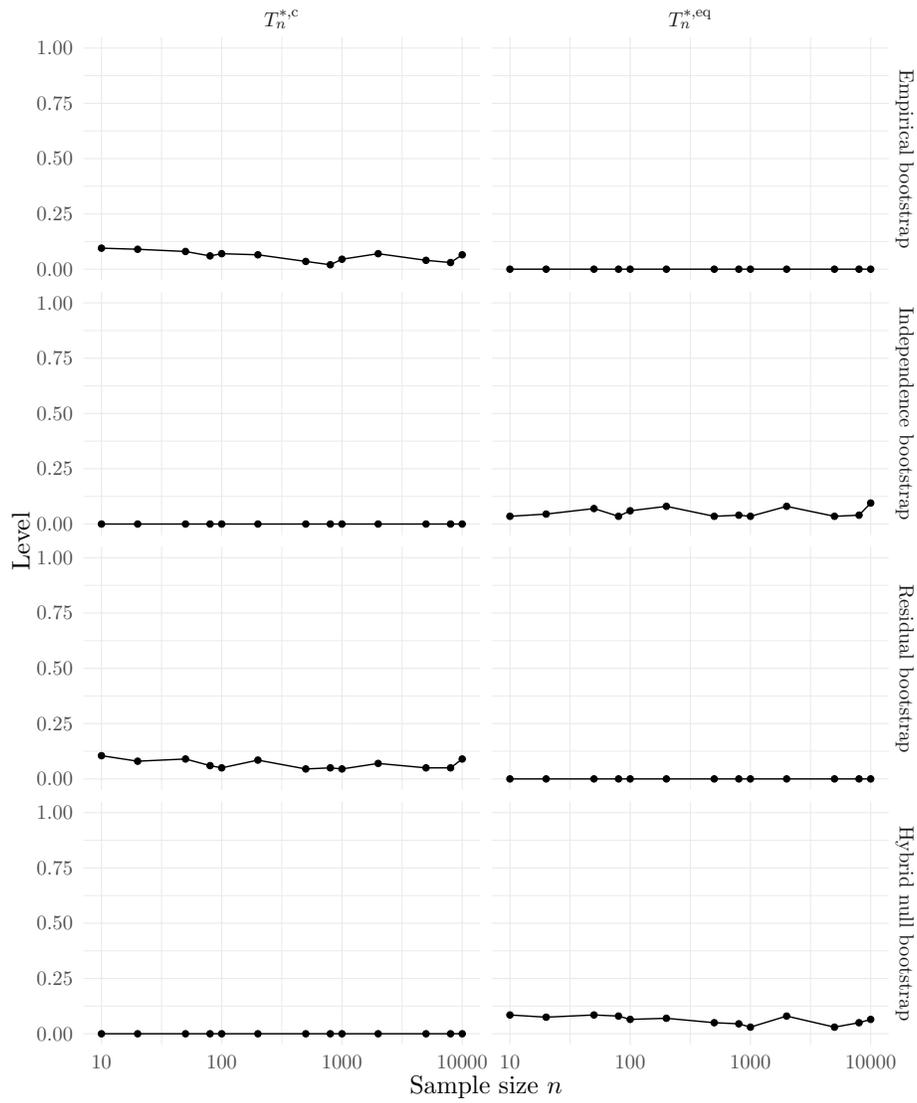}}
    \caption{
    Level as a function of the sample size in the regression test setting, for different combinations of bootstrap resampling schemes $R_n$ and bootstrap test statistic ($\Tncent$ or $\Tneq$).
    %
    }
    \label{fig:level_test_regressiontest}
\end{figure}

Combining the results of Figures~\ref{fig:power_test_regressiontest} and~\ref{fig:level_test_regressiontest}, we observe, as predicted by the theory, that only four combinations of resampling scheme and bootstrap test statistic show a high power and are relatively well-calibrated in terms of the level. The same combinations of resampling procedure and test statistic as in the independence testing setting also work here. The residual bootstrap performs similar to the empirical bootstrap (as it performs the same bootstrap procedure), but we include it for completeness. The hybrid null bootstrap performs similar to the independence bootstrap with the same combination of bootstrap test statistic.
Furthermore, we observe that as the values of the perturbation from independence increases, the power also rises. However, for very small values of the slope $b$, the power remains low for all sample sizes that we tested.
The other combinations have zero power and a level of zero, independent of the sample size, meaning that these combinations are forming inconsistent and invalid tests.

For completeness, we have also added the simulation results for the fixed design residual bootstrap resampling schemes in Appendix \ref{app:Fixed design regression}.


\subsection{Goodness-of-fit testing}\label{sec:simulation_GOF_general}

We want to illustrate the results from Section \ref{sec:GoF_Setting}. We consider samples $X_1, \dots, X_n \simiid H$. Under the null hypothesis, this distribution $H$ comes from a given parametrised family  of distributions $\{H_\theta:\theta\in\Theta\}$. In our simulation study we test whether or not our data-generating distribution comes from a family of normal distributions. In particular, for our simulation setup we use the setting of Example \ref{example:kolmogorov-smirnov GOF} with observations in $\RR$. For the estimators $\thetahatn$, we use the minimum distance estimators. 

\medskip

In the goodness-of-fit setting we use two types of bootstrap resampling schemes: the empirical $R_n= \HH_n$ and null bootstrap $R_n = H_{\thetahatn}$. According to our theory, the empirical and null bootstrap are both yielding valid tests when paired with the correct bootstrap test statistic. In particular, similar to the independence testing, we expect that $(\HHn, \Tncent)$ and $(H_{\thetahatn}, \Tneq)$ should work, with the precise forms of $\Tnst$ given in Section \ref{sec:GoF_Setting}. To calculate $\Tnst$, it is necessary to calculate a corresponding $\theta_n^*$. The choice of $\theta_n^*$ is of great importance for the success of the bootstrap procedure. The simulations also show this. In our minimum distance setting we choose $\theta_n^*$ to be $\thetastMD$ from Equation \eqref{EqBootstrapMDEstimator}. We illustrate that the choice of the `centered' $\thetastMD$ is correct, where we also show that the non-centered gives incorrect results. Note that it is only necessary to incorporate a centering correction term $R_n - H_{\hat\theta_n}$ $R_n$ in $\thetastMD$ whenever $R_n$ is not equal to the null bootstrap.

\medskip

To simulate samples under the null, we simply generate a random sample $X_1,\dots, X_n$ from the standard normal distribution. To simulate data under the alternative, we use the $t$-distribution, the log-normal distribution, mixtures of normal distributions, and Cauchy distributions. For specific parameters, see Table \ref{tab:distributions_GOF_simulation} in Appendix \ref{app:GOF}. The idea is that these distributions are increasingly `distant' to the family of normal distributions, to make the goodness-of-fit testing problem progressively harder. This is analogous to increasing the value of the slope in the slope regression test. 

\medskip

Combining the results of Figures~\ref{fig:power_test_MDestimator_GOFtest},~\ref{fig:power_test_MDcentestimator_GOFtest}, we observe, as predicted by the theory, that only four combinations of resampling scheme and bootstrap test statistic show a high power and are relatively well-calibrated in terms of the level. The same combinations of resampling procedure and test statistic as in the independence/regression testing setting also work here, but the choice of $\thetast$ is an additional consideration. In particular, if the bootstrap resampling scheme is different from the null bootstrap, the chosen $\thetast$ needs to be adjusted according to theory. For instance, in Figure \ref{fig:power_test_MDestimator_GOFtest}, where the `wrong' $\thetastMD$ is chosen for the empirical bootstrap, it is observed that the power of the tests shows worse performance compared to the case in Figure \ref{fig:power_test_MDcentestimator_GOFtest}, where the correct $\thetastMD$ is chosen.
What we also observe that as the distributions become increasingly `distant' from the normal family of distributions, the power also rises. However, for distributions that are very similar to the normal (e.g. the $t$-distribution with high degrees of freedom), the power remains low for all sample sizes that we tested.

\begin{figure}[htb]
    \centering
    \resizebox{1.0\textwidth}{!}{\input{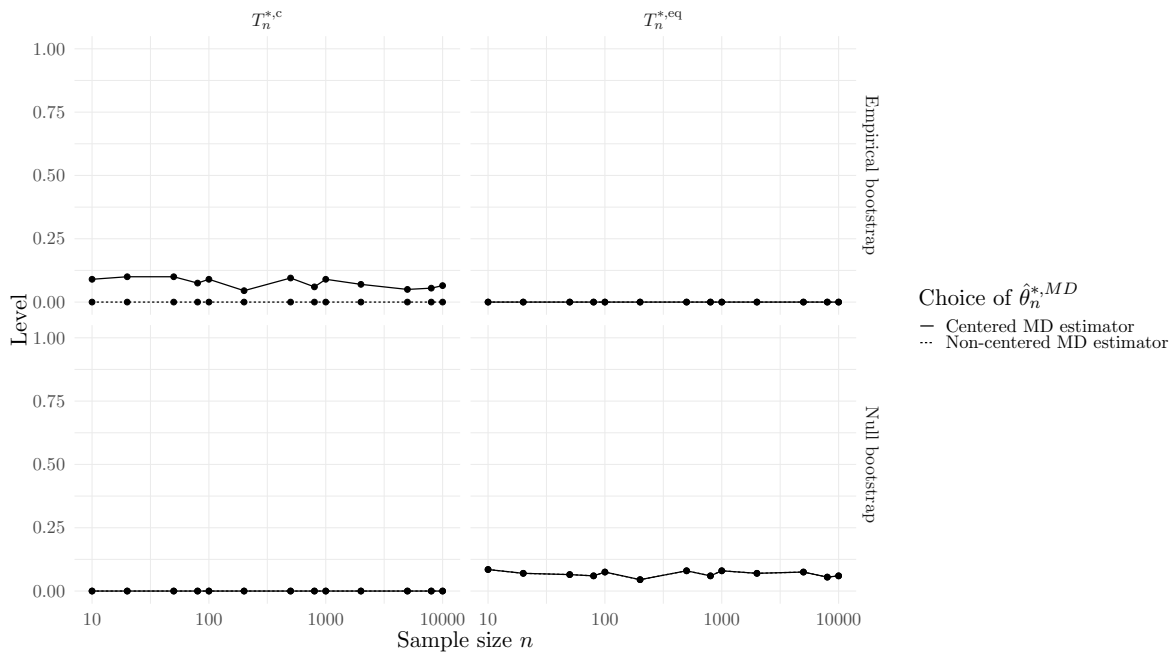}}
    \caption{
    Level as a function of the sample size in the goodness-of-fit test setting,
    for different combinations of bootstrap resampling schemes $R_n$, bootstrap test statistic ($\Tncent$ or $\Tneq$), and bootstrap estimators.
    The solid line refers to the `centered' $\theta_n^{*,MD}=\argmin_\theta \| \HHnst - H_{\theta}- R_n+H_{\hat\theta_n} \|$ and the dashed line is the `non-centered' $\theta_n^{*,MD}=\argmin_\theta \| \HHnst - H_{\theta} \|$. 
    }
    \label{fig:level_test_GOFtest}
\end{figure}

From Figure ~\ref{fig:level_test_GOFtest} we observe an empirical level of zero, independent of the sample size, for the incorrectly chosen combinations, meaning that these combinations form inconsistent and invalid tests. For the combinations chosen according to theory, the tests are relatively well-calibrated and show a level approximately equal to the true level.

\begin{figure}[hp]
    \centering
    \resizebox{0.95\textwidth}{!}{\input{GOF_power_MDestimator_without_Cauchy_plot}}
    \caption{
    Power of the goodness-of-fit test as a function of the sample size, for different data generating processes and different combinations of bootstrap resampling schemes $R_n$ and bootstrap test statistic ($\Tncent$ or $\Tneq$).
    Here, the bootstrap version estimator is the `non-centered' $\theta_n^{*,MD}=\argmin_\theta \| \HHnst - H_{\theta} \|$.
    }
    \label{fig:power_test_MDestimator_GOFtest}

  \centering
  \resizebox{0.90\textwidth}{!}{\input{GOF_power_MDcentestimator_without_Cauchy_plot}}
  \caption{
    Power of the goodness-of-fit test as a function of the sample size, for different data generating processes and different combinations of bootstrap resampling schemes $R_n$ and bootstrap test statistic ($\Tncent$ or $\Tneq$).
    Here, the bootstrap version estimator is the `centered'  $\theta_n^{*,MD}=\argmin_\theta \| \HHnst - H_{\theta}- R_n+H_{\hat\theta_n} \|.$
  }
  \label{fig:power_test_MDcentestimator_GOFtest}
\end{figure}

\bibliographystyle{abbrv}
\bibliography{main}{}

\newpage

\appendix

\appendixpage           
\addappheadtotoc        

\section{Minimum distance estimators}
\label{appendix:minimum-distance-estimators}

This section gives formal statements and proofs for the informal claims made in Section~\ref{SectionMinimumDistanceInformal}.
Let $\{H_\theta: \theta\in\Theta\}$ be a set of probability measures indexed by an open subset $\Theta\subset \RR^d$,
viewed as elements of a normed space $(E,\|\cdot\|)$. Let 
$\HHn$, $R_n$ and $\HHnst$ be random maps in $E$ and set
\begin{align*}
\hat{\theta}_n&=\argmin_{\theta\in\Theta} \|\HHn - H_\theta\|,\\
\theta_n^{*}&=\argmin_{\theta\in\Theta} \| \HH_n^*- H_\theta  - R_n + H_{\hat{\theta}_n} \|,\\
\theta_n^{*,u}&=\argmin_{\theta\in\Theta} \| \HH_n^*- H_\theta \|.
\end{align*}
All variables $\HH_n$, $R_n$, $\HHnst$, $\GG_n$,  $\MM_n$, $\MM_n^*$, $\hat\theta_n$, $\theta_n^*$ and $\theta_n^{*,u}$ in the following are maps from a common probability space into the range spaces as indicated,
possibly non-measurable. The limit variables $\GG_0$, $\MM$, $\hat h_0$ are assumed to be tight Borel measurable maps on a common probability space. We think of the starred maps as bootstrap versions, and consider their distributional convergence 
in a conditional setup, denoted as $| X_{1.n}\cvweakly$, and formally understood in the sense of Section~\ref{AppendixMeasurability}. In the intended application the variables $\HH_n$ and $\HH_n^*$ are the empirical measures of a sample of observations $X_{1.n}=(X_1,\ldots, X_n)$ and a bootstrap sample from $R_n=R_n(X_1,\ldots,X_n)$.

\begin{lemma}\label{lemma_argmin}
Suppose that $\inf_{\theta: \|\theta-\theta_0\|>\delta} \|H_\theta - H_{\theta_0}\| > 0 $, for all $\delta>0$,  and suppose that 
$\theta\mapsto H_\theta$ is continuously Fr\'echet differentiable in a neighbourhood of $\theta_0$ with derivative 
$\dot H_{\theta_0}: \RR^d\to E$ of rank $d$. 
If $\theta_{0,n}\to\theta_0$ and $\sqrt n\|\HHn - H_{\theta_{0,n}}\|=O_P(1)$, then  $\sqrt n\|\hat\theta_n-\theta_{0,n}\|=O_P(1)$.
\end{lemma}

\begin{proof}
The continuous differentiability of $\theta\mapsto H_\theta$ gives that \\ 
$H_\theta-H_{\theta_{0,n}}=\int_0^1\dot H_{\theta_{0,n}+s(\theta-\theta_{0,n})} (\theta-\theta_{0,n})\,ds$,
for $\theta$ and $\theta_{0,n}$ in a sufficiently small neighbourhood of $\theta_0$. By making the
neighbourhood smaller, if necessary, it can be ensured that $\|\dot H_{\theta_{0,n}+s(\theta-\theta_{0,n})} -\dot H_{\theta_0}\|<\epsilon$,
for arbitrarily small $\epsilon >0$, in view of the continuity of the derivative. 
It follows that $\|H_\theta-H_{\theta_{0,n}}-\dot H_{\theta_0}(\theta-\theta_{0,n})\|\le \epsilon \|\theta-\theta_{0,n}\|$. 
Because $\dot H_{\theta_0}$ is of full rank, it follows that there exist positive constants $\delta,  C$  such that
$\|H_\theta - H_{\theta_{0,n}}\|\ge C\|\theta - \theta_{0,n} \|$, whenever $\|\theta-\theta_0\|\le\delta$ and $\|\theta_{0,n}-\theta_0\|\le\delta$.

By the triangle inequality $\|H_{\hat{\theta}_n} - H_{\theta_{0,n}}\|\leq \|H_{\hat{\theta}_n} - \HHn\| + \|\HHn - H_{\theta_{0,n}}\|
\leq 2 \| \HH_n - H_{\theta_{0,n}}\|$, by the definition of $\hat\theta_n$. Because the right side tends to zero in probability and
$\|H_{\theta_{0,n}} - H_{\theta_0}\|\to 0$, we see that  $\|H_{\hat{\theta}_n} - H_{\theta_0}\|\to 0$ in probability. Since $\|H_\theta - H_{\theta_0}\|$ is bounded away 
from zero on the set $\{\theta: \|\theta-\theta_0\|>\delta\}$, we see that the probability of the event
$\|\hat\theta_n-\theta_0\|\le\delta$ tends to one. On this event we have $C\|\hat\theta_n-\theta_{0,n}\|
\le \|H_{\hat{\theta}_n} - H_{\theta_{0,n}}\|\le 2 \| \HH_n - H_{\theta_{0,n}}\|=O_P(n^{-1/2})$.
\end{proof}

\begin{lemma}
\label{lemma:joint_cv_H_theta}
Assume that $\{H_\theta : \theta\in\Theta\}$ satisfies the conditions of Lemma~\ref{lemma_argmin} and that
$\sqrt n(\HH_n-H_{\theta_{0,n}})\cvweakly \GG_0$ in $E$, for a tight variable $\GG_0 $
such that the stochastic process $h\mapsto \| \GG_0   - \dot{H}_{\theta_0}h \|$ possesses a unique point of
minimum $\hat h_0$. Then 
$\sqrt{n}(\HHn - H_{\theta_{0,n}} ,\hat\theta_n - \theta_{0,n}) \cvweakly (\GG_0  , \hat h_0)$,   in $E\times \RR^d$.
\end{lemma}

\begin{proof}
\label{proof:lemma:joint_cv_H_theta}
The lemma is a consequence of Proposition~\ref{theorem:joint_argmax_continuous_mapping}, applied with
$\hathn := \sqrt{n}(\hat\theta_n - \theta_{0,n})$, $\GG_n := \sqrt{n}(\HHn - H_{\theta_{0,n}})$ and the 
 stochastic processes $\{\MM_n(h): h\in\RR^d\}$ and $\{\MM(h): h\in \RR^d\}$ given by
\begin{align*}
    \MM_n(h)
    &= -\| \GG_n   - \sqrt{n}( H_{\theta_{0,n} + h/\sqrt{n}}    - H_{\theta_{0,n}}) \|,     \\
    \MM(h)
    &= -\| \GG_0  - \dot{H}_{\theta_0}h \|.
\end{align*}
The sequence $\hathn := \sqrt{n}(\hat\theta_n - \theta_{0,n})$ is bounded in probability by Lemma~\ref{lemma_argmin}.
For every compact set $K\subset \RR^d$, the processes $(\GG_n,\MM_n)$ converge in distribution in $E\times \ell^\infty(K)$ 
to the process $(\GG_0 , \MM)$, by the continuous mapping theorem and the differentiability of the map $\theta \mapsto H_\theta$.
The process $\MM$ is continuous in $h$. Thus the conditions of Proposition~\ref{theorem:joint_argmax_continuous_mapping}
are satisfied.
\end{proof}

\begin{lemma}\label{lemma_convergence_theta_rescaled}
Assume that $\{H_\theta: \theta\in\Theta\}$ satisfies the conditions of Lemma~\ref{lemma_argmin} and that
$\sqrt n(\HH_n^*-R_n) | X_{1.n}\cvweakly \GG_0 $ in $E$, outer almost surely or in outer probability,
for a tight variable $\GG_0 $
such that the stochastic process $h\mapsto \| \GG_0   - \dot{H}_{\theta_0}h \|$ possesses a unique point of
minimum $\hat h_0$,
and that $\hat{\theta}_n \to\theta_0$ outer almost surely or in outer probability. Then
\begin{equation*}
        \sqrt{n}(\HH_n^*-R_n,\theta_n^* - \hat\theta_n ) \, |\,  X_{1.n} \cvweakly (\GG_0 , \hat h_0),\qquad\text{outer a.s. or in prob}.
\end{equation*}
\end{lemma}

\begin{proof}
By the triangle inequality followed by the definition of $\theta_n^*$,
$$\|H_{\theta_n^*}-H_{\hat\theta_n}\|
\le\| \HH_n^* - R_n - H_{\theta_n^*} + H_{\hat{\theta}_n} \|+\|\HH_n^* - R_n\|
\le 2\|\HH_n^* - R_n\|.$$
By assumption the right side tends to zero in outer probability conditionally given $X_{1.n}$ and hence also unconditionally. Since $\hat\theta_n\to \theta_0$ by assumption,
it follows that also $\|H_{\theta_n^*}- H_{\theta_0}\|$ tends to zero in outer probability. Because  by assumption  $\inf_{\|\theta-\theta_0\|>\delta}\|H_\theta-H_{\theta_0}\|>0$ for every $\delta>0$, we can conclude that $\theta_n^*$ tends to $\theta_0$ in outer probability.

It is seen in the proof of Lemma~\ref{lemma_argmin} that there exist positive constants $\delta,  C$  such that
$\|H_\theta - H_{\theta'}\|\ge C\|\theta - \theta' \|$, whenever $\|\theta-\theta_0\|\le\delta$ and $\|\theta'-\theta_0\|\le\delta$.
Thus on the event where $\|\theta_n^*-\theta_0\|<\delta$ and
$\|\hat\theta_n-\theta_0\|<\delta$, we have 
$C\|\theta_n^*-\hat\theta_n\|\le \|H_{\theta_n^*}-H_{\hat\theta_n}\|\le 
2\| \HH_n^* - R_n\|$. It follows that 
$h_n^*:=\sqrt n(\theta_n^*-\hat\theta_n)$ is bounded in outer probability. If $\hat{\theta}_n \to\theta_0$ outer almost surely and $\sqrt n(\HH_n^*-R_n) | X_{1.n}\cvweakly \GG_0 $ in $E$, outer almost surely, then the conclusion can be strengthened with the same argument to 
$\Pr(\|h_n^*\|>M_n|  X_{1.n}\bigr)\to 0$ outer almost surely, for every $M_n\to\infty$.

Under suitable measurability conditions, the lemma is a consequence of the conditional version of  Proposition~\ref{theorem:joint_argmax_continuous_mapping}, applied with
$\GG_n^* := \sqrt{n}(\HHnst - R_n)$ and the 
 stochastic processes $\{\MM_n(h): h\in\RR^d\}$ and $\{\MM(h): h\in \RR^d\}$ given by
\begin{align*}
    \MM_n^*(h)
    &= -\| \GG_n^*   - \sqrt{n}( H_{{\hat\theta_n} + h/\sqrt{n}}- H_{\hat\theta_n}) \|,\\
    \MM(h)
    &= -\| \GG_0  - \dot{H}_{\theta_0}h \|.
\end{align*}
We have $(\GG_n^*,\MM_n)| X_{1.n}\cvweakly (\GG_0 ,\MM)$ in $E\times\ell^\infty(K)$, for every compact $K\subset \RR^d$, $h_n^*$ is conditionally tight, and the limit process possesses a unique point of maximum by assumption. Thus in the almost sure case, Proposition~\ref{theorem:joint_argmax_continuous_mapping} applies given almost every sequence $X_{1.n}$. The case of convergence in probability can be reduced to the  almost sure case by the characterisation of convergence in  probability as almost sure convergence along subsequences.

To avoid unnecessary measurability conditions, we sketch a full argument. As in Section~\ref{AppendixMeasurability}, write $\Expec_Z^\outer f(\GG_n^*,h_n^*)$ for the conditional expectation of the joint measurable cover function $f(\GG_n^*,h_n^*)^\outer$ given $X_{1.n}$, with $Z$ referring to the randomness in the bootstrap samples and $\Expec_Z$ the expectation on $Z$ given $X_{1.n}$.

For every bounded Lipschitz function $f: E\to[0,1]$, closed set $F\subset \RR^d$, compact set $K\subset\RR^d$ and Lipschitz function $\chi: \RR\to[0,1]$ with $1_{[0,\infty)}\le \chi$, we have 
\begin{align*}
\Expec_Z^\outer f(\GG_n^*)1_{h_n^*\in F}&\le 
\Expec_Z^\outer f(\GG_n^*)1\bigl\{\sup_{h\in F\cap K} \MM_n^*(h)\ge \sup_{h\in K} \MM_n^*(h)\bigr\}+ \Prob_Z^\outer(h_n^*\notin K)\\
&\le \Expec_Z^\outer f(\GG_n^*)\chi\Bigl(\sup_{h\in F\cap K} \MM_n^*(h)-\sup_{h\in K} \MM^*_n(h)\Bigr)+ \Prob_Z^\outer(h_n^*\notin K).
\end{align*}
The second term on the far right side can be made arbitrarily small by choice of $K$, while the first term converges to $\Expec f(\GG_0)\chi\bigl(\sup_{h\in F\cap K}\MM(h)- \sup_{h\in K} \MM(h)\bigr)$,  almost surely or in  probability, since the variable in the expectation is a Lipschitz function of $(\GG_n^*,\MM_n^*)$.
By the argument in the proof of Proposition~\ref{theorem:joint_argmax_continuous_mapping}, for a sequence of Lipschitz functions $\chi_m\downarrow 1_{[0,\infty)}$, the limiting expectations decrease to $$\Expec f(\GG_0)1_{[0,\infty)}\bigl(\sup_{h\in F\cap K}\MM(h)- \sup_{h\in K} \MM(h)\bigr)\le \Expec f(\GG_0)1_{\hat h_0\in F}+\Prob(\hat h_0\notin  K).$$ We conclude that for every $\epsilon>0$, there exist random variables 
$Y_n(\epsilon, F)\ge \Expec_Z^\outer f(\GG_n^*)1_{h_n^*\in F}$ such that 
$Y_n(\epsilon,F)\to y(\epsilon,F)$,  almost surely or in  probability, for a number $y(\epsilon,F)$ with $y(\epsilon,F)\le \Expec f(\GG_0)1_{\hat h_0\in F}+\epsilon$. 

For a given bounded Lipschitz function $g: \RR^d\to[0,1]$, the functions $g_m=\sum_{i=1}^m m^{-1}1_{F_{i,m}}$, for $F_{i,m}=\{x: g(x)\ge (i-1)/m\}$, satisfy $0\le g\le g_m\le 1$ and $|g-g_m|\le 1/m$. We have 
$\Expec_Z^\outer f(\GG_n^*)g(h_n^*)\le \sum_i m^{-1}Y_n(\epsilon,F_{i,m})
\to \sum_im^{-1} y(\epsilon,F_{i,m})=:y(\epsilon, g_m)$,  almost surely or in  probability, where $y(\epsilon, g_m)\le \Expec f(\GG_0)g_m(\hat h_0)+\epsilon$. As $m\to\infty$, the latter expression tends to $\Expec f(\GG_0)g(\hat h_0)+\epsilon$. We conclude that for every $\epsilon>0$, there exist variables $Y_n(\epsilon,f, g)\ge\Expec_Z^\outer f(\GG_n^*)g(h_n^*)$ such that $Y_n(\epsilon,f, g)\to y(\epsilon,f,g)$  almost surely or in   probability for a number $y(\epsilon,f,g)$ with
$y(\epsilon,f,g)\le \Expec f(\GG_0)g(\hat h_0)+\epsilon$.

Since $\Expec_Z^\outer f(\GG_n^*)(1-g)(h_n^*)\ge \Expec_Z^\outer f(\GG_n^*)-\Expec_Z^\outer f(\GG_n^*)g(h_n^*)$, applying the preceding with the functions $f$ and $1-g$ we find that $\Expec_Z^\outer f(\GG_n^*)g(h_n^*)\ge\Expec_Z^\outer f(\GG_n^*)-Y_n(\epsilon,f,1-g)\to\Expec f(\GG_0)-y(\epsilon,f,1-g)\ge \Expec f(\GG_0)g(\hat h_0)-\epsilon$, where the convergence is  almost surely or in  probability. Thus for every $\epsilon>0$, the sequence $\Expec_Z^\outer f(\GG_n^*)g(h_n^*)$ is sandwiched between two sequences of random variables that converge  almost surely or in  probability to a number between $\Expec f(\GG_0)g(\hat h_0)-\epsilon$ and $\Expec f(\GG_0)g(\hat h_0)+\epsilon$. 

This implies that
 $\Expec_Z^\outer f(\GG_n^*)g(h_n^*)\to \Expec f(\GG_0)g(\hat h_0)$, outer almost surely or in outer probability, for every pair of bounded Lipschitz functions $f: D\to\RR$ and $g: \RR^d\to\RR$. We conclude by applying Lemma~\ref{LemmaProductFunctionsConditional}.
\end{proof}

The following result is analogous to Lemma~\ref{lemma:joint_cv_H_theta}, but it is applicable to observations from a distribution $H$ that does not necessarily belong to the model $\{H_\theta: \theta\in\Theta\}$. It is restricted to Hilbertian distances. We write $\theta\to\partial\Theta$ for $\theta$ approximating the boundary of $\Theta$ in the one-point compactification of $\Theta$.

\begin{lemma}[Hilbert space]
\label{LemmaMDHilbertian}
Let $E$ be a Hilbert space. Assume that there exists
$\theta(H)\in\Theta$ such that $\|H-H_{\theta(H)}\|<\inf_{\theta: \|\theta-\theta(H)\|>\delta}\|H-H_\theta\|$, for all $\delta>0$, and suppose that $\theta\mapsto H_\theta$ is twice Fr\'echet differentiable at $\theta(H)$ with first derivative $\dot H_{\theta(H)}: \RR^d\to E$ of rank $d$. Assume that the $(d\times d)$ matrix
$V_H:=\langle \dot H_{\theta(H)},\dot H_{\theta(H)}^T\rangle-
\langle H-H_{\theta(H)},\ddot H_{\theta(H)}\rangle$ is positive definite.
If  
$\sqrt n(\HH_n-H)\cvweakly \GG_H$ in $E$, for a tight variable $\GG_H$, then
$\sqrt{n}\bigl(\HHn - H,\hat\theta_n - \theta(H)\bigr) \cvweakly (\GG_H, \hat h_H)$,   in $E\times \RR^d$, for $\hat h_H=2V_H^{-1}\langle\dot H_{\theta(H)},\GG_H\rangle$.
\end{lemma}

\begin{proof}
Define a stochastic process and map by $M_n(\theta)=\|\HHn-H_\theta\|$ and $M(\theta)=\|H-H_\theta\|$. By the triangle inequality $\sup_\theta |M_n(\theta)-M(\theta)|\le \|\HHn-H\|\to0$, in outer probability. By assumption $\theta(H)$ is a well-separated point of minimum of $M$. By a standard argument it follows that $\hat\theta_n\to \theta(H)$ in outer probability (see \cite{van2000asymptotic}, Theorem~5.7).

Redefine $M_n$ and $M$ as the square distances
$M_n(\theta)=\|\HHn-H_\theta\|^2$ and $M(\theta)=\|H-H_\theta\|^2$.
Then 
\begin{align*}
M_n(\theta)-M(\theta)&=\|\HH_n -H\|^2+2\langle \HHn-H,H-H_\theta\rangle,\\
\sqrt n(M_n-M)(\theta)-\sqrt n(M_n-M)(\theta(H))&=-2\langle \GG_n,H_\theta-H_{\theta(H)}\rangle,
\end{align*}
for $\GG_n=\sqrt n(\HHn-H)$. In view of the Cauchy-Schwarz inequality and the differentiability of $\theta \mapsto H_{\theta}$, we have
$\langle \GG_n,H_{\tilde\theta_n}-H_{\theta(H)}\rangle
=\langle \GG_n,\dot H_{\theta(H)}(\tilde\theta_n-\theta(H))\rangle+o_P\bigl(\|\tilde\theta_n-\theta(H)\|\bigr)$, for
any random sequence $\tilde\theta_n\to\theta(H)$. This verifies the stochastic expansion of Theorem~3.2.16 in \cite{van2023weak}, with $Z_n=-2\langle \GG_n,\dot H_{\theta(H)}\rangle$. 

The second condition of the latter theorem is a second-order Taylor expansion of $M$ at $\theta(H)$. 
By differentiating $M$ at its point of minimum, we find that $\langle H-H_{\theta(H)},\dot H_{\theta(H)}\rangle=0$. Combining this with the twice differentiability of $\theta\mapsto H_\theta$, we find
\begin{align*}
M(\theta)-M(\theta(H))&=    
\|H_{\theta(H)}-H_\theta\|^2
+2\langle H-H_{\theta(H)},H_{\theta(H)}-H_\theta\rangle\\
&=\|\dot H_{\theta(H)}(\theta-\theta(H))\|^2
-\langle H-H_{\theta(H)},(\theta-\theta(H))\ddot H_{\theta(H)}(\theta-\theta(H))\rangle\\
&\qquad\qquad\qquad\qquad\qquad\qquad\qquad\qquad+o\bigl(\|\theta-\theta(H)\|^2\bigr).
\end{align*}
 It follows that $M$ possesses a 
second-order Taylor expansion, with second derivative matrix
$V_H$ as given. 

Now Theorem~3.2.16 in \cite{van2023weak} gives that  $\sqrt n\bigl(\hat\theta_n-\theta(H)\bigr)=-V_H^{-1}Z_n+o_P(1)$.
The lemma follows by Slutsky's lemma.
\end{proof}

\begin{lemma}[Hilbert space]
\label{LemmaMDEquivalenceCentered}
Let $E$ be a Hilbert space.
Assume that $\{H_\theta: \theta\in\Theta\}$ satisfies the conditions of Lemma~\ref{LemmaMDHilbertian} where the map $\theta\mapsto H_\theta$ is twice continuously differentiable in a neighbourhood of $\theta(H)$, that $\sqrt n(\HH_n-H)\cvweakly \GG_H$ in $E$ and that
$\sqrt n(\HH_n^*-R_n) | X_{1.n}\cvweakly \GG_H$ in $E$,  in outer probability, for a tight random variable $\GG_H$ in $E$ and where $R_n\to H$ in outer probability.
Then, for $\GG_n^*=\sqrt n(\HH_n^*-R_n)$,
\begin{align*}
    \sqrt n(\theta_n^*-\hat\theta_n)&=\langle \dot H_{\theta(H)},\dot H_{\theta(H)}^T\rangle ^{-1}\langle \GG_n^*,\dot H_{\theta(H)}\rangle+o_P(1).
    \end{align*}
Moreover, in the case that $R_n=\HH_n$ and the matrix
$\langle \dot H_{\theta(H)},\dot H_{\theta(H)}^T\rangle -\langle \ddot H_{\theta(H)},H-H_{\theta(H)}\rangle$ is positive definite,
\begin{align*}
\sqrt n(\theta_n^{*,u}-\hat\theta_n)&
    =\bigl(\langle \dot H_{\theta(H)},\dot H_{\theta(H)}^T\rangle -\langle \ddot H_{\theta(H)},H-H_{\theta(H)}\rangle\bigr)^{-1}\langle \GG_n^*,\dot H_{\theta(H)}\rangle+o_P(1).
\end{align*}
In particular, if $H\in \{H_\theta: \theta\in\Theta\}$, then 
$\sqrt n(\theta_n^*-\theta_n^{*,u})\to 0$, in outer probability.
\end{lemma}

\begin{proof}
Lemma~\ref{LemmaMDHilbertian}  gives $\hat\theta_n\to \theta_0:=\theta(H)$ in outer probability, and then Lemma~\ref{lemma_convergence_theta_rescaled} and the continuous mapping theorem imply that $V_H\sqrt{n}(\theta_n^* - \hat\theta_n ) -\langle \GG_n^*,\dot H_{\theta_0}\rangle|X_{1.n}\cvweakly V_H\hat h_0-\langle \GG_0 ,\dot H_{\theta_0}\rangle$, in outer probability
which can be computed to be zero in case of a Hilbertian norm
since this is the stationary equation of $\hat h_0$.
Here, $V_H$ is defined as $\langle \dot H_{\theta(H)},\dot H_{\theta(H)}^T\rangle$.
This proves the first assertion.

The second assertion may be proved along the lines of the proof of Lemma~\ref{LemmaMDHilbertian}, or similarly to the argument in the preceding paragraph. Following the second path, we need to prove a version of 
Lemma~\ref{lemma_convergence_theta_rescaled} for the uncentered bootstrap values $\theta_n^{*,u}$. We start by noting that $\HHnst\to H$, and hence
$\|\HHnst - H_\theta\|\to \|H-H_\theta\|$, in outer probability. Under the conditions of Lemma~\ref{LemmaMDHilbertian}, this shows that $\theta_n^{*,u}\to \theta_0:=\theta(H)$ in outer probability. 

Next we decompose the square criterion $\|\HHnst -H_\theta\|^2$ as
\begin{align*}
\|\HHnst-\HHn+ H_{\hat\theta_n}-H_\theta\|^2
+\|\HH_n-H_{\hat\theta_n}\|^2
+2\langle \HHnst-\HHn+ H_{\hat\theta_n}-H_\theta,\HH_n-H_{\hat\theta_n}\rangle.
\end{align*}
The second term does not depend on $\theta$, nor does the term 
$\langle \HHnst-\HHn,\HH_n-H_{\hat\theta_n}\rangle$. It follows that $\theta_n^{*,u}$ minimises $\theta\mapsto \MM_n(\theta)$ given by
$$\MM_n(\theta)=\|\HHnst-\HHn+ H_{\hat\theta_n}-H_\theta\|^2
+2\langle H_{\hat\theta_n}-H_\theta,\HH_n-H_{\hat\theta_n}\rangle.$$
Because $\hat\theta_n$ is the point of minimum of $\theta\mapsto\|\HHn-H_\theta\|^2$, it satisfies the stationary equation $\langle \dot H_{\hat\theta_n}, \HH_n-H_{\hat\theta_n}\rangle=0$. Together with the twice continuous differentiability of $\theta\mapsto H_\theta$, this shows that
the inner product is equal to 
$(1/2)(\theta-\hat\theta_n)^T\langle \ddot H_{\theta_0},H-H_{\theta_0}\rangle (\theta-\hat\theta_n)+E_n\|\theta-\hat\theta_n\|^2$, for $|E_n|\to0$, in outer probability. 
The first term of $M_n(\theta)$ is bounded below by
$(1-c^{-1})\|H_{\theta}-H_{\hat\theta_n}\|^2 - (c-1)\|\HHnst-\HH_n\|^2$, for every $c>1$, where
$\|\HHnst -H_\theta\|^2
= \|\dot H_{\theta_0}(\theta-\theta_n)\|^2
+ F_n\|\theta-\hat\theta_n\|^2$,
for $|F_n|\to 0$, in outer probability. Because the matrix $V_H=\langle \dot H_{\theta_0},\dot H_{\theta_0}^T\rangle-\langle \ddot H_{\theta_0},H-H_{\theta_0}\rangle$ is positive definite, it follows that there exist constants $C,D>0$ such that $\MM_n(\theta)\ge C\|\theta-\hat\theta_n\|^2- D\|\HHnst-\HH_n\|^2$, with outer probability tending to one. From the fact that $\MM_n(\theta_n^{*,u})\le \MM_n(\hat\theta_n)=\|\HHnst-\HHn\|^2$,
we conclude that the sequence $\sqrt n(\theta_n^{*,u}-\hat\theta_n)$ is tight conditionally given $X_{1.n}$.

Now redefine $\MM_n^*$ and define a process $\MM$ by
\begin{align*}
\MM_n^*(h)&=n\|\HHnst-\HHn-(H_{\hat\theta_n+h/\sqrt n}-H_{\hat\theta_n})\|^2
-2 n\langle H_{\hat\theta_n+h/\sqrt n}-H_{\hat\theta_n},\HH_n-H_{\hat\theta_n}\rangle,\\
\MM(h)&=\|\GG_H-\dot H_{\theta_0} h\|^2-\langle h^T\ddot H_{\theta_0}h, H-H_{\theta_0}\rangle.
\end{align*}
By definition $h_n^{*,u}=\sqrt n(\theta_n^{*,u}-\hat\theta_n)$ is the point of minimum of $\MM_n^*$.
For every compact set $K\subset \RR^d$, the sequence 
of processes $(\GG_n^*,\MM_n^*)$ converges in $E\times \ell^\infty(K)$ conditionally in distribution given
$X_{1.n}$ to the process $(\GG_H,\MM)$. Therefore, by the conditional argmax theorem, it follows that $(\GG_n^*,h_n^{*,u})$ converges conditionally in distribution to $(\GG_H,\hat h_0)$, for $\hat h_0=\argmin_h \MM(h)$. By the continuous mapping theorem
$V_H \hat h_n^{*,u}-\langle \GG_n^*,\dot H_{\theta_0}\rangle$ converges conditionally in distribution to $V_H \hat h_0-\langle \GG_H,\dot H_{\theta_0}\rangle=0$.
\end{proof}

The following proposition extends Theorem~3.2.2 in \cite{van2023weak} to include joint convergence of points of maximum (also see \cite{KimPollard1990}).

\begin{proposition}[Joint argmax continuous mapping]
\label{theorem:joint_argmax_continuous_mapping}
Let $(\GG_n,\MM_n)$ be random variables in $E\times \RR^H$, for a metric space $H$ such that
$(\GG_n, \MM_n) \cvweakly (\GG, \MM)$  in $E\times\ell^\infty(K)$, for every compact $K\subset H$ and a tight random element $\GG$ in $E$ and
stochastic processes $\MM$ with upper semi-continuous sample paths that possess a unique point of maximum $\hat{h}$, 
which is tight as a map into $H$.
If $\hat h_n$ is a uniformly tight sequence of variables in $H$ such that $\MM_n(\hathn) \geq \sup_h \MM_n(h) - o_P(1)$,
then $(\GG_n, \hathn) \cvweakly (\GG, \hat{h})$.
\end{proposition}

\begin{proof}
\label{proof:theorem:joint_argmax_continuous_mapping}
Let $F\subset H$ be closed and $K\subset H$ be compact.
If $\hathn\in F\cap K$, then $\sup_{h\in F\cap K}\MM_n(h)\ge \sup_{h\in H}\MM_n(h)-\hat\epsilon_n$, by the definition of $\hathn$, for a sequence
$\hat\epsilon_n\to 0$ in probability. Therefore, for every closed set $\tilde{F}\subset E$, 
\begin{align*}
P^\outer(\GG_n \in \tilde{F}, \hathn \in F)
&\le  P^\outer\Bigl(\GG_n \in \tilde{F}, \sup_{h\in F\cap K}\MM_n(h)\ge \sup_{h\in K} \MM_n(h)-\hat\epsilon_n\Bigr)+P^\outer(\hathn\notin K).
\end{align*}
In view of the continuous mapping theorem and Slutsky's lemma, 
the sequence of variables $\bigl(\GG_n,\sup_{h\in F\cap K}\MM_n(h)-\sup_{h\in K} \MM_n(h)+\hat\epsilon_n\bigr)$
tends in distribution in $E\times\RR$ to the variable $\bigl(\GG,\sup_{h\in F\cap K}\MM(h)-\sup_{h\in K} \MM(h)\bigr)$.
Therefore, the Portmanteau lemma gives that the limsup as $n\to\infty$ of the first term on the right of the 
preceding display is bounded above by
\begin{align*}
P^\outer\Bigl(\GG \in \tilde{F}, \sup_{h\in F\cap K}\MM(h)\ge \sup_{h\in K} \MM(h)\Bigr)
&\le P(\GG \in \tilde{F}, \hat{h} \in F\cap K) + P(\hat h\notin K).
\end{align*}
The last inequality follows, because in view of upper semi-continuity, $\MM$ attains its maximum over the compact set $F\cap K$ at some point $\bar h$ of this set. If the maximum value is larger than the maximum value over $K$, as in the event in the left side of the display, and the point of global maximum
$\hat h$ is contained in $K$, then $\bar h=\hat h$, by the assumed uniqueness of the latter value and hence $\hat h\in F\cap K$.

The terms on the far right of the preceding displays can be made arbitrarily small by choice of $K$. Thus we conclude that
$\limsup_{n\to\infty}P^\outer(\GG_n \in \tilde{F}, \hathn \in F)\le P(\GG \in \tilde{F}, \hat{h} \in F)$, for every closed sets $\tilde F$ and $F$.
Since $\GG$ and $\hat{h}$ are tight, the joint law of $(\GG, \hat{h})$ is tight too and hence separable (see e.g. \cite[page 15]{van2023weak}).
Finally we apply Lemma~\ref{lem_coupled_weak_conv} to conclude that $(\GG_n, \hathn) \cvweakly (\GG,\hat{h})$.
\end{proof}

\begin{lemma}[Joint portmanteau]\label{lem_coupled_weak_conv}
 Let $(X_n,Y_n):\Omega_n \to D \times E$ be an arbitrary sequence of maps in metric spaces $D$ and $E$
such that $\limsup_{n\to\infty} P^\outer(X_n \in F_1, Y_n \in F_2) \leq L(F_1 \times F_2) $ for all closed sets $F_1\subset D$, $F_2 \subset E$ and a separable Borel measure $L$ on $D \times E$. Then $(X_n,Y_n) \cvweakly L$. 
\end{lemma}

\begin{proof}
By Corollary 1.4.5 in \cite{van2023weak}, it suffices to show that $\Expec^\outer f(X_n)g(Y_n)\to \int f\otimes g\,dL$, for every pair
of Lipschitz functions $f: D\to[0,1]$ and $g: E\to [0,1]$. The functions $f_m$ and $g_m$ defined by
$f_m(x)=\sum_{i=1}^m m^{-1}1_{F_i}(x)$, for $F_i=\{ x \in D : f(x)\geq (i-1)/m\}$, and similarly for $g_m$ with sets $F_i'$,
satisfy $0 \leq f \leq f_m \leq 1$ and $| f_m - f | \leq 1/m$ and similarly for $g_m$ and $g$. It follows that
\begin{align*}
    \Expec^\outer[f(X_n)g(Y_n)]
    \leq \Expec^\outer[f_m(X_n)g_m(Y_n)] 
    = \Expec^\outer\Bigl[\sum_{i=1}^m \sum_{j=1}^m \frac{1}{m^2} 1_{F_i\times F_j'}(X_n,Y_n)\Bigr].
\end{align*}
Since $F_i,F_j'$ are closed sets, the assumption gives that the limsup as $n\to\infty$ of the right side, for fixed $m$, is bounded above by
  \begin{align*}
\sum_{i=1}^m \sum_{j=1}^m \frac{1}{m^2} L(F_i \times F_j')
&= \int \sum_{i=1}^m \sum_{j=1}^m \frac{1}{m^2}1_{F_i}(x)1_{F_j'}(y)\,dL(x,y)
= \int f_m\otimes g_m \,dL.
    \end{align*}
 As $m \to \infty$ the right side tends to $\int f\otimes g \,dL$, by the dominated convergence theorem.
It follows that $\limsup_{n\to\infty}  \Expec^\outer[f(X_n)g(Y_n)] \leq  \int f\otimes g \,dL$.

Application of the assumption with $F_1=D$, shows that $\limsup_{n\to\infty} P^\outer(Y_n \in F_2) \leq L_2(F_2)$, for every closed $F_2\subset E$,
where $L_2$ is the marginal of $L$ on the second coordinate. By the Portmanteau theorem, it follows that $Y_n\cvweakly Y$
and hence $\Expec^\outer g(Y_n)\to \int g\,dL_2$.

The argument of the first paragraph applied to the function $(x,y)\mapsto (1-f(x))g(y)$ gives that 
$\limsup_{n\to\infty}  \Expec^\outer[g(Y_n)-f(X_n)g(Y_n)]  \le  \int (1-f)\otimes g \,dL$. Here $\Expec^\outer[g(Y_n)-f(X_n)g(Y_n)]\ge\Expec^\outer g(Y_n)-\Expec^\outer[f(X_n)g(Y_n)]$ (see Lemma~1.2.2(ii) in \cite{van2023weak}). Combined with the convergence $\Expec^\outer g(Y_n)\to \int g\,dL_2$, this gives
$\liminf_{n\to\infty}  \Expec_\inner[f(X_n)g(Y_n)] \geq  \int f\otimes g \,dL$. 

Combined the results of the first and third paragraphs give that $\Expec^\outer[f(X_n)g(Y_n)] \to  \int f\otimes g \,dL$.
\end{proof}

\section{Measurability}
\label{AppendixMeasurability}
The theory of empirical processes employs outer expectations to circumvent non-measurability relative to the Borel $\sigma$-field when considering convergence in distribution or the Delta-method (see \cite{van2023weak}). While elegant, these outer expectations depend on the definition of the variables $X_1,\ldots, X_n$ on an underlying probability space, and for general theorems it is often useful to employ a canonical definition as coordinate projections on a product space. Measurability details become extra involved when considering also bootstrap values $X_1^*,\ldots, X_n^*$. In this section we give precise definitions.

The observations $X_1,\ldots, X_n$ are assumed to be i.i.d.\ variables in a measurable space $(\Xc, \Ac)$ with distribution $H$. The canonical definition is to define these variables as the coordinate projections on the product probability space $(\Xc^n,\Ac^n, H^n)$ (if $x=(x_1,\ldots, x_n)\in\Xc^n$, then by definition $X_i(x)=x_i$). Given $X_1,\ldots, X_n$, the bootstrap values are an i.i.d.\ sample $X_1^*,\ldots, X_n^*$ from a probability measure $R_n=R_n(X_1,\ldots, X_n)$. We could formalise this by assuming that $(x_1,\ldots, x_n, B)\mapsto R_n(x_1,\ldots, x_n)(B)$ is a Markov kernel from $(\Xc^n,\Ac^n)$ into $(\Xc,\Ac)$, define a probability measure $\Pr$ on $(\Xc^n,\Ac^n)\times (\Xc^n,\Ac^n)$ by $\Pr(A\times B)=\int_A R_n(x)^n(B)\,dH^n(x)$, and define $X_1,\ldots, X_n$ and $X_1^*,\ldots, X_n^*$ as the coordinate projections on the first $n$ and last $n$ coordinates in $\Xc^n\times\Xc^n$, respectively. A disadvantage of this construction is that outer expectations of functions involving only $X_1,\ldots,X_n$ but viewed as maps on $(\Xc^n\times\Xc^n,\Ac^n\times\Ac^n,\Pr)$ are not guaranteed to be equal to the outer expectations in the original setting involving only $(\Xc^n,\Ac^n, H^n)$. 

To avoid this, we instead define the bootstrap values on a product probability space $(\Xc^n,\Ac^n,H^n)\times(\Zc,\Cc,Q)$,
where the second factor $(\Zc,\Cc,Q)$ is used to define the extra randomness involved in creating the bootstrap values. As before, the observations $X_1,\ldots, X_n$ are defined as  the coordinate projections on the first $n$ coordinates:
$X_i(x_1,\ldots, x_n, z)= x_i$, for $i=1,\ldots, n$. The bootstrap values $X_i^*$ are defined as measurable maps $X_i^*: (\Xc^n\times \Zc,\Ac^n\times \Cc)\to (\Xc,\Ac)$ such that, for every given $(x_1,\ldots, x_n)\in\Xc^n$ and $A_1,\ldots,A_n\in\Ac$,
$$Q\bigl(z\in\Zc: X_1^*(x_1,\ldots,x_n,z)\in A_1,\ldots,X_n^*(x_1,\ldots,x_n,z) \in A_n\bigr)=\prod_{i=1}^n R_n(x_1,\ldots, x_n)(A_i).$$
We give explicit constructions for $R_n$ equal to the empirical, the independence and the parametric bootstraps in Examples~\ref{ExampleDefEmpiricalBootstrap}-\ref{ExampleDefParametricBootstrap} below.

Given an arbitrary map $T: \Xc^n \times \Zc\to \RR$, the outer expectation relative to $H^n \otimes Q$ is defined as
$$\Expec^\outer T
= (H^n \otimes Q)^\outer T
:= \inf_U \int\!\!\int U(x,z)\,dH^n(x)\,dQ(z),$$ 
where the infimum is taken over all measurable maps $U:(\Xc^n \times \Zc,\Ac^n\times \Cc)\to\RR$ with $U\ge T$. The two expressions $\Expec^\outer T$ and $(H^n \otimes Q)^\outer T$ on the left are two different notations for the outer expectation, the second one longer but making the underlying measure explicit. It can be shown that in case $T$ does not depend on the $z$-coordinate, then the extra integral relative to $z$ can be removed: if $T(x_1,\ldots, x_n,z)= T_1(x_1,\ldots, x_n)$, for a map $T_1 :\Xc^n\to\RR$, then $\Expec^\outer T= \Expec^\outer T_1=(H^n)^\outer T_1$, where the outer expectation of $T_1$ is computed for $T_1$ as a map on the probability space $(\Xc^n,\Ac^n,H^n)$ relative to the measure $H^n$. (This is the \emph{perfectness} of a coordinate projection on a product space; see Lemma~1.2.5 in \cite{van2023weak}.) This is the advantage of using a product probability space. 

It can be shown that there is a minimal map $U$ as in the definition, called the \emph{least measurable majorant} of $T$ and denoted $T^\outer$, for which the infimum in the display is taken: $\Expec^\outer T=\Expec T^\outer$. This least measurable majorant is a measurable map $T^\outer: (\Xc^n\times \Zc,\Ac^n\times \Cc)\to\RR$ with $T^\outer\ge T$ and $T^\outer \le U$ $H^n \otimes Q$-almost surely, for every $U$ as in the definition. 

Because $T^\outer$ is a jointly measurable map on the product space $(\Xc^n\times \Zc,\Ac^n\times \Cc)$, a conditional expectation can be defined in the usual way as
$$\Expec_Z T^\outer:=\Expec \bigl(T^\outer(X_1,\ldots, X_n,Z) | X_1,\ldots, X_n\bigr).$$
We shall use the left side of the equation as short-hand for the right side.
This could be taken as a definition of a ``conditional outer expectation'' of the possibly non-measurable map $T$. Alternatively, we could consider the map $z\mapsto T(x_1,\ldots,x_n,z)$ for fixed $x=(x_1,\ldots, x_n)$ as a map $T^x: \Zc\to\RR$ and consider its outer expectation
$Q^\outer (T^x)$ relative to the measure $Q$, the infimum over all expectations $\int U\,dQ$ of measurable maps $U: \Zc\to\RR$ with $U(z)\ge T(x_1,\ldots, x_n,z)$, for every $z$. This might also be taken as a definition of ``conditional outer expectation'', and it would be particularly attractive if the map $z\mapsto T(x_1,\ldots,x_n,z)$ were measurable, for every fixed $x_1,\ldots,x_n$. Because $T^\outer(x,z)\ge T(x,z)= T^x(z)$, for every $z$, and $z\mapsto T^\outer(x,z)$ is measurable, for every $x$, the first proposal $\Expec_Z T^\outer$, evaluated with $X_1=x_1,\ldots, X_n=x_n$, gives a larger value than $Q^\outer (T^x)$, in general.

In analogy to outer expectations we can define inner expectations $\Expec_\inner T$ and largest measurable minorants $T_\inner$ (satisfying $\Expec_\inner T=-\Expec^\outer(-T)$ and $T_\inner= -(-T)^\outer$). For the conditional variants, we then have four quantities, which are ordered as 
\begin{align}
&\Expec \bigl(T_\inner(X_1,\ldots, X_n,Z) | X_1=x_1,\ldots, X_n=x_n\bigr)
\le Q_\inner (T^x)
\le Q^\outer (T^x)\label{EqFourQuantities}\\
&\qquad\qquad\qquad\qquad\qquad\le\Expec \bigl(T^\outer(X_1,\ldots, X_n,Z) | X_1=x_1,\ldots, X_n=x_n\bigr),\quad\text{a.s.}\nonumber
\end{align}
All of these are well defined, for fixed $x_1,\ldots, x_n$, and could be taken as definition of a conditional expectation in the case of non-measurable variables. If the maps $z\mapsto T(x_1,\ldots,x_n,z)$ are measurable,
for given $x_1,\ldots,x_n$, then the middle two expressions coincide, while all four expressions are identical if $T$ is jointly measurable.

In the theory of convergence in distribution, maps are not assumed Borel measurable, but their limits are, and this forces a converging sequence to be \emph{asymptotically measurable}. This causes that the differences between the four quantities in the preceding display are usually negligible in asymptotic arguments. Asymptotic measurability of a sequence of maps $T_n$ with values in a metric space $D$ is defined as the property 
that $\Expec^\outer h(T_n)-\Expec_\inner h(T_n)\to 0$, for every bounded continuous function $h: D\to\RR$. The sequence $T_n$ is said to converge
in distribution to a Borel measurable random element $T$, denoted $T_n\cvweakly T$, if $\Expec^\outer h(T_n)\to \Expec h(T)$, for every bounded continuous function $h: D\to \RR$, which implies that
also $\Expec_\inner h(T_n)\to \Expec h(T)$ for every such $h$ (evaluate the outer expectations for $-h$). Thus a sequence $T_n$ that converges in distribution is necessarily asymptotically measurable.

A similar observation is valid for conditional convergence in distribution of bootstrap processes. To define conditional convergence in distribution, it is convenient to metrise convergence in distribution. (A ``pointwise'' definition is possible too, see Lemma~\ref{LemmaProductFunctionsConditional} or Chapter~1.13 in \cite{van2023weak}.) It can be shown that $T_n\cvweakly T$ for a tight Borel measurable random element $T$ if and only if $\sup_{h\in \BLD}\bigl|\Expec^\outer h(T_n) - \Expec h(T)\bigr|\to 0$, where the supremum is taken over all functions $h: D\to[-1,1]$ with $|h(x)-h(y)|\le d(x,y)$, for every $x,y\in D$ (see \cite{Dudley1968,Dudley1990} or Chapter~1.12 in \cite{van2023weak}). Thus a reasonable definition for conditional convergence in distribution 
of a sequence of maps  $T_n^*: \Xc^n\times\Zc\to D$ given $X_1,\ldots, X_n$ is
\begin{equation}
\label{EqConditionalWeakConvergence}
\sup_{h\in \BLD}\bigl| \Expec_Z h(T_n^*)^\outer -\Expec h(T)\bigr|\to 0,
\end{equation}
where the convergence can be outer almost surely or in outer probability.
Here $\Expec_Z h(T_n^*)^\outer$ is the ordinary conditional expectation
$\Expec \bigl(h(T_n)^\outer | X_1,\ldots, X_n\bigr)$ of the jointly measurable variable $h(T_n^*)^\outer: (\Xc^n\times\Zc,\Ac^n\times\Cc)\to\RR$ given $X_1,\ldots, X_n$, as defined previously. An alternative would be to replace this quantity (at $X=x$) by $Q^\outer h(T_n^*)^x$, which takes the outer expectation of $h(T_n^*)^x: \Zc\to\RR$ relative to the $z$-coordinate, for fixed $x$. The latter outer expectations are smaller in general, but one then usually adds the condition that the variables $h(T_n^*)$ are asymptotically measurable relative to $H^n \otimes Q$, which blurs the difference.

For the empirical bootstrap (and general multiplier processes), the maps
$z \mapsto h(T_n^*)^x{(z)}$ are measurable and conditional weak convergence was defined in \cite{van1996}
as 
$$\sup_{h\in \BLD}\bigl| \Expec_Z h(T_n^*) -\Expec h(T)\bigr|\to 0,$$
together with, for every bounded continuous function $h: D\to\RR$,
\begin{equation}
\Expec_Z h(T_n^*)^\outer-\Expec_Z h(T_n^*)_\inner \to 0.
\label{EqAsConditionalMeasurability}
\end{equation}
The variables in both displays are functions of $X_1,\ldots, X_n$ and
both limits can be imposed outer almost surely or in outer probability. By inspection of the proofs in \cite{van1996}, it is seen that the convergence in the second display is actually uniform in $h\in \BLD$, whence $\Expec_Z h(T_n^*)$ in the first display can be replaced by the bigger quantity $\Expec_Z h(T_n^*)^\outer$, which gives \eqref{EqConditionalWeakConvergence}. The point is that under ``conditional asymptotic measurability'' all four quantities in \eqref{EqFourQuantities} are asymptotically the same. If $h(T_n^*)$ is not measurable, then this function must be replaced by a measurable function before taking an expectation, but under asymptotic measurability \eqref{EqAsConditionalMeasurability}, it does not matter which of the four quantities is chosen.

The simplest solution is to define conditional weak convergence by the single equation \eqref{EqConditionalWeakConvergence}. With this definition the continuous mapping and Delta-method theorems remain valid, as usual. (Another solution for convergence in probability, explained in \cite{bucher2019note}, is to consider joint convergence of the original and the bootstrap processes. This is technically convenient, but looses the intuition that the bootstrap is a conditional process.)

The Delta-method turns convergence in distribution of a sequence $\sqrt{n}(\HHnst - R_n)$ into convergence in distribution of the transformed processes $\sqrt{n} \bigl(\phi(\HHnst) - \phi(R_n) \bigr)$, for a given differentiable map $\phi$. For a conditional version in the notation of the present paper, we assume
\begin{equation}
        \sup_{h \in \BLD} \Bigl|
        \Expec \Bigl[h \bigl(\sqrt{n}(\HHnst - R_n) \bigr)^\outer 
        \Big| X_1, \dots, X_n \Bigr]
        - \Expec \big[ h(\GG_{R(H)}) \big] \Bigr|
        \to 0,
    \label{eq:cond_cv_Gnboot_precise_explicit}
    \end{equation}
where the convergence can be in outer probability or outer almost surely. We denote this by $\sqrt n(\HHnst-R_n)| X_1,\ldots, X_n\cvweakly \GG_{R(H)}$ in $D$. The desired conclusion is that 
$\sqrt{n} \bigl(\phi(\HHnst) - \phi(R_n) \bigr)|X_1,\ldots, X_n \cvweakly \phi_{R(H)}'( \GG_{R(H)})$ in $E$, which we understand to have the precise meaning
\begin{equation}
        \sup_{h \in \BLoneE} \Big|
        \Expec \Big[ h \Big(
        \sqrt{n} \bigl(\phi(\HHnst) - \phi(R_n) \bigr)
        \Big)^\outer \Big| X_1, \dots, X_n \Big]
        - \Expec \Big[ h \big( \phi_{R(H)}'( \GG_{R(H)} ) \big) 
        \Big] \Big|
        \to 0.
\label{eq:cond_cv_Gnboot_precise_explicit_delta_method}
\end{equation}

\begin{lemma}[Conditional Delta-method]
\label{lemma:Functional_Delta_Method_boot}
Let  $D$ and $E$ be normed spaces. Let $\phi: D_\phi\subset D\to E$ be Hadamard differentiable at $R(H)$ tangentially to a subspace $D_0\subset D$. Let $R_n$ and $\HHnst$ be maps with values in $D_\phi$ such that
the sequence  $\sqrt n(\HHnst-R_n)| X_1,\ldots, X_n\cvweakly \GG_{R(H)}$ in outer probability, for a tight random element $\GG_{R(H)}$ taking its value in $D_0$. 
\begin{enumerate}
\item If the sequence $\sqrt n\bigl(R_n-R(H)\bigr)$ is asymptotically tight in $D$ with limit points in $D_0$, then 
$\sqrt{n} \bigl(\phi(\HHnst) - \phi(R_n) \bigr)|X_1,\ldots, X_n \cvweakly \phi_{R(H)}'( \GG_{R(H)})$ in $E$, in outer probability.
\item If $\phi: D_\phi\subset D\to E$ is uniformly Hadamard differentiable at $R(H)$, then the condition that the sequence $\sqrt n\bigl(R_n-R(H)\bigr)$ is asymptotically tight in $D$ can be relaxed to convergence $R_n\to R(H)$ in $D$ in outer probability. 
\item If $\phi: D_\phi\subset D\to E$ is uniformly Hadamard differentiable at $R(H)$, and $R_n\to R(H)$ in $D$ outer almost surely and $\sqrt n(\HHnst-R_n)| X_1,\ldots, X_n\cvweakly \GG_{R(H)}$, outer almost surely, then 
$\sqrt{n} \bigl(\phi(\HHnst) - \phi(R_n) \bigr)|X_1,\ldots, X_n \cvweakly \phi_{R(H)}'( \GG_{R(H)})$ in $D$, outer almost surely.
\end{enumerate}
\end{lemma}

\begin{proof}
This essentially combines Theorems~3.10.11 and~3.10.13 in \cite{van2023weak}. Because the latter theorems were written for more special bootstrap processes, we provide a complete proof.

Without loss of generality we can assume that the derivative $\phi_{R(H)}': D\to E$ is defined and continuous on the whole space $D$. (There always is an extension of the derivative, possibly taking values in a bigger space. See \cite{van2023weak}, page~525 for details.) Because  for every  given bounded Lipschitz function $h: E\to\RR$, the
function $h\circ\phi_{R(H)}': D\to \RR$ is bounded Lipschitz of norm
at most $\|\phi_{R(H)}'\|$ bigger,  \eqref{eq:cond_cv_Gnboot_precise_explicit} implies
$$\sup_{h\in \BLoneE}
\Bigl|\Expec_Z h\Bigl(\phi_{R(H)}'\bigl(\sqrt n(\HHnst-R_n)\bigr)\Bigr)^\outer
-\Expec h\bigl(\phi_{R(H)}'(\GG_{R(H)})\bigr)\Bigr|\to 0,$$
in outer probability or outer almost surely, corresponding
to which of the two assumptions is made. Next
\begin{align}
&\sup_{h\in \BLoneE)}
\Bigl|\Expec_Z h\Bigl(\sqrt n\bigl(\phi(\HHnst)-\phi(R_n)\bigr)\Bigr)^\outer
-\Expec_Z h\Bigl(\phi_{R(H)}'\bigl(\sqrt n(\HHnst -R_n)\bigr)\Bigr)^\outer\Bigr|
\label{eq:decomposition_Delta_method_epsilon}
\\
&\qquad\qquad
\le \epsilon +2\Prob_Z\Bigl(\bigl\|\sqrt n\bigl(\phi(\HHnst)-\phi(R_n)\bigr)
-\phi_{R(H)}'\bigl(\sqrt n(\HHnst-R_n)\bigr)\bigr\|^\outer>\epsilon\Bigr),
\nonumber
\end{align}
for every $\epsilon > 0$.
The lemma is proved once it has been shown that the conditional probability on the right converges to zero in outer probability or outer almost surely.

In view of \eqref{eq:cond_cv_Gnboot_precise_explicit},
$\Expec h \bigl(\sqrt{n}(\HHnst - R_n) \bigr)^\outer \to 
\Expec h(\GG_{R(H)})$, for every bounded Lipschitz function $h$. This implies that $\sqrt{n}(\HHnst - R_n)\cvweakly \GG_{R(H)}$, also unconditionally.
In particular it is asymptotically tight, and hence so is joint sequence
$\sqrt n\bigl(\HHnst - R_n, R_n-R(H)\bigr)$ and also their sum
$\sqrt n\bigl(\HHnst-R(H)\bigr)$. The limit points concentrate on $D_0$, by assumption. By the unconditional Delta-method (for instance Theorem~3.10.4 in \cite{van2023weak}),
\begin{align*}
\sqrt n\bigl(\phi(\HHnst)-\phi(R(H))\bigr)
&=\phi_{R(H)}'\bigl(\sqrt n(\HHnst-R(H))\bigr)+o_P^*(1),\\
\sqrt n\bigl(\phi(R_n)-\phi(R(H))\bigr)
&=\phi_{R(H)}'\bigl(\sqrt n(R_n-R(H))\bigr)+o_P^*(1),
\end{align*}
where $o_P^*(1)$ designates a term that converges to $0$ in outer probability.
Subtract these equations to conclude that 
the sequence $\sqrt n\bigl(\phi(\HHnst)-\phi(R_n)\bigr)
-\phi_{R(H)}'\bigl(\sqrt n(\HHnst-R_n)\bigr)$ converges 
(unconditionally) to zero in outer probability.
Thus, the conditional probability on the right in the last display of the preceding paragraph converges to zero in outer mean. This concludes the
proof of the first assertion of the lemma.

In the remainder of the proof we assume that $\phi$ is uniformly Hadamard differentiable. Fix $\epsilon>0$ and choose a compact set 
$K\subset D_0$ such that $\Pr(\GG_{R(H)}\notin K)<\epsilon$. By the uniform
Hadamard differentiability of $\phi$, there exist $\delta,\eta >0$
such that for every $H'\in D_\phi$, $\|H'-R(H)\|<\eta$, $t<\eta$,
$H'+th\in D_\phi$, and $d(h, K)<\delta$:
$$\left\|\frac{\phi(H'+th)-\phi(H')}t-\phi_{R(H)}'(h)\right\|<\epsilon.$$
Applying this with
$H'= R_n$,
$h = \sqrt{n}(\HHnst - R_n)$,
$t = 1 / \sqrt{n}$
gives that whenever
$h = \sqrt{n}(\HHnst - R_n)
\in K^\delta$ (meaning that $d(h, K)<\delta$), $t < \eta$
and $\|R_n - R(H)\| < \eta$,
we have
\begin{align*}
    \left\|\sqrt{n} (\phi(\HHnst)-\phi(R_n))
    - \phi_{R(H)}'(\sqrt{n}(\HHnst - R_n))
    \right\|
    < \epsilon.
\end{align*}
\noindent
Consequently and applying a union bound,
for $n \ge 1/\eta^2$,
the right side of
\eqref{eq:decomposition_Delta_method_epsilon}
is bounded by
$$\epsilon+2\Expec_Z  1_{D-K^\delta}\bigl(\sqrt n(\HHnst-R_n)\bigr)^\outer
+1\bigl\{\|R_n-R(H)\|\ge \eta\bigr\}^\outer.$$
The last term converges to zero in outer probability or outer almost surely by assumption. The function $h(z)= \delta^{-1}\bigl(d(z,K)\wedge \delta\bigr)$ 
is bounded and Lipschitz and satisfies
$1_{D-K^\delta}\le h\le 1_{D-K}$. Hence the conditional
expectation in the middle term is bounded by
$$\Expec_Z h\bigl(\sqrt n(\HHnst-R_n)\bigr)^\outer\to \Expec h(\GG_{R(H)})
\le \Pr(\GG_{R(H)}\notin K) <\epsilon.$$
Thus, the conditional probability in the right side of the display 
converges to zero
almost surely. This concludes the proof of  the second and third assertions of the lemma.
\end{proof}

\begin{corollary}
\label{lemma:Continuous_mapping_theorem_boot}
Under the conditions of Lemma~\ref{lemma:Functional_Delta_Method_boot}, or more generally, if Equation~\eqref{eq:cond_cv_Gnboot_precise_explicit_delta_method}
    holds, then
    \begin{equation*}
        \sup_{h \in \BLoneR} \bigg|
        \Expec \Big[ h \Big( \big\|
        \sqrt{n} \big(\phi(\HHnst) - \phi(R_n) \big)
        \big\|_E \Big)^\outer
        \Big| X_1, \dots, X_n \Big]
        - \Expec \Big[ h \big(
        \big\| \phi_{R(H)}'( \GG_{R(H)} )  \big\|_E \big) \Big] 
        \bigg| \to 0,
    \end{equation*}
outer almost surely or in outer probability. 
\end{corollary}

\begin{proof}
For $h\in \BLoneR$, the map $x\mapsto h(\|x\|_E)$ is contained in $\BLoneE$.
\end{proof}

\begin{lemma}
\label{LemmaProductFunctionsConditional}
Let $(X_n,Y_n)$ be arbitrary maps with values in the product $D\times E$ of metric spaces $D$ and $E$ such that $\Expec_Z^\outer f(X_n)g(Y_n)\to \Expec f(X)g(Y)$,  almost surely or in probability, for a tight Borel measurable random element $(X, Y)$ in $D\times E$ and every bounded Lipschitz maps $f: D\to\RR$ and $g: E\to\RR$. Then
$$\sup_{h\in \BLDE}|\Expec_Z h(X_n,Y_n)^\outer-\Expec h(X,Y)\bigr|\to0,\qquad\text{outer a.s. or in outer prob.}$$
\end{lemma}

\begin{proof}
Let $\Hc$ be the set of all finite linear combinations of functions $h: D\times E\to \RR$ of the form $h(x,y)=f(x)g(y)$, for $f: D\to \RR$ and $g: E\to\RR$ bounded Lipschitz functions. We first show that $\Expec_Z^\outer h(X_n,Y_n)\to \Expec h(X,Y)$, almost surely or in probability, for every $h\in\Hc$. 

Take any finite set of pairs of bounded Lipschitz functions $(f_j,g_j)$. Because $\bigl(\sum_jf_j\times g_j\bigr)^\outer\le \sum_j(f_j\times g_j)^\outer$, it follows that
$\Expec_Z^\outer (\sum_jf_j(X_n) g_j(Y_n))\le\sum_j\Expec_Z(f_j(X_n) g_j(Y_n))^\outer$, which by assumption tends to $\sum_j\Expec f_j(X)g_j(Y)$, almost surely or in probability.
By applying the assumed convergence of the outer expectations to the function $-f\times g=(-f)\times g$, obtain that
also $\Expec_{Z,\inner}f(X_n)g(Y_n)\to \Expec f(X)g(Y)$, and conclude in the same way that $\Expec_{Z,\inner} (\sum_jf_j(X_n) g_j(Y_n))\ge 
\sum_j\Expec_{Z} (f_j(X_n) g_j(Y_n))_\inner$, which tends to $\sum_j\Expec f_j(X)g_j(Y)$,
almost surely or in probability. Because the limit is the same, it follows
that $\Expec_Z^\outer (\sum_jf_j(X_n) g_j(Y_n))\to \sum_j\Expec f_j(X)g_j(Y)$, almost surely or in probability. 

By the assumed tightness of the limit $(X,Y)$,
for every $\epsilon>0$ there exists a compact set $K\subset D\times E$ such that
$\Prob\bigl((X,Y)\in K\bigr)>1-\epsilon$. Because $\Hc$ is a vector space that is closed under taking products, contains the constant function and separates points of $D\times E$, its set of restrictions to $K$ are dense in $C(K)$, by the Stone-Weierstrass theorem
(see \cite{Jameson}, page~266). Hence for $h\in \BLDE$ and any $\epsilon>0$,
there exists $h_\epsilon\in \Hc$ with $\sup_{(x,y)\in K}|h-h_\epsilon|(x,y)<\epsilon$. 
In view of the compactness of $K$, there exists $\delta>0$ such that $\sup_{(x,y)\in K^\delta}|h-h_\epsilon|(x,y)<\epsilon$, where $K^\delta$ is the set of points $(x,y)\in D\times E$ within distance $\delta$ of $K$.
For a Lipschitz function $\chi$ with  $1_K\le \chi\le 1_{K^\delta}$, we then have
$|\Expec_Z^\outer (h-h_\epsilon)\chi(X_n,Y_n)|\le \epsilon$ and $|\Expec (h-h_\epsilon)\chi(X,Y)|\le \epsilon$, while $|\Expec_Z^\outer (h-h_\epsilon)(1-\chi)(X_n,Y_n)|\le \Expec_Z^\outer (1-\chi)(X_n,Y_n)$, which converges  almost surely or in probability to $\Expec (1-\chi)(X,Y)|\le \epsilon$. It follows that $|\Expec_Z^\outer h(X_n,Y_n)-\Expec_Z^\outer h_\epsilon(X_n,Y_n)|\le \epsilon+E_n$, for random variables $E_n$ that converge almost surely or in probability to a constant smaller than $\epsilon$. Since
$\Expec_Z^\outer h_\epsilon(X_n,Y_n)\to \Expec h_\epsilon(X,Y)$
and the argument is true for every $\epsilon>0$, it follows that
$\Expec_Z^\outer h(X_n,Y_n)\to \Expec h(X_n,Y_n)$ almost surely or in probability.

Finally we show that this convergence is uniform in $h\in\BLDE$. The restrictions of the
functions $h\in\BLDE$ to the compact set $K$ are equicontinuous and uniformly bounded. Therefore the set of restrictions is totally bounded for the uniform norm, by the Ascoli-Arzela theorem.
For every $\eta>0$ there exist finitely many functions $h_1,\ldots, h_m$ so that
$\sup_{(x,y)\in K}\min_i|h-h_i|(x,y)<\eta$. There exists $\delta>0$ so this remains true if the supremum is taken over $K^\delta$ (see Problem~1.3.13 in \cite{van2023weak}). Then arguing as before, we see that 
$$\sup_{h\in\BLDE}\bigl|\Expec_Z^\outer h(X_n,Y_n)-\Expec h(X,Y)\bigr|\le 4\epsilon+\max_{i}\bigl|\Expec_Z^\outer h_i(X_n,Y_n)-\Expec h_i(X,Y)\bigr|.$$
The right side tends to $4\epsilon$, almost surely or in probability. This being true for every $\epsilon>0$, shows that the left side tends to zero.
\end{proof}

\begin{example}[Empirical bootstrap]
\label{ExampleDefEmpiricalBootstrap} 
The empirical bootstrap measure can be represented as $\HHnst=\sum_{i=1}^nM_{n,i}\delta_{X_i}$, where $(M_{n,1},\ldots, M_{n,n})$ is a multinomial vector with parameters $n$ and $(1/n,\ldots, 1/n)$, independent of $X_1,\ldots, X_n$ (the numbers of times that the $X_i$ are resampled). The vector $M_n$ can be defined on the factor $(\Zc,\Cc,Q)$, while $X_1,\ldots, X_n$ are the coordinate projections on the
first $n$ coordinates of $\Xc^n\times\Zc$. By constructing $M_n$ as the sum of $n$ independent multinomial variables $m_n^{(1)},\ldots, m_n^{(n)}$ with parameters $1$ and $(1/n,\ldots, 1/n)$, the individual bootstrap values can also be defined, as $X_i^*=X_j$ if $m_n^{(i)}=e_j$, for $e_j$ the $j$th unit vector in $\RR^n$. More formally, for $\omega=(x_1,\ldots,x_n,z)\in\Xc^n\times \Zc$, we define $X_i(\omega)=x_i$ and
$X_i^*(\omega)=x_j$ if $m_n^{(i)}(z)=e_j$.
\end{example}

\begin{example}[Independence bootstrap]
\label{ExampleDefIndependenceBootstrap} 
The independence bootstrap is the product of two empirical bootstraps and hence can be formally defined by duplicating the construction in Example~\ref{ExampleDefEmpiricalBootstrap}.
\end{example}

\begin{example}[Parametric bootstrap]
\label{ExampleDefParametricBootstrap} 
If $(\Xc,\Ac)$ is a Polish space with its Borel $\sigma$-field, and $(\theta,A)\mapsto H_\theta(A)$ is a Markov kernel, then there exists a jointly measurable map $\psi: \Theta\times[0,1]\to\Xc$ such that $\psi(\theta,U)$ is distributed according to $H_\theta$ if $U$ is a standard uniform variable. For $\Xc=\RR$, the quantile transformation $\psi(\theta,U)=F_\theta^{-1}(U)$, for $F_\theta$ the cumulative distribution function of $H_\theta$, gives a concrete definition, and for $\Xc=\RR^d$ one can employ conditional quantile functions. For a general Polish space the claim follows, because such a space is measurable-isomorph to $\RR$ or to a countable subset of $\RR$.

Given a map $\psi$, we can choose $(\Zc,\Cc,Q)$ any probability space that is rich enough to carry $n$ i.i.d.\ uniform variables $U_1,\ldots, U_n$, and given an estimator $\hat\theta_n$ define parametric bootstrap values by $X_i^*=\psi(\hat\theta_n,U_i)$. More formally,  for $\omega=(x_1,\ldots, x_n,z)\in \Xc^n\times\Zc$, the original observations are defined as $X_i(\omega)=x_i$ and the bootstrap values as $X_i^*(\omega)=\psi\bigl(\hat\theta_n(x_1,\ldots, x_n),U_i(z)\bigr)$.
\end{example}

\section{Supplementary simulations}

\subsection{Fixed design regression}
\label{app:Fixed design regression}

\begin{figure}[tbp]
  \centering
  \resizebox{0.9\textwidth}{!}{\input{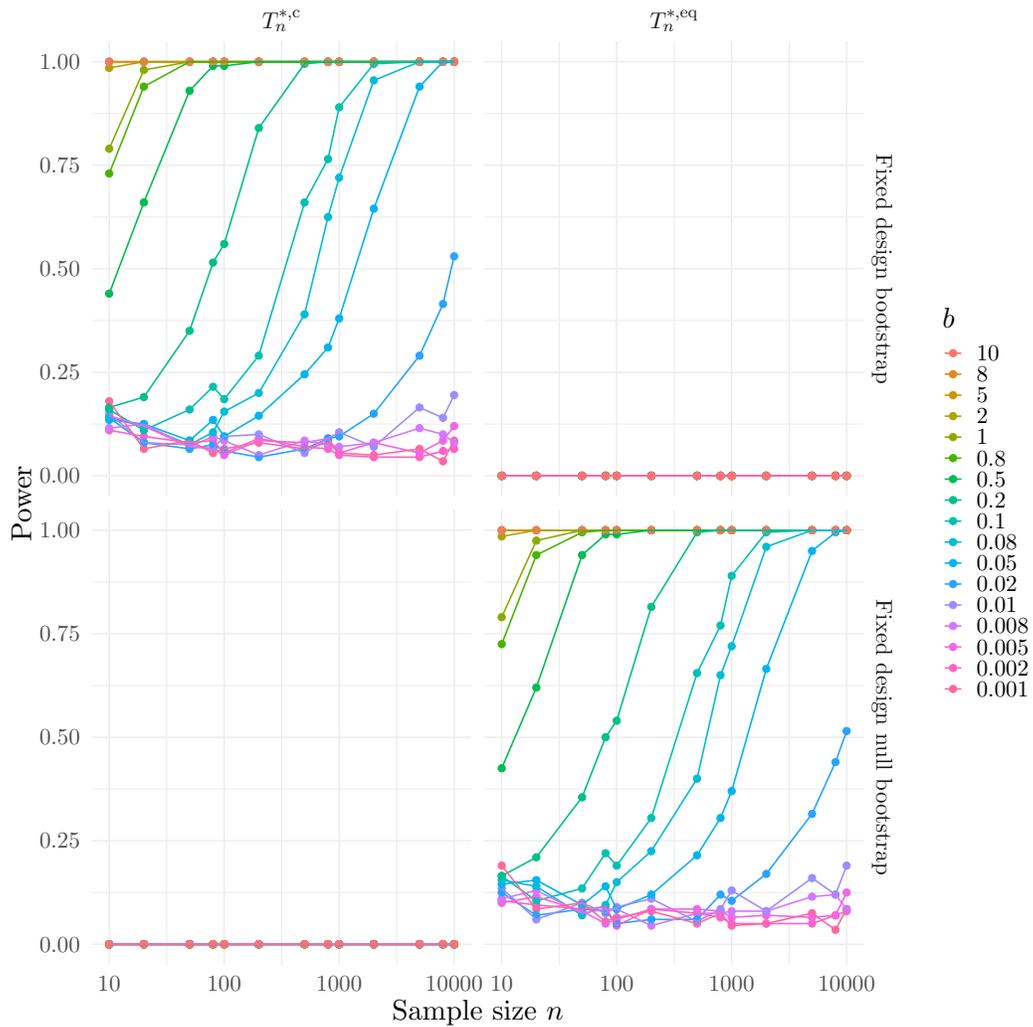}}
  \caption{
    Power in the regression setting as a function of the sample size, for different values of the coefficient $b$ and different combinations of fixed design bootstrap resampling schemes and bootstrap test statistics ($\Tncent$ or $\Tneq$).
  }
  \label{fig:power_test_fixeddesignregressiontest}
\end{figure}

\begin{figure}[tbp]
  \centering
  \resizebox{0.7\textwidth}{!}{\input{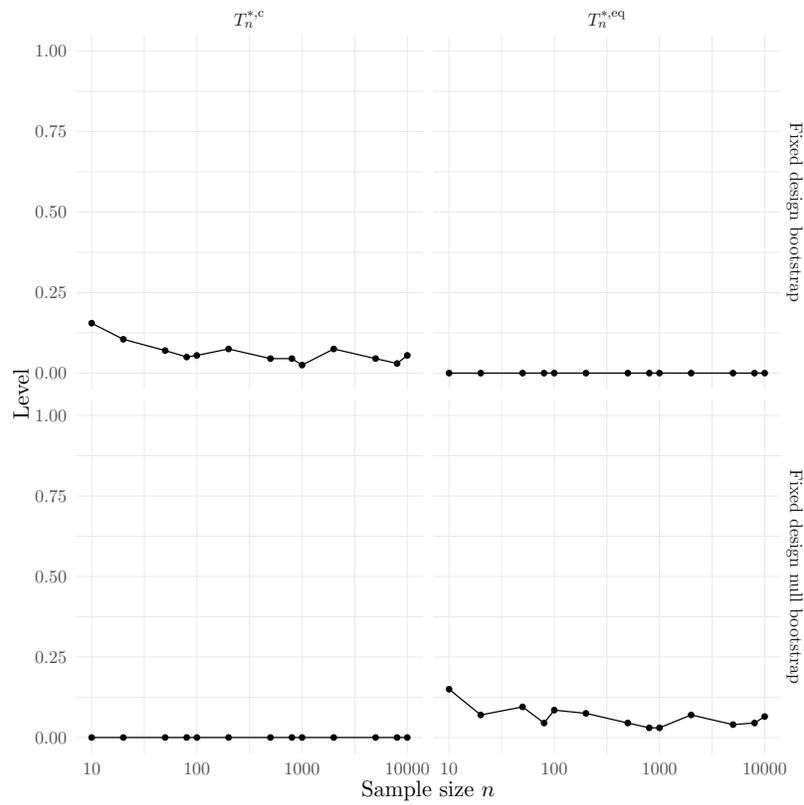}}
  \caption{
    Level in the regression setting as a function of the sample size, for different combinations of fixed design bootstrap resampling schemes and bootstrap test statistics ($\Tncent$ or $\Tneq$).
  }
  \label{fig:level_test_fixeddesignregressiontest}
\end{figure}

Here we include the power and level simulation for the fixed design residual bootstrap scheme in Figures \ref{fig:power_test_fixeddesignregressiontest} and \ref{fig:level_test_fixeddesignregressiontest}. The simulation procedure is similar to the regression setting, described in Section \ref{SectionSimulation}. Here we use resampling scheme (iv) from Section \ref{sec:residual_bootstraps}. Furthermore, we also include a fixed design null bootstrap, which is equal to resampling scheme (iv), except that $\hat{b}_n$ is set to zero. This is also the reason we call it a `null' bootstrap. From the theory in Section \ref{sec:residual_bootstraps} we know that the fixed design residual bootstrap, paired with $\Tncent$ will work. The simulation results also show this. The fixed design null bootstrap resampling procedure should be paired with $\Tneq$ to work.

\begin{figure}[tbhp]
  \centering
  \resizebox{0.92\textwidth}{!}{\input{GOF_power_MDestimator_plot}}
  \caption{
    Power of the goodness-of-fit test as a function of the sample size for different data generating processes (including the Cauchy distribution) and different combinations of bootstrap resampling schemes and bootstrap test statistics ($\Tncent$ or $\Tneq$).
    Here, the bootstrap version estimator is the `non-centered' $\theta_n^{*,MD}=\argmin_\theta \| \HHnst - H_{\theta} \|$.
    }
  \label{fig:power_test_MDestimator_GOFtest_including_Cauchy}

  \vspace{0.8cm}

  \centering
  \resizebox{0.92\textwidth}{!}{\input{GOF_power_MDcentestimator_plot}}
  \caption{
    Power of the goodness-of-fit test as a function of the sample size for different data generating processes (including the Cauchy distribution) and different combinations of bootstrap resampling schemes and bootstrap test statistics ($\Tncent$ or $\Tneq$). 
    Here, the bootstrap version estimator is the `centered'  $\theta_n^{*,MD}=\argmin_\theta \| \HHnst - H_{\theta}- R_n+H_{\hat\theta_n} \|.$
  }
  \label{fig:power_test_MDcentestimator_GOFtest_including_cauchy}
\end{figure}

\subsection{Goodness-of-fit}\label{app:GOF}

We have performed an additional goodness-of-fit power simulation that includes the Cauchy distribution. The setup is the same as in in Section \ref{sec:simulation_GOF_general}. The results are shown in Figures \ref{fig:power_test_MDestimator_GOFtest_including_Cauchy} and \ref{fig:power_test_MDcentestimator_GOFtest_including_cauchy}. The precise distribution parameters used in the goodness-of-fit simulations are given in Table \ref{tab:distributions_GOF_simulation}.

\medskip

Note that the power for Cauchy alternatives behaves in a different way compared to the other distributions for the empirical bootstrap. This is not surprising since the null hypothesis is that the distribution is normal with an unknown mean $\theta$. In general, one estimates $\theta = \theta(H) = \Expec_H[X]$ by $\thetahat = \bar X_n$.
When $H$ is a Cauchy distribution $\theta(H)$ is not defined, and therefore the theoretical results cannot be applied.
It is interesting to note that the power still seems to tend to $1$ for the parametric bootstrap (null bootstrap) for the Cauchy case, even though the corresponding Lemma~\ref{LemmaGoodnessNull} cannot be applied either.

\begin{table}[ht]
\centering
\begin{tabular}{|l|l|}
\hline
\textbf{Distribution} & \textbf{Parameters} \\
\hline
Normal & $\Nc(0,1)$ \\
t-Distribution & $t(\nu = \{60, 40, 20, 10, 5\})$ \\
Log-normal & $\mathrm{LogNormal}(\mu=0, \sigma=\{0.3, 0.6, 0.8, 1.0\})$ \\
Mixtures & $\Nc(\mu=\{-1.5, -2, -3, -4\}, \sigma=1)$ and $\Nc(\mu=\{1.5, 2, 3, 4\}, \sigma=1)$ \\
Cauchy & $\mathrm{Cauchy}(\mu=0, \gamma=\{1, 2, 3, 4\})$ \\
\hline
\end{tabular}
\caption{Distributions used for goodness-of-fit simulations with their corresponding parameters.}
\label{tab:distributions_GOF_simulation}
\end{table}

\subsection{Which bootstrap is better?}\label{app:comparison}

\begin{figure}[htbp]
  \centering
  \resizebox{\textwidth}{!}{\input{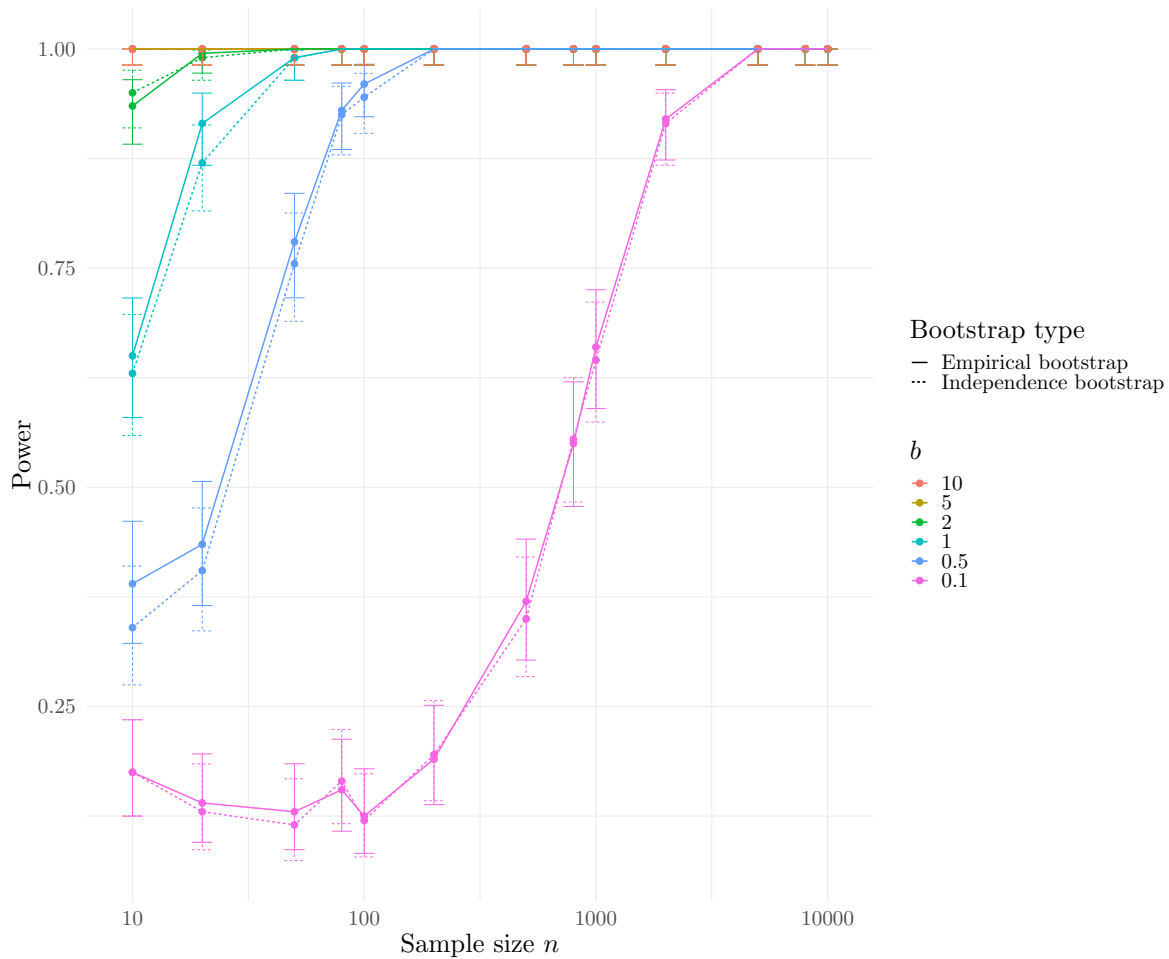}}
  \caption{
    Power of the independence test as a function of the sample size. The empirical bootstrap is used with $\Tncent$ and the independence bootstrap is used with $\Tneq$. The value of $b$ shows the value used in the data-generating process $Y_i = b X_i + \epsilon_i$. The error bars indicate the 95$\%$ confidence interval according to the binomial test for the simulation error.
  }
  \label{fig:independence_comparison_plot}
\end{figure}

Figures \ref{fig:independence_comparison_plot} and \ref{fig:regression_comparison_plot}
displays comparisons of the power of the two theoretically valid bootstrap schemes, respectively in the independence and in the regression settings.
In both settings, the theoretically valid combinations are: the empirical bootstrap paired with the centred test statistic;
and the independence bootstrap paired with the equivalent test statistic.

\medskip

These powers are estimated using the Monte Carlo technique, i.e. as the percentage of rejections observed out of $N = 200$ simulations.
In both figures, we display the $95\%$ confidence intervals around each (estimated) power using error bars, computed using the \texttt{R} function \texttt{binom.test()}.
Indeed, the number of rejections follows a binomial distribution with parameters $N$ and $p$ being the (true) power of the test.

\medskip

In Figure~\ref{fig:regression_comparison_plot}, the empirical bootstrap has a significantly higher power than the independent bootstrap for small sample size $n = 10$ and intermediate values of $b$ ($0.5$, $1$, $2$). For larger sample sizes, the differences are smaller and therefore not significant compared to the simulation error.
For $n = 10$ and very large or small values of $b$, the power of both bootstrap schemes are respectively close to $100\%$ and $5\%$, and therefore cannot be distinguished.

\medskip

In Figure~\ref{fig:independence_comparison_plot}, none of the differences seems to be statistically significant. Nevertheless, to confirm the intuition given in the regression setting, we choose the most promising combinations of parameters ($n = 20, b = 1$) which showed the largest difference in power between both methods.
In that case, we performed $N = 2 \, 000$ simulations, again, to obtain tighter confidence intervals. 
Performing a proportion test using the \texttt{R} function \texttt{prop.test()} gives a significant difference between the estimated powers $0.845$ and $0.797$ for the empirical and independence bootstrap, respectively, with a corresponding p-value of $8.918\cdot 10^{-5}$.

\medskip

This confirms the hypothesis that the empirical bootstrap outperforms the independence bootstrap, although the improvement in power is quite small.
We conjecture that this also holds for the other combinations of $n$ and $b$, in both settings, with even smaller differences between both schemes.

\begin{figure}[htbp]
  \centering
  \resizebox{1\textwidth}{!}{\input{Regression_comparison_plot}}
  \caption{
    Power of the regression test as a function of the sample size. 
    The empirical bootstrap is used with $\Tncent$ and the independence bootstrap is used with $\Tneq$. The value of $b$ shows the value used in the data-generating process $Y_i = b X_i + \epsilon_i$. The error bars indicate the 95$\%$ confidence interval according to the binomial test for the simulation error.
  }
  \label{fig:regression_comparison_plot}
\end{figure}

\end{document}